\documentclass{article}

\usepackage[english]{babel}

\makeatletter
\let\@fnsymbol\@arabic
\makeatother

\usepackage{geometry}
\usepackage{amsmath}
\usepackage{graphicx}

\usepackage{verbatim}
\usepackage{subcaption}
\usepackage{amsgen,amsmath,amstext,amsbsy,tikz,amssymb}
\usepackage{fancyhdr}
\usepackage[colorlinks=true, urlcolor=blue,  linkcolor=black, citecolor=black]{hyperref}

\usepackage{rotating}

\usepackage{booktabs}
\usepackage{wasysym}
\usepackage{algorithm}
\usepackage{algpseudocode}
\usepackage{multirow}
\usepackage{nicefrac}
\usepackage{hhline}

\usepackage[sort]{cite}

\DeclareMathOperator{\af}{\mathcal{A}\MRkern \mathcal{F}}
\newcommand{\MRkern}{%
  \mkern-3.65mu
  \mathchoice{}{}{\mkern0.2mu}{\mkern0.5mu}%
}

\newcommand{\AF}[1]{\mathcal{A}\MRkern \mathcal{F}_{\hspace{-1pt} {#1}}}

\graphicspath{{figures/}}

\usepackage{lipsum}
\usepackage{amsfonts}
\usepackage{epstopdf}
\ifpdf
  \DeclareGraphicsExtensions{.eps,.pdf,.png,.jpg}
\else
  \DeclareGraphicsExtensions{.eps}
\fi

\newenvironment{@abssec}[1]{%
       \vspace{.05in}\footnotesize
       \parindent .2in
         {\upshape\bfseries #1. }\ignorespaces 
     }
     {\par\vspace{.1in}}
\newcommand\keywordsname{Key words}
\newcommand\MSCcodesname{MSC codes}
\newenvironment{keywords}{\begin{@abssec}{\keywordsname}}{\end{@abssec}}
\newenvironment{MSCcodes}{\begin{@abssec}{\MSCcodesname}}{\end{@abssec}}

\title{Augmentations of Forman's Ricci Curvature and their Applications in Community Detection}

\author{Lukas Fesser*$^,$\thanks{Harvard University, Cambridge, MA.}
\and Sergio Serrano de Haro Iváñez*$^,$\thanks{University of Oxford, Oxford, UK.}
\and Karel Devriendt\thanks{Max Planck Institute for Mathematics in the Sciences, Leipzig, Germany.}
\and Melanie Weber\footnotemark[1]
\and Renaud Lambiotte\footnotemark[2]}

\date{}

\usepackage{amsopn}

\begin{document}

\renewcommand{\thefootnote}{$\ast$} 
\footnotetext{Joint first authors.}
\renewcommand{\thefootnote}{\arabic{footnote}}

\maketitle

\begin{abstract}
The notion of curvature on graphs has recently gained traction in the networks community, with the Ollivier-Ricci curvature (ORC) in particular being used for several tasks in network analysis, such as community detection. In this work, we choose a different approach and study augmentations of the discretization of the Ricci curvature proposed by Forman (AFRC). We empirically and theoretically investigate its relation to the ORC and the un-augmented Forman-Ricci curvature. In particular, we provide evidence that the AFRC frequently gives sufficient insight into the structure of a network to be used for community detection, and therefore provides a computationally cheaper alternative to previous ORC-based methods. Our novel AFRC-based community detection algorithm is competitive with an ORC-based approach.
\end{abstract}

\begin{keywords}
Discrete Curvature, Network Analysis, Community Detection
\end{keywords}

\begin{MSCcodes}
        05C82, 05C10, 53C21, 68R10, 05C75
\end{MSCcodes}

\section{Introduction} 

Curvature is a key notion in differential geometry, and Ricci curvature is one of several ways of formalizing it. Intuitively, Ricci curvature measures the deformation of a small ball as we push it along a geodesic (the shortest path between two points), and quantifies how much a shape locally differs from standard, flat Euclidean space. As Ricci curvature encapsulates essential information in differential geometry and related fields (e.g. it is a term of Einstein’s field equations in general relativity), there have been multiple attempts to profitably extend its definition to other domains.\\

\noindent The structural insight gained via discrete Ricci curvature has proven valuable for applications in Network Science and Machine Learning. 
Discrete Ricci curvature relates to many classical network characteristics and network analysis tasks~\cite{weber2017characterizing} and, as a result, has been applied for the analysis of network data in many domains, including biological networks~\cite{tannenbaum2015ricci,weber2017curvature},  chemical networks~\cite{leal2021forman} and financial transaction networks~\cite{sandhu2016ricci}, among others. In Machine Learning, discrete curvature has been utilized to characterize and mitigate the oversquashing effect in Graph Neural Networks~\cite{topping2021understanding} and to identify embedding spaces in Representation Learning~\cite{weber2020neighborhood}.\\

\noindent The two curvature versions studied in this paper are respectively due to Ollivier ~\cite{ollivier2007Riccimetricspaces,ollivier2009RicciMarkovchains} and Forman~\cite{forman2003bochner}. Ollivier-Ricci curvature (ORC) is a direct discretization of the original definition, substituting spheres and geodesics by, respectively, node neighborhoods and shortest edge-paths. The construction of Forman-Ricci curvature (FRC) takes a different approach: it follows from Bochner’s method ~\cite{bochner1,bochner2,bochner3}, which in the differential case is used to write the Laplacian as the sum of a Laplacian-like term and a term that depends on the Ricci curvature. Forman’s work consists in discretizing both the Laplacian and Bochner’s method, which produces a term that, by analogy, can be interpreted as a discrete Ricci curvature. Thus, ORC and FRC are two measures with solid geometric foundations yet presenting significant differences, which makes them two potential ways of gaining a better understanding of the geometry of networks.\\

\noindent Here, we focus on geometric approaches to community detection, which have previously received much attention \cite{kovacs2021inherent, xue2017reliable, yang2021hidden}, particularly the use of ORC \cite{gosztolai2021unfolding,ni2019community,sia2019ollivier,sia2022inferring,tian2023mixed}.  Community structure, i.e. the mesoscale organization of a network into groups of nodes with similar connections, is widely observed in real-world networks and its detection is a classical problem in the study of networks and graphs \cite{abbe2018community, porter_2009_communities}. Community detection using ORC relies on the observation that positively curved edges tend to be well connected, and hence naturally lie inside communities. Negatively curved edges, on the other hand, can be interpreted as bridges between communities. ORC-based community detection algorithms have been found to perform well on both synthetic and real-world networks. These methods are, however, computationally expensive as the definition of ORC involves the computation of network flows of shortest paths via the Wasserstein distance, which has cubic complexity and is often approximated via the (quadratic) Sinkhorn algorithm~\cite{cuturi2013sinkhorn}.\\

\noindent 
One advantage of FRC over ORC is its significantly smaller computational complexity: in its simplest form, the FRC of an edge only depends on the degrees of its end points. This is often too simple of a measure in practice, and there have been several proposals to enrich the simple FRC by taking into account the cycles that contain an edge \cite{samal2018comparative, serrano2022comparative, weber2018coarse}\footnote{This operation is reminiscent of measures for the density around an edge, which have been used in Network Science, e.g. to test the strength of weak ties \cite{onnela2007structure}.}.
In this article, we consider the augmented Forman-Ricci curvature (AFRC) introduced in \cite{serrano2022comparative} and study how different augmentations affect the balance between expressivity and computational complexity in the context of community detection. Based on this evaluation, we identify classes of graphs in which AFRC-based community detection can be recommended compared with ORC-based detection.\\


\noindent The main contributions in this article are
\begin{itemize}
    \item We make a detailed empirical and theoretical investigation of the augmentations of Forman-Ricci curvature.
    \item We discuss the relation of augmented Forman-Ricci curvature to Ollivier-Ricci curvature and un-augmented Forman-Ricci curvature.
    \item We provide evidence that AFRC distinguishes communities in a manner comparable to ORC for many model and real-world networks.
    \item We introduce a novel AFRC-based community detection algorithm, with comparable performance to existing ORC-based methods while being computationally cheaper.
    \item We develop the first codebase for fast and efficient computations with augmentations of the Forman-Ricci curvature.\\
\end{itemize}


\noindent The rest of this paper is structured as follows: Section \ref{sec:Definitions} formally introduces ORC, FRC, its augmentations, and possible approximations. Section \ref{sec:CorrCurv} studies how the different curvatures are related in real-world and artificial networks. Section \ref{sec:CurvGaps} studies the curvature gap -- a quantitative indicator for how well a curvature can distinguish between edges within and between communities -- for the ORC, FRC, and AFRC. Section \ref{sec:Algorithm} introduces our AFRC-based community detection algorithm and demonstrates its performance on several model and real-world networks. Section \ref{sec:Conclusion} concludes the article with a summary of the main findings.

\section{Curvature definitions}\label{sec:Definitions}

\noindent Ricci curvature is a central tool in differential geometry that characterizes the local geometric properties of geodesic spaces by relating the local rate of volume growth to geodesic dispersion.  Several discretizations of Ricci curvature have been proposed, including notions by Forman~\cite{forman2003bochner}, Ollivier~\cite{ollivier2007Riccimetricspaces,ollivier2009RicciMarkovchains}, Erbar and Maas~\cite{erbar2012ricci}, Bakry and \'{E}mery~\cite{bakry2014analysis}, and Devriendt and Lambiotte~\cite{devriendt2022discrete}.  In this paper, we focus on edge-level Ollivier and Forman curvature.  For the latter, we analyze \emph{augmented} Forman curvature, which explicitly encodes the impact of cycles that contain the respective edge. Variants of augmented Forman curvature have been previously studied in~\cite{jost2021characterizations,samal2018comparative,weber2018coarse}, among others.\\

\noindent \textbf{Notation.} We reserve the letters $v$ and $u$ to denote the vertices (equiv., nodes) of a network (equiv., graph), and the symbols $e$ and $e'$ to denote its edges. A node $u$ is a neighbor of $v$ if $(v,u)$ is an edge of the graph. A cycle (of length $r$) is a sequence of different vertices $\gamma=v_1v_2\dots v_r$ such that all cyclically consecutive vertices are connected by an edge in the graph.

\subsection{Ollivier-Ricci curvature}
Ollivier introduced a notion of discrete Ricci curvature for metric spaces in \cite{ollivier2007Riccimetricspaces}; this setting includes graphs with the shortest-path metric between vertices. Aside from a metric, Ollivier's curvature also requires a choice of probability distributions $\mu_v$ -- one for each vertex $v$ in the case of a graph -- which can be thought of as open balls around points in the space. A common choice for graphs is the uniform distribution over the neighbors $\text{neigh}(v)$ of a vertex \cite{sia2019ollivier, tian2023mixed}
\begin{align*}
    \mu_v (u) := \left\{\begin{array}{cl}
    \frac{1}{\deg(v)} & \text{if } u \in \text{neigh}(v) \\[0.5em]
    0 & \text{else}
    \end{array}\right.,
\end{align*}
where $\deg(v)=\vert\text{neigh}(v)\vert$ is the number of neighbors of $v$. With this data, the \emph{Ollivier-Ricci curvature} (ORC) $\mathcal{O}(e)$ of an edge $e=(v_1,v_2)$ is defined as
\begin{equation*}
\mathcal{O}(e) := 1 - \frac{W_1(\mu_{v_1},\mu_{v_2})}{d(v_1,v_2)},
\end{equation*}

\noindent where $d$ is the shortest-path distance between vertices (in this case $=1$) and $W_1$ denotes the Wasserstein distance (also called earth mover's distance in discrete settings) between measures on a metric space; see \cite{ollivier2007Riccimetricspaces, ollivier2009RicciMarkovchains} for details.\\

\noindent Geometrically, the ORC compares the distance between the end points of an edge to the distance between the balls centered at said end points. If $v_1$ and $v_2$ share many neighbors or have many edges between their neighborhoods, the distance between $\mu_{v_1}$ and $\mu_{v_2}$ will be small and so $\mathcal{O}(e)\approx 1$. On the other hand, if the neighborhoods of $v_1$ and $v_2$ are poorly connected, the distance between $\mu_{v_1}$ and $\mu_{v_2}$ will be large and $\mathcal{O}(e)$ will be small and potentially negative. In the context of network community structure, this explains why we expect larger ORC within than between communities: the neighborhoods of two vertices at the ends of an edge are better connected if the edge is contained in a single community. In particular, the threshold value 0 has been previously used in ORC-based community detection to discriminate between inter- and intra-community edges, see \cite{sia2019ollivier} for a more detailed discussion.

\subsection{Forman-Ricci curvature and augmentations}
\begin{figure}[h!]
    \centering
    \includegraphics[width=0.8\textwidth]{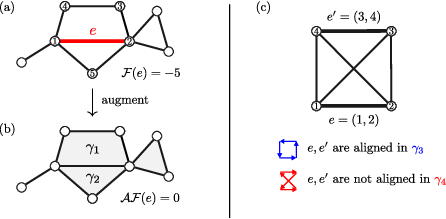}
    \caption{(a) Example of Forman-Ricci curvature. The edge $e=(1,2)$ has FRC equal to $\mathcal{F}(e)=4-\deg(1)-\deg(2)=-5$. (b) Example of augmented Forman-Ricci curvature. The edge $e$ is contained in the cycles $\gamma_1=1234$ and $\gamma_2=125$ and the AFRC of $e$ is equal to $\af(e)=2+2-3-1=0$. This augmentation arises from including faces corresponding to $\gamma_1$ and $\gamma_2$ and then applying Forman's formula to the resulting 2D complex. This example is worked out in detail in Appendix \ref{appendix: detailed example}. (c) Let the vertices be ordered as $1<2<3<4$ and thus $e=(1,2)<(3,4)=e'$. The edges $e,e'$ are contained in cycles $\gamma_3=1234$ and $\gamma_4=1243$. In $\gamma_3$ both edges are traversed from the small to large vertex ($1\rightarrow 2$ and $3\rightarrow 4$), so they are aligned and $\Gamma_{ee'}=1$. In $\gamma_2$ edge $e$ is traversed from small to large ($1\rightarrow 2$) while edge $e'$ is traversed from large to small ($4\rightarrow 3$), so the edges are not aligned in $\gamma_4$ and $\Gamma_{e'e}=1$.} \label{fig: example of AFRC and cycles in AFRC}
\end{figure}

Forman developed a notion of discrete Ricci curvature for weighted CW complexes in \cite{forman2003bochner} in the context of algebraic topology. A CW complex is a type of topological space that includes undirected graphs as a one-dimensional case. A direct translation of Forman's definition to unweighted graphs without multi-edges, as derived in \cite{sreejith2016forman}, gives the \textit{Forman-Ricci curvature} (FRC) $\mathcal{F}(e)$ of an edge $e=(v_1,v_2)$
\begin{align}\label{eq:FRC}
    \mathcal{F}(e):= 4 - \deg(v_1) - \deg(v_2).
\end{align}
While this definition is well-supported by Forman's work and is easy to calculate, it is too simple a measure for many practical and theoretical purposes. For instance, $\mathcal{F}(e)$ does not depend on the number of triangles in which $e$ is contained, which is an important local geometric property of graphs \cite{jost2013ollivier}.
\\
~
\\
To overcome the limitations of the simple formula for FRC, definition \eqref{eq:FRC} can be extended to take into account more of the graph structure around a given edge. The idea is to start from the graph and add solid 2D faces in the spaces enclosed by cycles; see Figure \ref{fig: example of AFRC and cycles in AFRC} for an illustration. This ``augmented graph" is still a CW complex -- the new faces are two-dimensional cells with one-dimensional boundaries given by the cycles -- and one can apply Forman's original definition to get a discrete edge curvature that also depends on the cycles that contain an edge. The following formula for the \emph{augmented Forman-Ricci curvature} (AFRC) $\af(e)$ of an edge $e$ was derived in \cite{serrano2022comparative} based on Forman's definition \cite[Thm. 2.2]{forman2003bochner}:
\begin{equation}\label{eq:AFRC}
\af(e) = 2 + \Gamma_{ee} - \sum_{e'\sim e} \left\vert \Gamma_{ee'}+\Gamma_{e'e}-1\right\vert - \sum_{e'\not\sim e}\vert\Gamma_{ee'}-\Gamma_{e'e}\vert,
\end{equation}
where $e\sim e'$ (resp. $e\not\sim e'$) denotes two distinct edges that share (resp. do not share) a vertex, and where $\Gamma$ is a matrix (defined in \eqref{eq:definition counts} below) of size equal to the number of edges and whose entries count the cycles that contain $e$ and $e'$. The counting goes as follows: fix an arbitrary order on the vertices of the graph. Two edges $e,e'$ are called \emph{aligned} in a cycle $\gamma$ if the cycle traverses both edges from smallest to largest vertex (with respect to the order) or both from largest to smallest vertex. The entries of $\Gamma$ are then defined as
\begin{align}
\Gamma_{ee'} &:= \begin{cases}
\text{\# cycles in which $e$ and $e'$ are aligned} \quad\text{(if $e<e'$)}\\
\text{\# cycles in which $e$ and $e'$ are not aligned} \quad\text{(if $e>e'$)}\label{eq:definition counts}
\end{cases}
\\
\Gamma_{ee} &:= \text{\# cycles that contain $e$}.\nonumber
\end{align}
Here, the order $e<e'$ on edges follows from the order on the vertices\footnote{If we think of an edge as a two-letter word with letters in the (ordered) alphabet of vertices, the induced order on the edges is simply the familiar alphabetical order on the words. See Appendix \ref{appendix: detailed example} for an example.}, and we note that the value of $\af(e)$ is independent of the choice of vertex order. Figure \ref{fig: example of AFRC and cycles in AFRC} illustrates (non-)alignment of edges in a cycle and the resulting distinction between $\Gamma_{ee'}$ and $\Gamma_{e'e}$. This distinction is necessary when translating Forman's definition, since it relies on \textit{oriented} 2D faces. We note that when calculating $\af$ for all edges, the count matrix $\Gamma$ only needs to be constructed once, by iterating over all cycles in the graph. Appendix \ref{appendix: detailed example} contains a detailed example with the above definitions.\\

\noindent In definition \eqref{eq:AFRC} we may choose to only take a subset of all cycles into account; in particular, we will consider the following class of curvatures:
$$
\AF{n}(e) := \af(e)\text{~only taking into account cycles of length $\leq n$.}
$$
Restricting AFRC to cycles up to a certain length is supported by the idea that Ricci curvature should be a local property. This is also satisfied by other notions of discrete curvature such as ORC and Bakry-\'{E}mery curvature, which only depend on cycles of length $n\leq 5$. In this article we mainly work with $\AF{n}$ with $n\in \lbrace 3,4\rbrace$, which are the first augmentations of FRC that give a non-trivial improvement over FRC in a number of applications while still being computationally efficient; in addition to being more expensive, we have found $\AF{5}$ to be less effective in community detection -- see Appendix Figures \ref{fig:AppC_1} and \ref{fig:AppC_2}. We believe the poorer performance is due to the higher number of 5-cycles that contain between-community edges, which blurs the curvature difference between inter- and intra-community edges.
\\
~
\\
For certain choices of cycles, definition \eqref{eq:AFRC} for the AFRC reduces to known extensions of FRC. For instance, the curvature $\mathcal{F}(e)+3\cdot\#\lbrace \text{3-cycles $\ni e$}\rbrace = \AF{3}(e)$ was proposed in \cite{samal2018comparative} as a variant of Forman's curvature in graphs and networks. In \cite{tee2021enhanced, weber2018coarse}, it was observed that in graphs where any two cycles share at most one edge -- these graphs are called quasi-convex in \cite{forman2003bochner} and independent short-cycle graphs in \cite{tee2021enhanced} --, Forman's curvature simplifies to $\mathcal{F}(e)+\sum_{n\geq 3}(6-n)\cdot \text{\#$\lbrace \text{$n$-cycles}\ni e\rbrace$}$. This expression can provide a simple estimate for edges whose neighborhoods are expected to be ``almost quasi-convex'', but in general it differs from AFRC. In \cite{jost2021characterizations}, Jost and M\"{u}nch studied different choices of weighted cycles when augmenting a graph.
\\
~
\\
Finally, we present a rewriting of equation \eqref{eq:AFRC} as an explicit augmentation of FRC:
\begin{equation}\label{eq:AFRC2}
\af(e) = \mathcal{F}(e) - \Gamma_{ee} + 2\cdot\sum_{e'\sim e}\min\lbrace (\Gamma_{ee'}+\Gamma_{e'e}),1\rbrace - \sum_{e'\not\sim e}\vert \Gamma_{ee'}-\Gamma_{e'e}\vert.    
\end{equation}
Here the second term again counts the number of cycles that contain $e$, the third term counts the number of edges $e'\sim e$ that share a cycle with $e$ and the fourth term again counts the aligned and unaligned cycles through $e$ and $e'\not\sim e$. 
\\
~
\\
Equation \eqref{eq:AFRC2} reveals how the addition of cycles affects the value of AFRC compared to FRC. Including cycles of length $3$, we find that $\AF{3}(e)=\mathcal{F}(e)+3\cdot\#\lbrace \text{$3$-cycles $\ni e$}\rbrace$ and thus that AFRC is higher for edges whose end points share many neighbors. In networks with community structure, we thus expect $\AF{3}$ to be large within communities and small (possibly negative) between communities, similar to ORC. Including cycles of length $4$ and above, however, we expect a different behavior. Now the positive term in \eqref{eq:AFRC2} due to the neighbors of $e$ is limited by the degrees of the end points of $e$, while the negative terms ($\Gamma_{ee}$ and the sum over non-neighbors) can become increasingly negative as more cycles are included. As a result, we can expect edges within communities to have \emph{lower} $\AF{n\geq 4}$ values compared to edges between, opposite to what is observed for ORC and $\AF{3}$; in later experiments on networks with community structure, we confirm this analysis. Finally, we remark that un-augmented FRC (\ref{eq:FRC}) does not distinguish edges within and between communities when different communities have the same typical degrees.\\ 

\noindent We note that, to the best of our knowledge, there is no natural threshold to discriminate edges according to their AFRC (unlike the curvature value 0 for ORC). In fact, our experiments show that the potential curvature threshold varies drastically between networks -- see e.g. Figure \ref{fig:curvature_switch}; we discuss how to choose an adequate case-by-case value in Section \ref{sec:Algorithm}. 

\section{Correlations between curvatures}\label{sec:CorrCurv} Since ORC and AFRC are discretizations of the same smooth concept, their comparison is a natural and important question.\\

\noindent On the theoretical side, \cite{tee2021enhanced} found evidence that ORC and $\AF{5}$ are equivalent in quasi-convex graphs and \cite{jost2021characterizations} found that the maximal Forman curvature over all weighted augmentations of a graph is equal to ORC. In general, the relation between different notions of discrete curvature remains an important open question.\\

\noindent On the experimental side, Samal et al \cite{samal2018comparative} extensively studied the correlation between ORC and both FRC and $\AF{3}$ in model and real-world networks. In their results, ORC tended to greatly correlate with both curvatures, and it generally correlated higher with $\AF{3}$. We extend this study here by including correlations of these curvatures with $\AF{4}$. Tables \ref{table:1} and \ref{table:2} show the Pearson correlation coefficient between the different curvature measures in, respectively, model and real-world networks.\\

\noindent \textit{\textbf{Model networks:}} We study four different network models: the Erd\H{o}s-R\'{e}nyi (ER) graph, which is a completely homogeneous model and should not exhibit community structure; a bipartite version (BG), which does not contain any 3-cycles and thus $\mathcal{F}(e)=\AF{3}(e)$ for all edges; and the stochastic block model (SBM) and hierarchical bipartite graph (HBG), which do have ground-truth communities. Details on these network models can be found in Appendix \ref{appendix: studied networks}.
\\

\noindent Table \ref{table:1} shows that the FRC vs $\AF{3}$ and $\AF{4}$ correlations exhibit the expected behavior for model networks: the fewer edges -- and thus cycles -- present in a network, the more similar the curvatures are (as they are equal in a network with no cycles). In particular, note the correlation of 1 between FRC and $\AF{3}$ for the bipartite networks -- as expected from networks with no 3-cycles. The only outlier of this behavior is the HBG, as the correlation FRC vs $\AF{4}$ seems to increase with the number of edges.\\

\noindent Regarding ORC vs the different versions of AFRC, all exhibit a similar behavior on the ER and BG graphs, being quite correlated in sparser networks and losing this correlation as the connection probability $p$ increases -- although in general ORC appears best correlated to $\AF{4}$. This same behavior can be observed in the SBM for the ORC vs FRC correlation, whereas the ORC vs $\AF{3}$ and $\AF{4}$ are notably distinct. We observe how the correlation ORC vs $\AF{3}$ increases positively with the community size, whereas ORC vs $\AF{4}$ grows negatively with it. This behavior agrees with our discussion at the end of Section \ref{sec:Definitions}, and is illustrated later in Figure \ref{fig:curvature_switch}, and Appendix Figures \ref{fig:AppC_1}-\ref{fig:AppC_2}. Finally, in the HBG the ORC vs $\AF{4}$ correlation initially shows a similar behavior to that of the SBM, but we see that they become uncorrelated as the inter-community connection probability $q$ grows (while the intra-community probability $p$ is held fixed). We will further explore this last effect in Section \ref{sec:CurvGaps}, where we show that ORC does not appear to detect communities in HBGs with noisy structure - and so it is natural for ORC and $\AF{4}$ not to be correlated when only one depends on the communities of the network.\\

\begin{table}[h!]\label{tab:CurvCorrModel}
\centering
\begin{tabular}{ |l|c|c|c|c|c| } 
\hline

\textbf{Network} & OR, FR & OR, $\AF{3}$ & OR, $\AF{4}$ & FR, $\AF{3}$ & FR, $\AF{4}$\\
\hline
ER($1000$, $0.003$) & 0.84 & 0.87 & 0.87 & 0.99 &  0.98\\ 
\hline
ER($1000$, $0.007$) & 0.42 & 0.48 & 0.59 & 0.98 &  0.94\\ 
\hline
ER($1000$, $0.01$) & -0.01 & 0.12 & 0.30 & 0.98 & 0.89\\ 
\hhline{|=|=|=|=|=|=|}
SBM($10$, $5$, $0.7$, $0.05$) & 0.27 & 0.79 & 0.76 & 0.64 &  0.59\\ 
\hline 
SBM($10$, $10$, $0.7$, $0.05$) & 0.07 & 0.83 & -0.66 & 0.49 & 0.31\\ 
\hline 
SBM($10$, $15$, $0.7$, $0.05$) & 0.08 & 0.90 & -0.84 & 0.28 & 0.19\\ 
\hline 
SBM($10$, $20$, $0.7$, $0.05$) & -0.03 & 0.90 & -0.91 & 0.33 & 0.19\\ 
\hhline{|=|=|=|=|=|=|} 
BG($50$, $0.03$) & 0.74 & 0.74 & 0.84 & 1.0 & 0.91\\
\hline
BG($50$, $0.07$) & 0.01 & 0.01 & 0.58 & 1.0 & 0.69\\
\hline
BG($50$, $0.1$) & 0.29 & 0.29 & 0.27 & 1.0 & 0.79\\
\hhline{|=|=|=|=|=|=|}
HBG($50$, $0.5$, $0.05$) & -0.01 & -0.01 & -0.72 & 1.0 & 0.55\\
\hline
HBG($50$, $0.5$, $0.1$) & -0.11 & -0.11 & -0.25 & 1.0 & 0.61\\
\hline
HBG($50$, $0.5$, $0.15$) & -0.01 & -0.01 & 0.00 & 1.0 & 0.74\\
\hline
\end{tabular}
\caption{Comparison of Ollivier-Ricci curvature (OR) with Forman-Ricci curvature (FR) and augmentations of the Forman-Ricci curvature ($\AF{3}$ and $\AF{4}$) of edges in artificial networks. In this table, we list the Pearson correlation between the edge curvatures.}
\label{table:1}
\end{table}

\begin{table}[h!]
\centering
\begin{tabular}{ |l|c|c|c|c|c| } 
\hline
\textbf{Network} & OR, FR & OR, $\AF{3}$ & OR, $\AF{4}$ & FR, $\AF{3}$ & FR, $\AF{4}$\\
\hline
Dolphin & 0.08 & 0.75 & -0.15 & 0.25 & 0.77\\ 
\hline
Power Grid & 0.48 & 0.72 & 0.61 & 0.82 & 0.89\\ 
\hline
Word Adjacency & -0.04 & 0.18 & -0.24 & 0.89 & 0.89\\ 
\hline
Southern Women & 0.47 & 0.47 & -0.01 & 1.0 & 0.71\\ 
\hline
Corporate Interlocks & 0.48 & 0.48 & -0.05 & 1.0 & 0.64\\ 
\hhline{|=|=|=|=|=|=|}
Karate Club & 0.56 & 0.81 & 0.58 & 0.62 & 0.84 \\ 
\hline
College Football & -0.09 & 0.94 & -0.76 & -0.07 & 0.37\\ 
\hline
School 1 & -0.01 & 0.84 & -0.26 & 0.38 & 0.94\\ 
\hline
School 2 & 0.02 & 0.81 & -0.27 & 0.36 & 0.92\\ 
\hline
\end{tabular}
\caption{Comparison of Ollivier-Ricci curvature (OR) with Forman-Ricci curvature (FR) and augmentations of the Forman-Ricci curvature ($\AF{3}$ and $\AF{4}$) of edges in real-world networks. The two sections of the table correspond respectively to networks without and with community metadata. In this table, we list the Pearson 
correlation between the edge curvatures.}
\label{table:2}
\end{table}

\noindent \textit{\textbf{Real-world networks:}} We study a variety of networks, including two bipartite (Southern Women and Corporate Interlocks) and four with community structure metadata (second half of Table \ref{table:2}). All networks have between 10 and 250 nodes, except for the US Power Grid, which has around 5000. Detailed descriptions of the networks can be found in Appendix \ref{appendix: studied networks}.\\

\noindent As seen in Table \ref{table:2}, all networks without a baseline community structure (except for the Power Grid) exhibit similar curvature behaviors. The main observation is that ORC and $\AF{4}$ are slightly negatively correlated, an effect we only observe in the models \textit{with} an underlying community structure. On the other hand, the Power Grid network behaves akin to the ER, with all curvatures moderately correlated. We believe this is related to the network being substantially sparser (average degree $\sim1.3$) and with few 3- and 4-cycles (651 and 979, respectively), which causes all curvatures to behave similarly.\\

\noindent Regarding the networks with ground-truth community structure, all but the Karate Club show a positive (resp. negative) correlation between ORC and $\AF{3}$ (resp. $\AF{4}$), as expected. We attribute the different behavior of the Karate Club to the fact that its ground-truth community structure is slightly blurred in the context of network analysis, which is reflected by the relatively low ratio $\hat{p} / \hat{q} \approx 6.57$, where $\hat{p}$ and $\hat{q}$ are estimates of the probabilities of edges forming within and between communities, respectively. By comparison, this ratio is about $19.51$ in the College Football network. We obtain these estimates by dividing the number of edges within (between) communities by the number of possible edges. For the Karate Club and School networks, we note that $\AF{4}$ is much more highly correlated to FRC than ORC, while $\AF{3}$ exhibits the converse behavior. This could be an indicator that in these cases the simpler augmentation $\AF{3}$ is better for community detection, as ORC is expectedly effective in detecting community structure whereas FRC is not.\\
~
\\
To give a summary of our correlation results, we have experimentally verified the discussion of Section \ref{sec:Definitions}. In particular, we have shown how the amount of edges and cycles of a network have different effects on the various versions of Forman-Ricci curvature, and we have also provided evidence for the expected behaviors of the different curvatures in networks with an underlying community structure. The high correlation between ORC and different versions of AFRC are a first indication that the latter have a potential use in community detection; our results also suggest that the effectiveness of the different augmentations may depend on the characteristics of the network.

\section{Curvature gaps}\label{sec:CurvGaps}
In this section we compare the performance of the studied curvatures in community detection -- both for model and real-world networks with an underlying community structure.\\
 
\begin{figure}[H]
  \begin{subfigure}[t]{.48\textwidth}
    \centering
    \includegraphics[width=\linewidth]{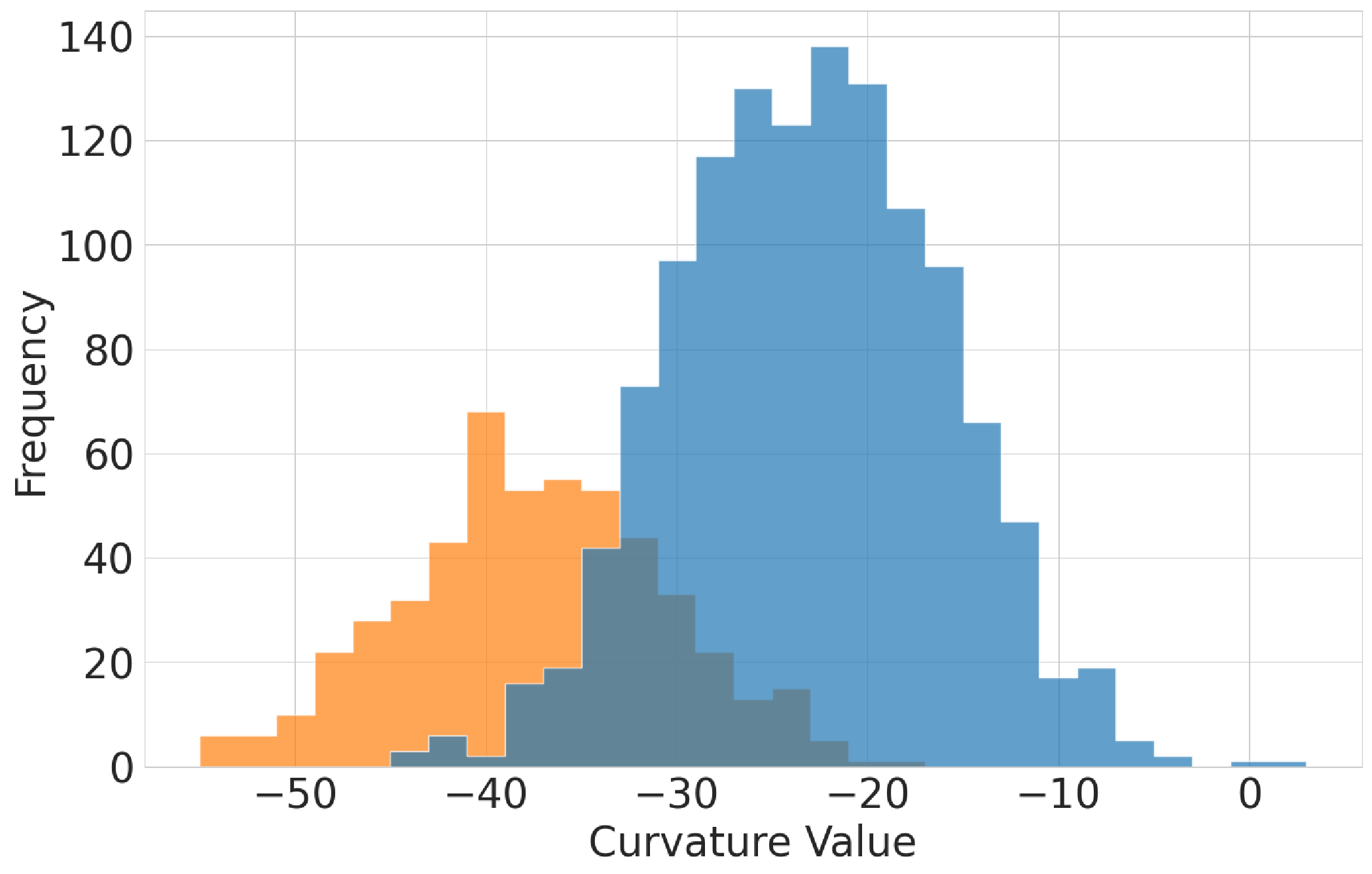}
    \caption{}
  \end{subfigure}
  \hfill
  \begin{subfigure}[t]{.48\textwidth}
    \centering
    \includegraphics[width=\linewidth]{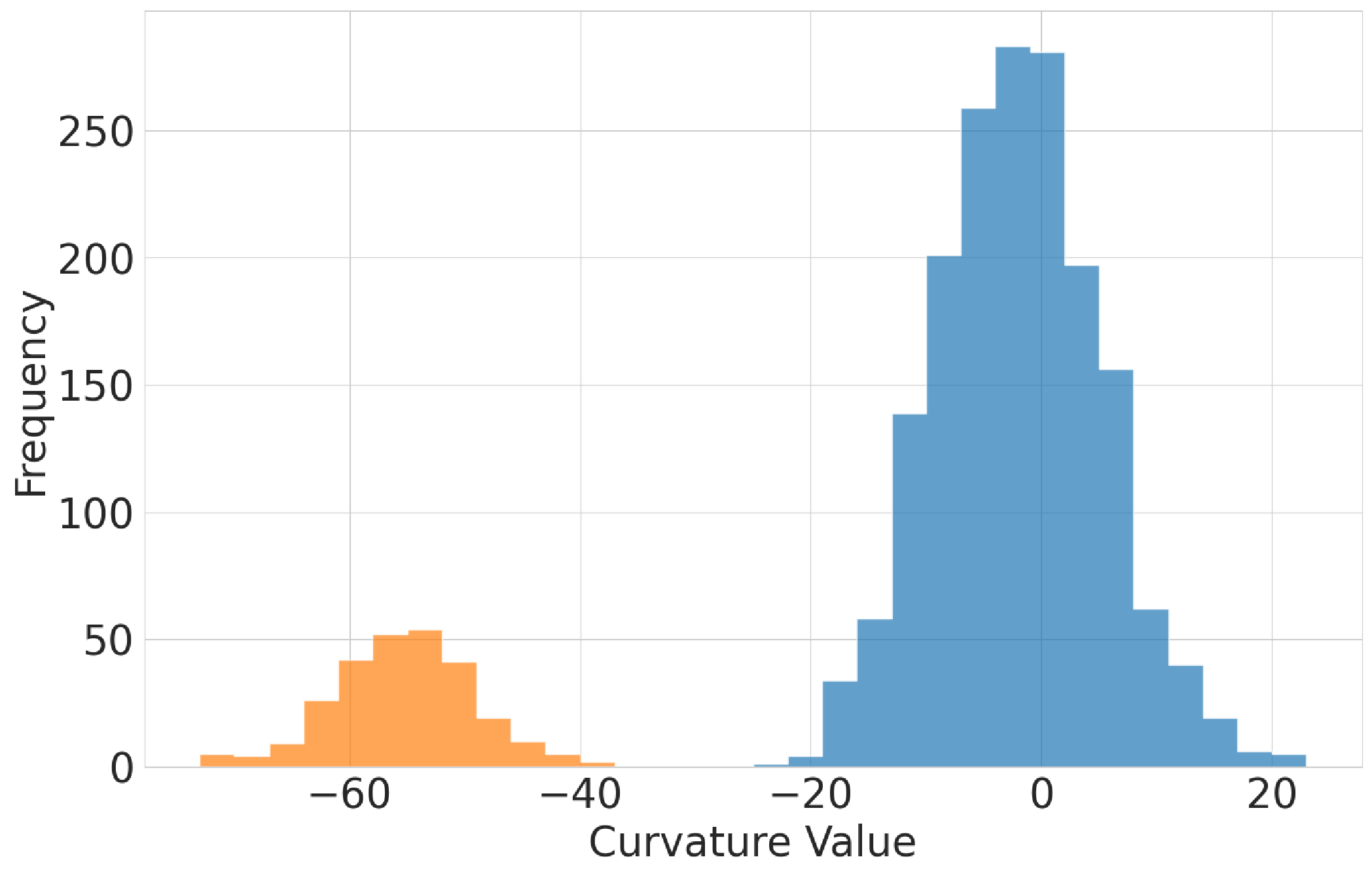}
    \caption{}
  \end{subfigure}

  \medskip

    \begin{subfigure}[t]{.48\textwidth}
    \centering
    \includegraphics[width=\linewidth]{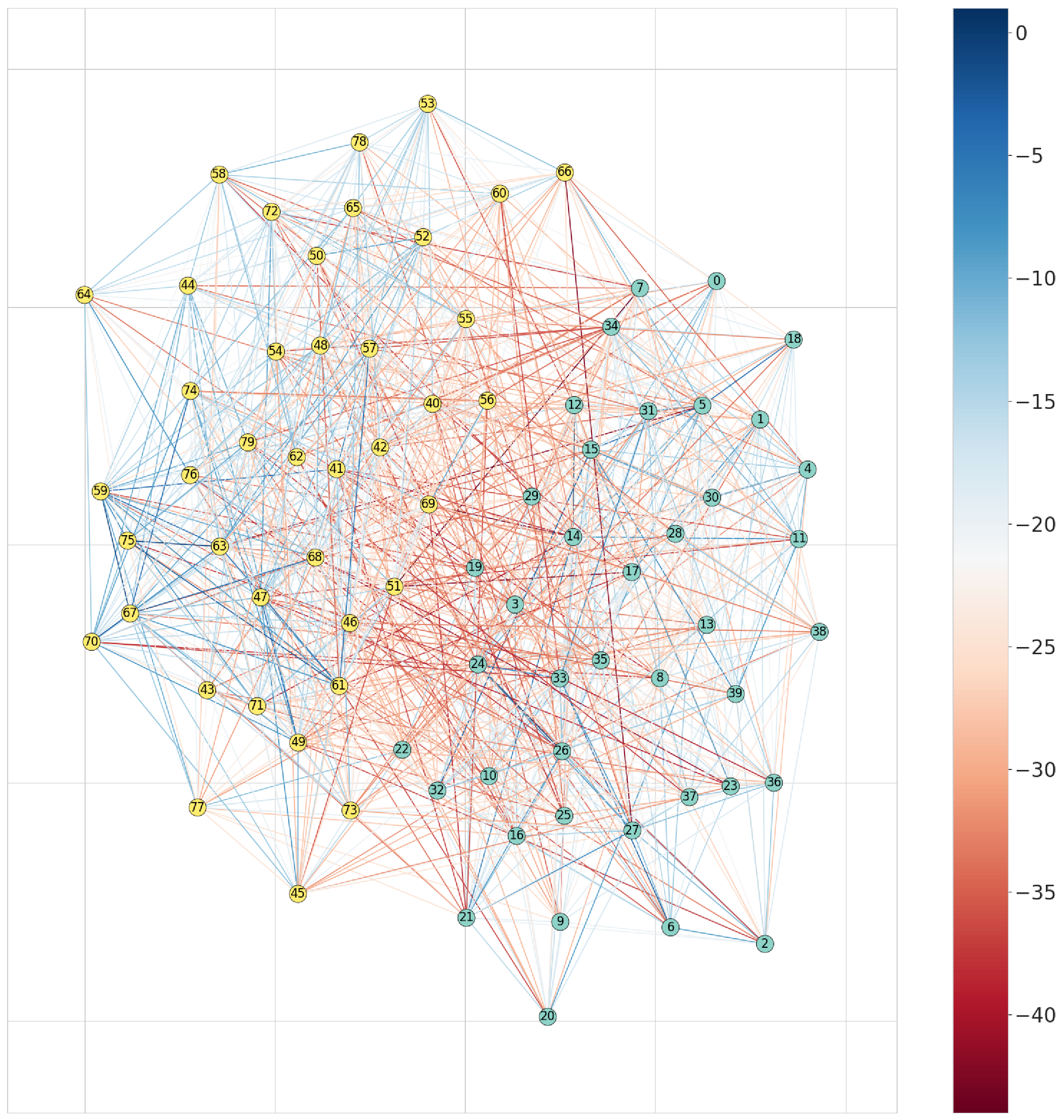}
    \caption{}
  \end{subfigure}
  \hfill
  \begin{subfigure}[t]{.48\textwidth}
    \centering
    \includegraphics[width=\linewidth]{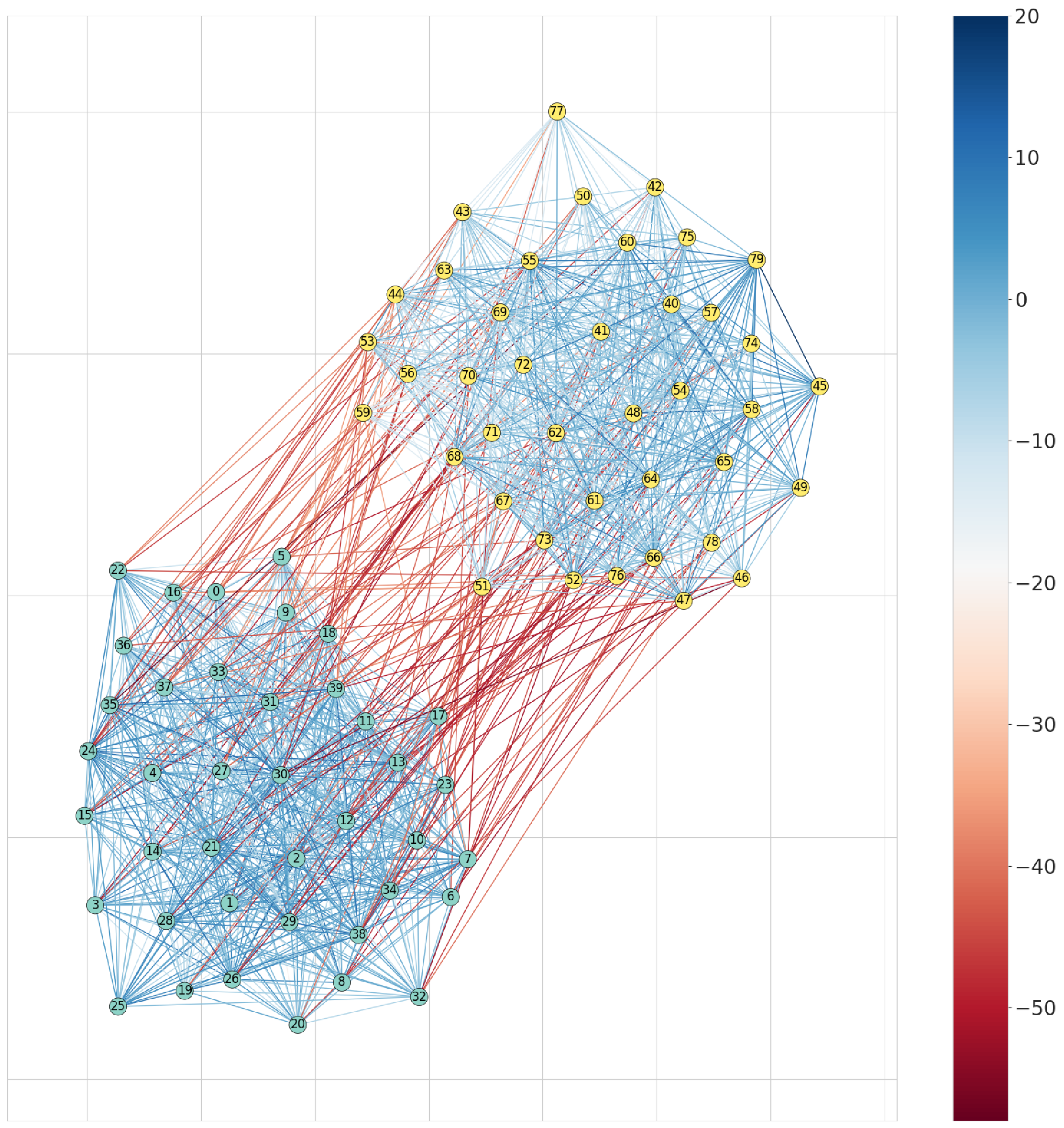}
    \caption{}
  \end{subfigure}

  \caption[Caption]{Distributions of the Augmented Forman-Ricci Curvature ($\AF{3}$) across edges in stochastic block models with weak (\textbf{(a)}, SBM(2, 40, $0.5$, $0.2$)) and strong (\textbf{(b)}, SBM(2, 40, $0.7$, $0.1$)) community structure. The curvature gap captures the (normalized) distance between the means of the edges between (orange) and within (blue) communities. \textbf{(c)} and \textbf{(d)} shows the corresponding networks with edges colored according to curvature.}
  \label{fig:various_curvature_gaps}
\end{figure}

\begin{figure}[H]
  \begin{subfigure}[t]{.24\textwidth}
    \centering
    \includegraphics[width=\linewidth]{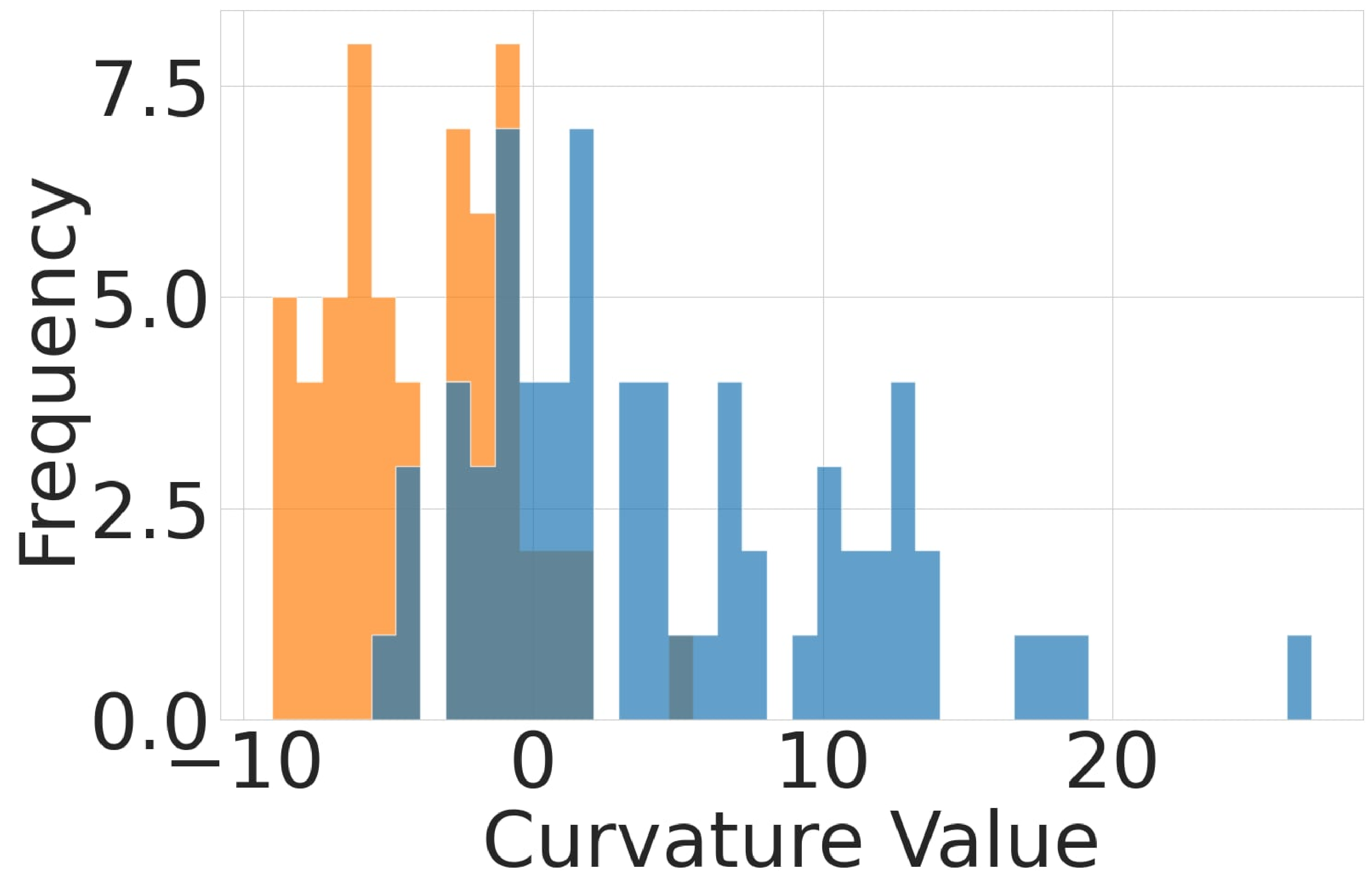}
  \end{subfigure}
  \hfill
  \begin{subfigure}[t]{.24\textwidth}
    \centering
    \includegraphics[width=\linewidth]{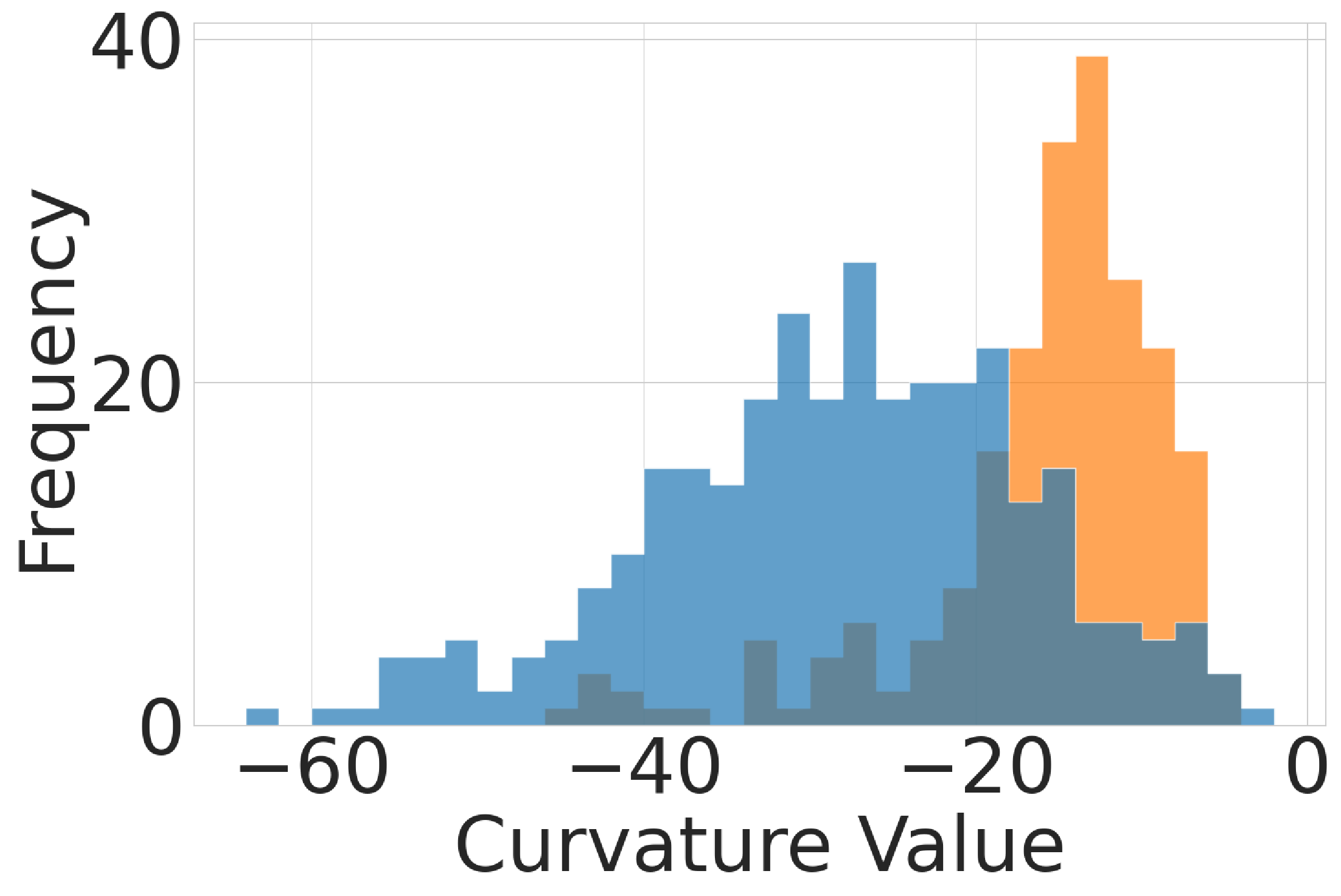}
  \end{subfigure}
    \begin{subfigure}[t]{.24\textwidth}
    \centering
    \includegraphics[width=\linewidth]{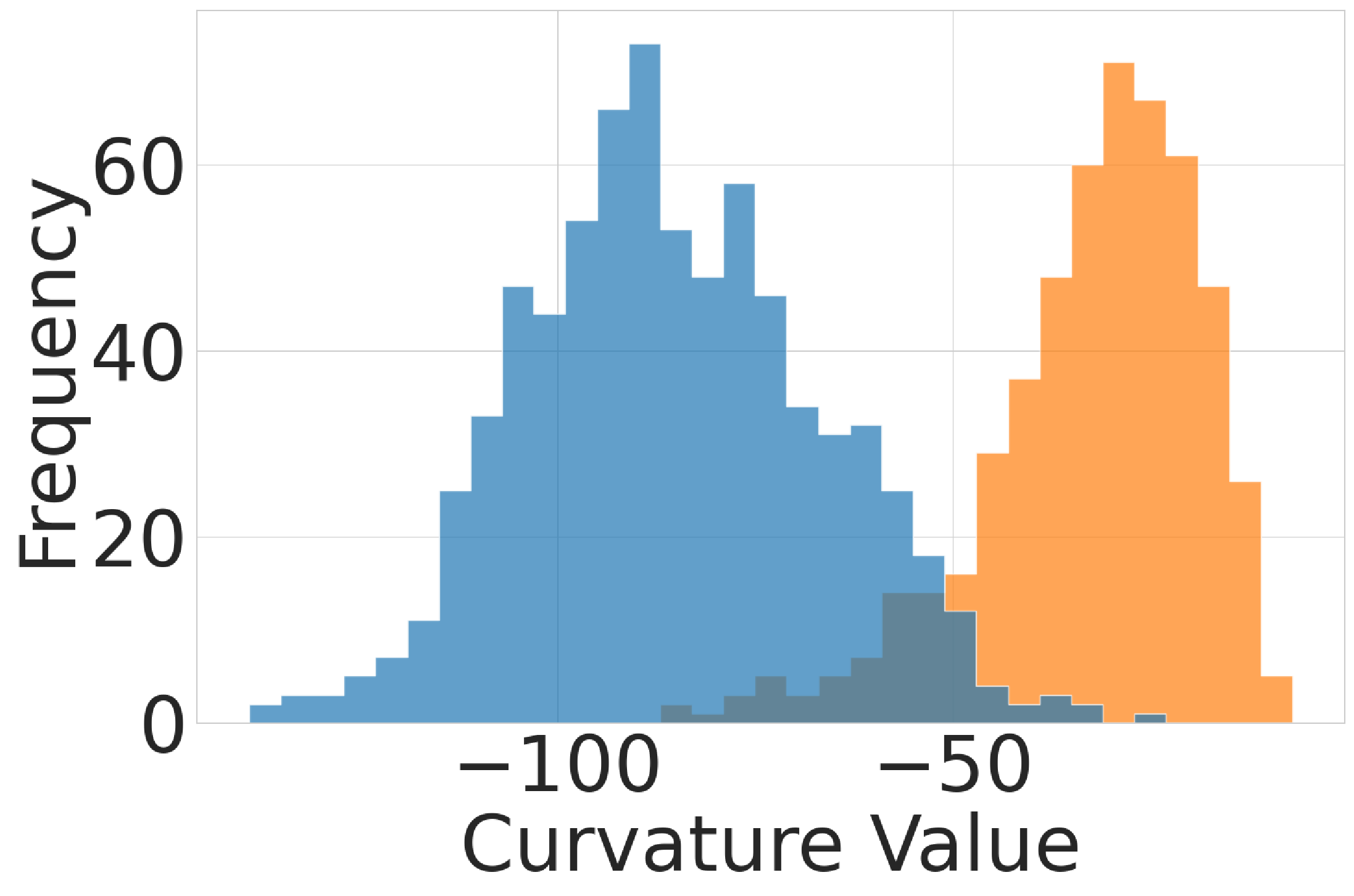}
  \end{subfigure}
  \hfill
  \begin{subfigure}[t]{.24\textwidth}
    \centering
    \includegraphics[width=\linewidth]{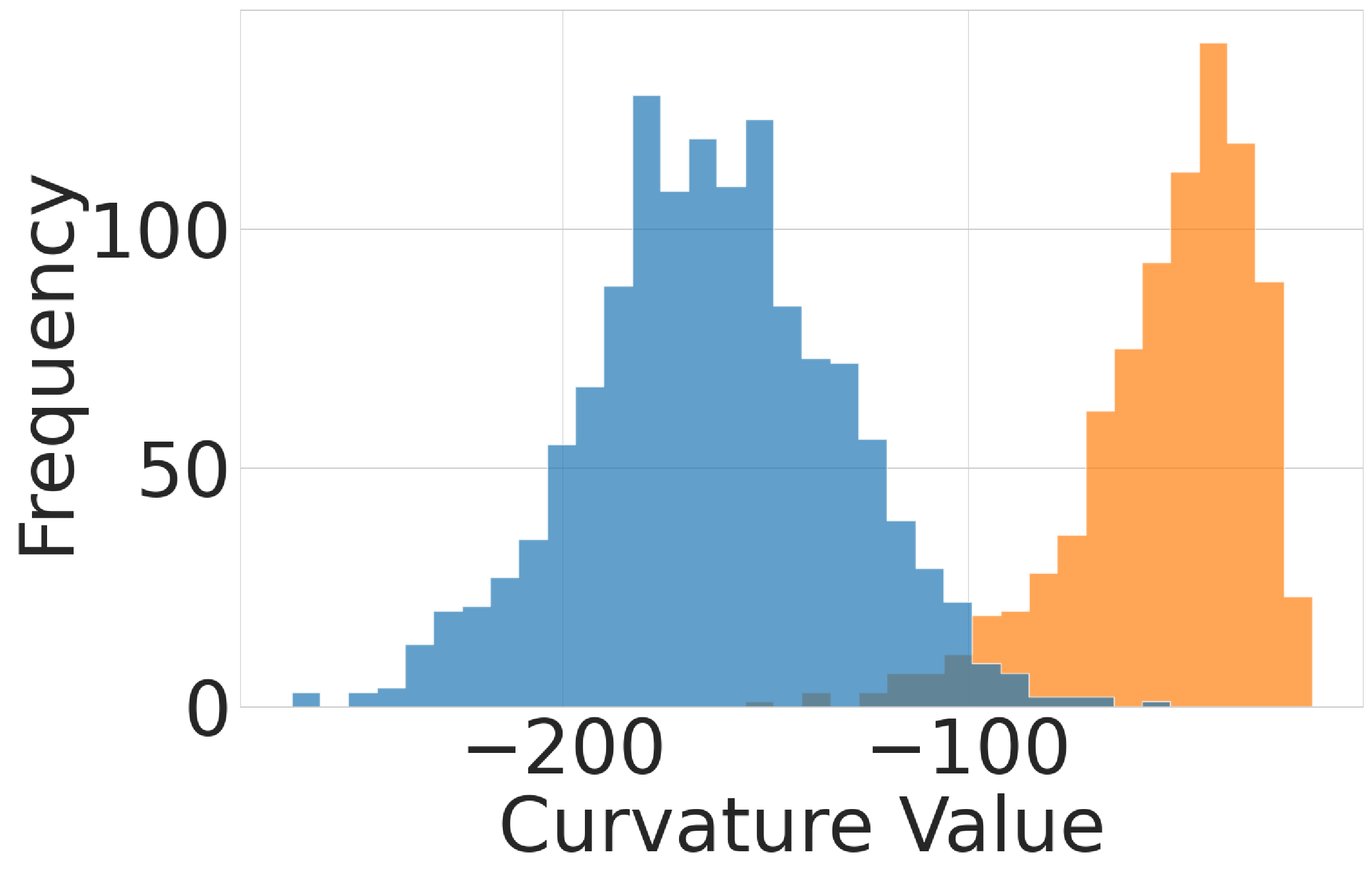}
  \end{subfigure}

  \medskip

  \begin{subfigure}[t]{.24\textwidth}
    \centering
    \includegraphics[width=\linewidth]{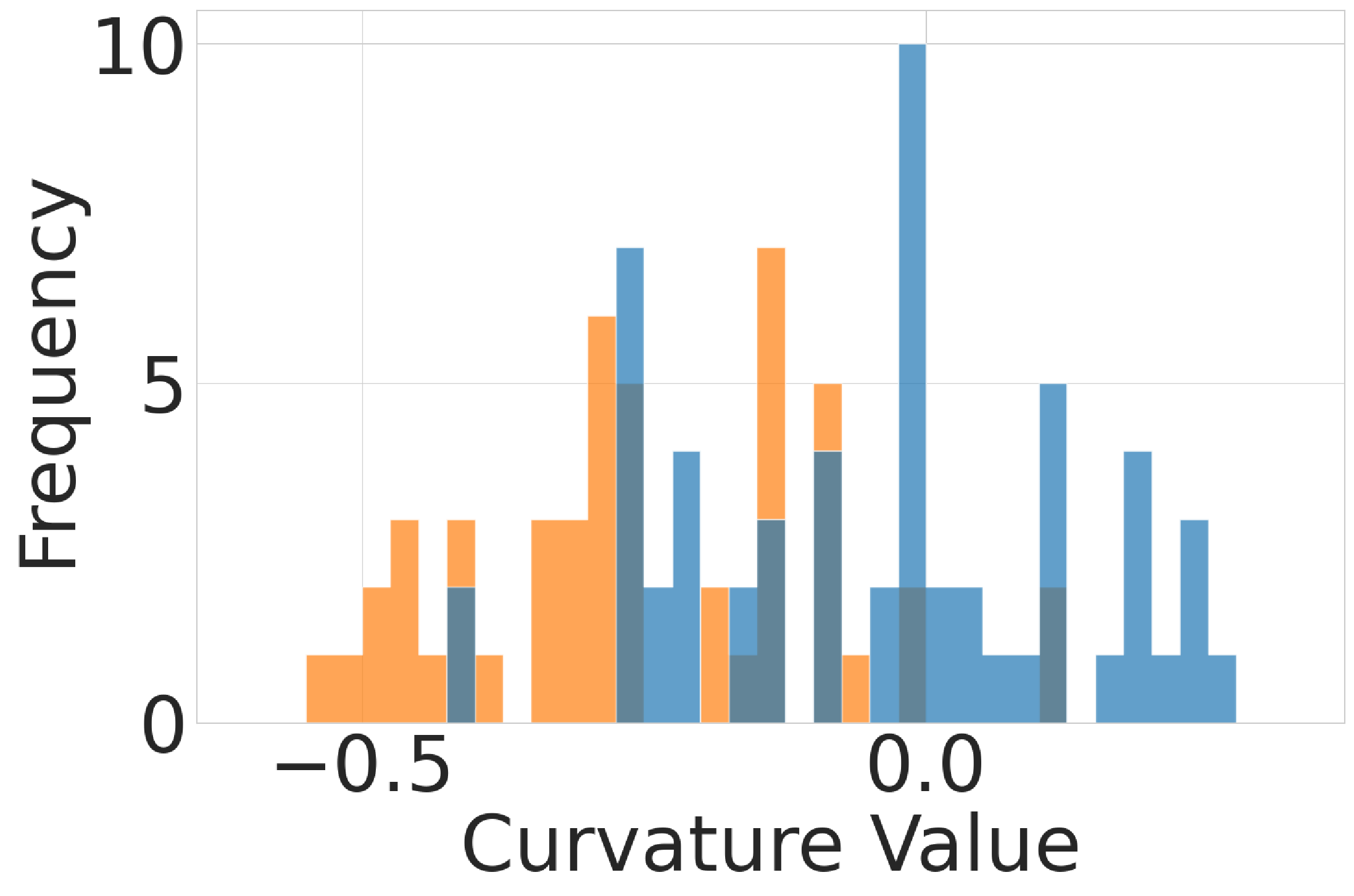}
  \end{subfigure}
  \hfill
  \begin{subfigure}[t]{.24\textwidth}
    \centering
    \includegraphics[width=\linewidth]{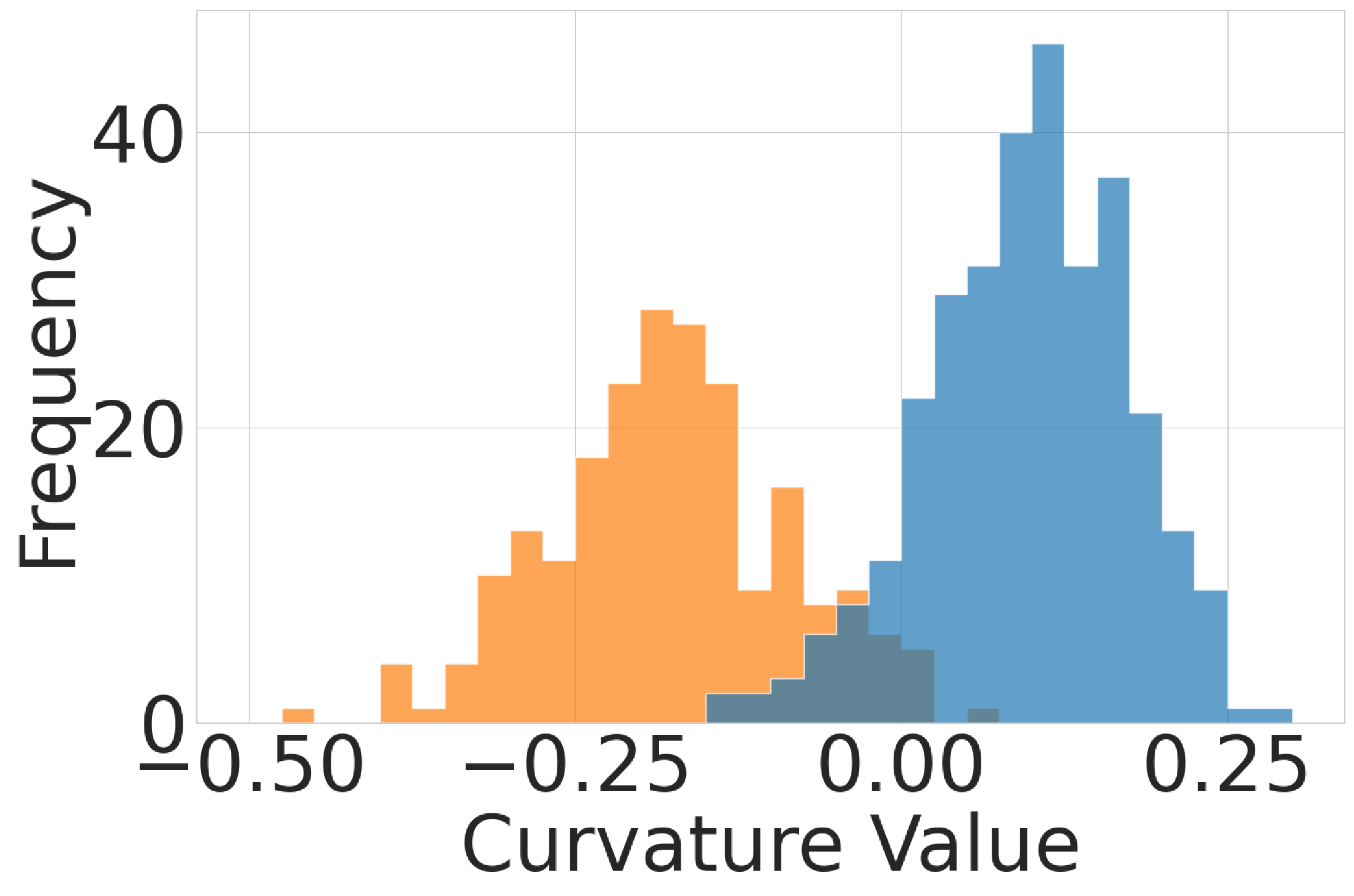}
  \end{subfigure}
    \begin{subfigure}[t]{.24\textwidth}
    \centering
    \includegraphics[width=\linewidth]{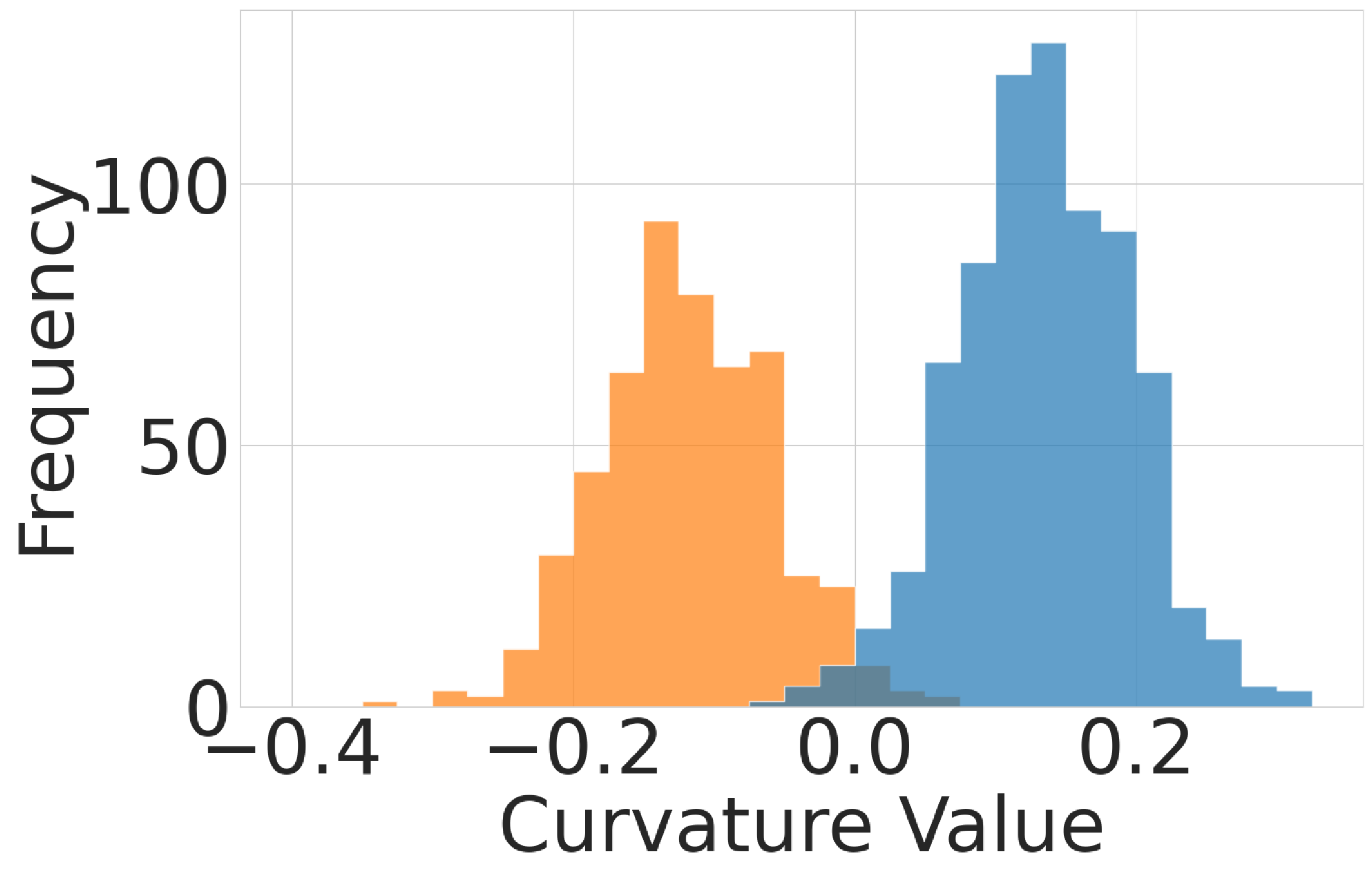}
  \end{subfigure}
  \hfill
  \begin{subfigure}[t]{.24\textwidth}
    \centering
    \includegraphics[width=\linewidth]{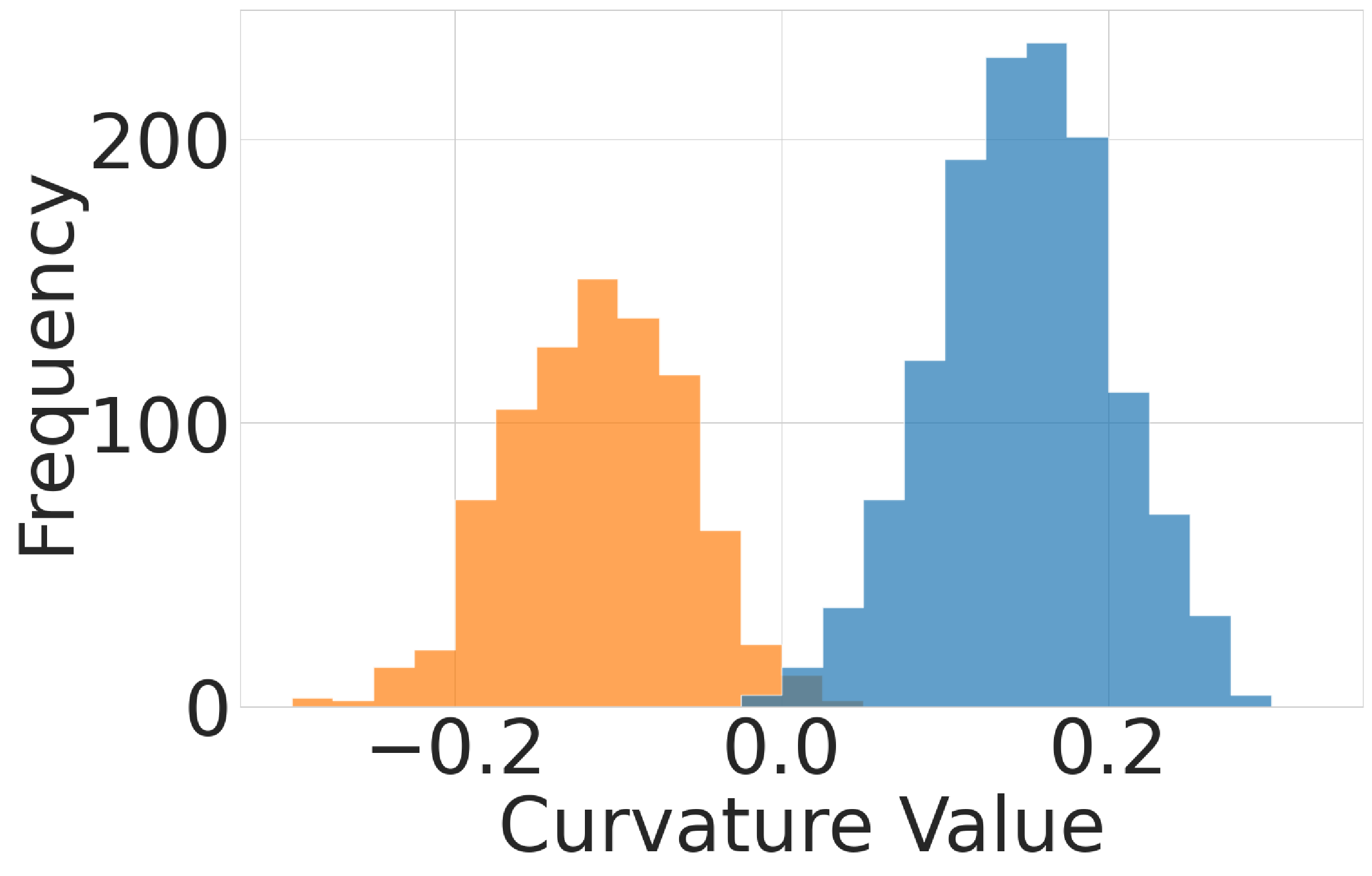}
  \end{subfigure}
  
  \caption{Distributions of $\AF{4}$ \textbf{(top)} and the Ollivier-Ricci Curvature \textbf{(bottom)} for SBM(10, $k$, $0.7$, $0.05$) with $k\in\{5, 10, 15, 20\}$ (from left to right). As the community structure becomes more pronounced with larger community sizes, the curvature gaps increase. Note also the ``switch" in the order of the within-community edges (blue) and the between-community edges (orange) in the AFRC histograms. }\label{fig:curvature_switch}
\end{figure}

\noindent The essential observation for ORC-based community detection, as discussed in Section \ref{sec:Definitions} and used for instance in \cite{sia2019ollivier}, is that the ORC of edges within a community tends to be higher than the ORC of edges between communities. As a consequence, the histogram of curvature values in a network with community structure tends to be bimodal with two separate peaks or clusters. Figure \ref{fig:various_curvature_gaps} shows this effect for $\AF{3}$ in two graphs with a clear community structure. Figure \ref{fig:curvature_switch} furthermore shows that the clusters in the histogram become more pronounced as the community structure becomes stronger. One way to quantify the separation between the two clusters is the \textit{curvature gap}, introduced in \cite{gosztolai2021unfolding}. Given a partition of a network into communities, let $\kappa_\text{within}$ denote the average curvature of all edges contained in a community, $\kappa_\text{between}$ the average curvature of edges between communities, and $\sigma_{\text{within}},\ \sigma_{\text{between}}$ the standard deviation of the respective curvatures. The curvature gap of this partition of the network

\noindent is defined as

$$ \Delta \kappa := \frac{1}{\sigma} \left|\kappa_\text{within} - \kappa_\text{between} \right|,$$

\noindent where $\sigma = \sqrt{\frac{1}{2} \left( \sigma_\text{within}^2 + \sigma_\text{between}^2\right)}$. Intuitively, the resulting number can be interpreted as `the number of standard deviations between the mean curvatures'. Since the curvature gap reflects how well the curvature distinguishes edges within and between communities, we will use it as a proxy for how well any given curvature detects the community structure. Tables \ref{table:3} and \ref{table:4} show the curvature gap obtained from each curvature measure in model and real-world networks, respectively. Importantly, these curvature gaps also depend on how clear the community structure is in a graph (i.e. the ratio

\begin{figure}[H]
  \begin{subfigure}[t]{.44\textwidth}
    \centering
    \includegraphics[width=\linewidth]{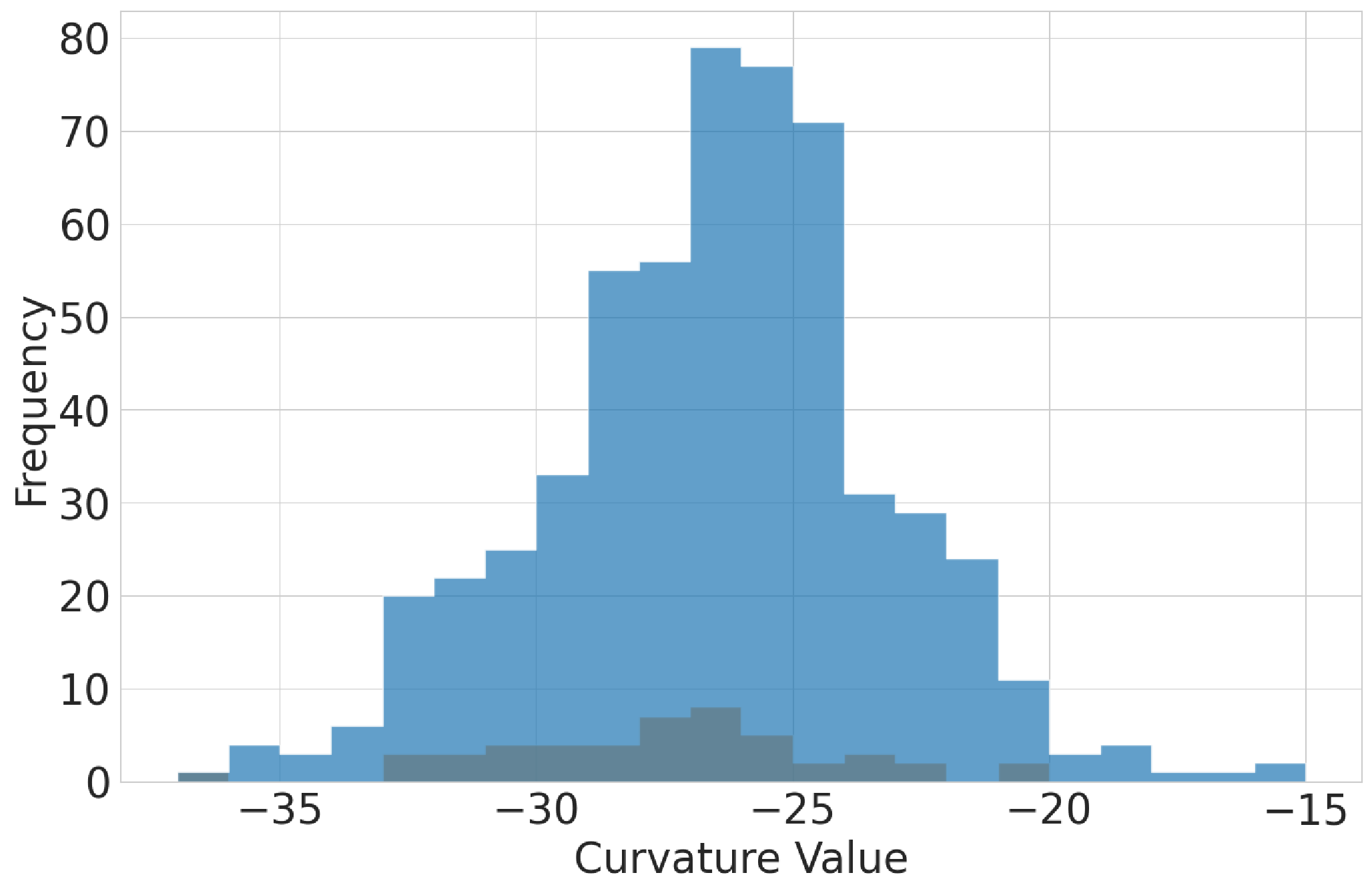}
    \caption{}
    \label{fig:hbg_curvature_gaps_a}
  \end{subfigure}
  \hfill
  \begin{subfigure}[t]{.44\textwidth}
    \centering
    \includegraphics[width=\linewidth]{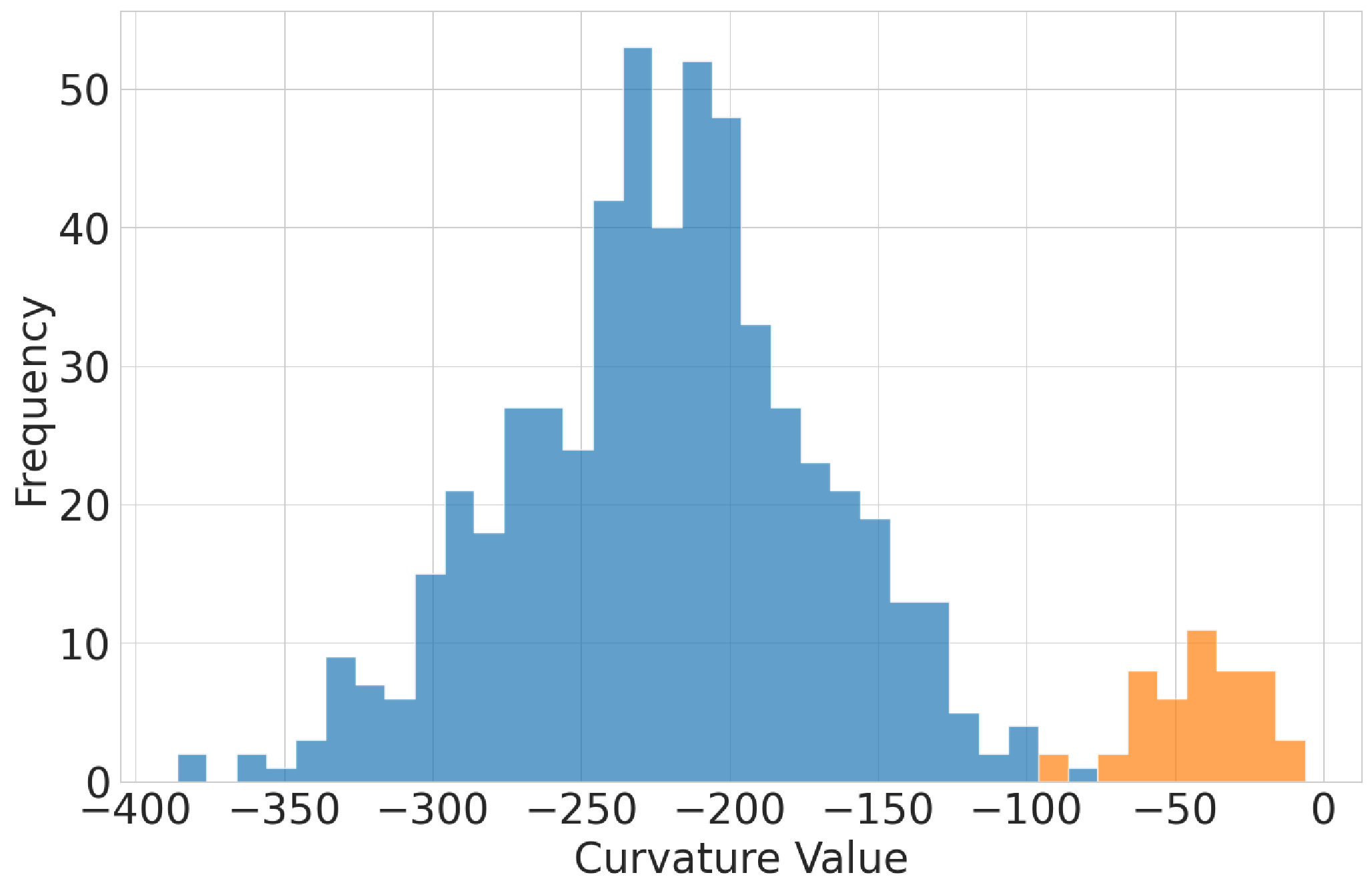}
    \caption{}
    \label{fig:hbg_curvature_gaps_b}
  \end{subfigure}
    
    \caption{Distributions of $\AF{3}$ in \textbf{(a)} and $\AF{4}$ in \textbf{(b)} for a HBG(40, 0.7, 0.05).}
  \label{fig:hbg_curvature_gaps}
\end{figure}

\noindent $p/q$ mentioned before, or the size vs the number of the communities); this dependency is illustrated in Figure \ref{fig:curvature_switch} and Appendix Figures \ref{fig:AppC_1}-\ref{fig:AppC_4}.
\\

\noindent \textbf{\textit{Model networks:}} We study the previously introduced SBM and HBG models. The latter is particularly interesting because its lack of triangles means that $\AF{3}$ no longer augments the FRC. In this case, $\AF{4}$ is the first non-trivial augmentation.\\

\noindent As discussed in Section \ref{sec:Definitions}, FRC should not be able to detect communities, which is confirmed by the very poor curvature gap observed in Table \ref{table:3} for both models. Similarly as expected, $\AF{3}$ performs poorly detecting communities in bipartite networks. This is also illustrated in Figure \ref{fig:hbg_curvature_gaps}.\\

\noindent On the other hand, ORC attains good curvature gaps. As noted, we observe that the gap depends on how clear the community structure is: it decreases when we reduce the number of inside-community edges or increase the number of between-community edges; moreover, the gap increases with the size of the communities of the SBM. Finally, $\AF{4}$ (and $\AF{3}$ for the unipartite model) performs similarly to ORC: it shows a relevant gap in all the studied models, which also depends on the clearness of the community structure. We observe how ORC and $\AF{3}$ tend to perform similarly and outperform $\AF{4}$ in the studied SBMs (although the latter still attains curvature gaps of a similar magnitude), whereas in the studied HBGs $\AF{4}$ attains relevant curvature gaps even in the cases where ORC performs poorly --- see Appendix Figures \ref{fig:AppC_3}, \ref{fig:AppC_4} for some plots of this phenomenon.\\

\noindent \textit{\textbf{Real-world networks:}} Next, we study the curvature gap in a range of real-world network with a known community structure. We remark that all networks studied in this section are unipartite.\\

\noindent As with the model networks, FRC performs poorly in all of the studied networks, and both ORC and the augmentations of FRC offer more relevant curvature gaps; $\AF{3}$ offers the best overall performance in our network sample, and $\AF{4}$ performs worst (see Table \ref{table:4}).  Similarly to Section \ref{sec:CorrCurv}, the higher the $\hat{p} / \hat{q}$ ratio of the ground-truth communities, the bigger the curvature gap: the Football network ($\hat{p} / \hat{q}\approx 19.51$) presents a substantial curvature gap for both ORC and the FRC augmentations, whereas for the Karate ($\hat{p} / \hat{q}\approx 6.57$) and School networks ($\hat{p} / \hat{q}\approx 5.89$ and $\approx 6.76$, resp.) the curvature gaps are smaller -- notably, $\AF{4}$ offers a curvature gap comparable (and in some cases even smaller) to that of FRC. We remark this as an example of how further augmentations of AFRC do not always improve the effectiveness of the community detection. We believe this could be the same effect we observed for $\AF{5}$: in networks with low $\hat{p}/\hat{q}$, the inter- and intra-community edges are contained in a number of 4-cycles similar enough that the two curvature distributions overlap significantly.\\
~
\\
To summarize the experimental results on the curvature gap, we have seen that plain FRC is not able to detect community structure -- which was theoretically expected -- whereas both its augmentations and ORC do. The performance of all curvatures depends on how strongly pronounced the community structure is, but the characteristics of each network have different effects on the various curvature versions. In particular, we have found that ORC and $\AF{3}$ tend to outperform $\AF{4}$ in unipartite networks, which suggests $\AF{3}$ as the better (in terms of performance/cost) alternative for curvature based-community detection for this type- of network. On the other hand, $\AF{4}$ seems the most effective in the study of bipartite networks. In the following section we further discuss the performance in community detection of these curvature measures, both in terms of accuracy and computational cost.

\begin{table}[t]
\centering
\begin{tabular}{ |l|c|c|c|c| } 
\hline

\textbf{Network} & OR & FR & $\AF{3}$ & $\AF{4}$ \\
\hline
SBM ($10$, $5$, $0.7$, $0.05$) & 2.14 & 0.58 & 1.93 & 1.34\\ 
\hline 
SBM ($10$, $10$, $0.7$, $0.05$) & 3.30 & 0.35 & 3.54 & 1.88\\ 
\hline 
SBM ($10$, $15$, $0.7$, $0.05$) & 4.10 & 0.33 & 4.02 &
3.31\\ 
\hline 
SBM ($10$, $20$, $0.7$, $0.05$) & 4.92 & 0.23 & 4.61 & 4.33\\ 
\hline 
SBM ($10$, $20$, $0.5$, $0.05$) & 2.66 & 0.16 & 2.16 & 2.01\\ 
\hline 
SBM ($10$, $20$, $0.3$, $0.05$) & 0.97 & 0.07 & 0.79 & 0.44\\ 
\hline 
SBM ($10$, $20$, $0.1$, $0.05$) & 0.06 & 0.18 & 0.17 & 0.15\\

\hhline{|=|=|=|=|=|}
HBG ($50$, $0.5$, $0.05$) & 3.83 & 0.19 & 0.19 & 4.77\\
\hline
HBG ($50$, $0.5$, $0.1$) & 0.40 & 0.12 & 0.12 & 3.05\\
\hline
HBG ($50$, $0.5$, $0.15$) & 0.03 & 0.16 & 0.16 & 1.62\\
\hline
HBG ($50$, $0.4$, $0.1$) & 0.62 & 0.06 & 0.06 & 1.69\\
\hline
HBG ($50$, $0.3$, $0.1$) & 0.83 & 0.14 & 0.14 & 0.76\\
\hline
HBG ($50$, $0.2$, $0.1$) & 0.27 & 0.07 & 0.07 & 0.09\\ 
\hline 
\end{tabular}
\caption{Comparison of curvature gaps attained by Ollivier-Ricci curvature (OR), Forman-Ricci curvature (FR), and augmentations of the  Forman-Ricci curvature ($\AF{3}$ and $\AF{4}$) of edges in artificial networks with community structure.}
\label{table:3}
\end{table}

\begin{table}[t]
\centering
\begin{tabular}{ |l|c|c|c|c| } 
\hline
\textbf{Network} & OR & FR & $\AF{3}$ & $\AF{4}$ \\
\hline
Karate & 1.01 & 0.31 & 0.94 & 0.43\\ 
\hline
College Football & 2.95 & 0.34 & 2.87 & 2.91\\ 
\hline
School 1 & 1.16 & 0.74 & 1.78 & 0.57\\ 
\hline
School 2 & 1.04 & 0.85 & 1.72 & 0.65\\ 
\hline
\end{tabular}
\caption{Comparison of curvature gaps attained by Ollivier-Ricci curvature (OR), Forman-Ricci curvature (FR), and augmentations of the  Forman-Ricci curvature ($\AF{3}$ and $\AF{4}$) of edges in real-world networks.}
\label{table:4}
\end{table}

\section{Community-detection via sequential edge-deletion}\label{sec:Algorithm}

\noindent In this section, we further analyze the performance in community detection of our curvature measures, both in terms of accuracy and computational cost. Our algorithm proceeds by analogy with the Girvan-Newman method and with the ORC-based community detection method proposed in \cite{sia2019ollivier}, in that we incrementally delete edges to disconnect the network. We provide pseudo-code of the algorithm in Figure \ref{alg:cap}.\\

\noindent The process can be divided into the following steps, with only minimal differences between using $\AF{3}$ and $\AF{4}$: given a simple, unweighted, and undirected graph $G(V, E)$, (i) calculate $\AF{3}$ ($\AF{4}$) for all edges in the network, (ii) based on the $\AF{3}$ ($\AF{4}$) distribution, choose a curvature threshold $\Delta$, (iii) remove the most negative $\AF{3}$ (positive $\AF{4})$ edge, (iv) re-compute the edge $\AF{3}$ ($\AF{4}$) only for those edges affected by the removal, (v) check if all $\AF{3}$ ($\AF{4}$) are above (below) $\Delta$, otherwise repeat steps (iii) and (iv) until the condition is satisfied. It is clear that our algorithm must terminate in at most $|E|$ steps, as there will be no more edges to delete -- for a sensible threshold such as the low-point between the modes in the $\AF{3}$ ($\AF{4}$) histogram, it will in fact terminate much sooner. To choose the threshold $\Delta$, we fit a Gaussian mixture model with two modes to the curvature distribution. For normal distributions $\mathcal{N}_1(\mu_1, \sigma_1)$ and  $\mathcal{N}_2(\mu_2, \sigma_2)$ determined in this way, we then compute $\Delta=\frac{\sigma_2}{\sigma_1+\sigma_2}\mu_1 + \frac{\sigma_1}{\sigma_1+\sigma_2}\mu_2$, based on the minimum-error decision boundary between these Gaussians \cite[Ch. 2.6]{duda2012pattern}.\\



\begin{figure}[h!]
\hrule
\vspace{0.3em}
\begin{algorithmic}
\Require A graph $G(V, E)$ and a threshold $\Delta$
\Ensure{A list of tuples $(v, l)$, with vertices $v$ and labels $l$}
\State For all $e \in E$, compute $\AF{3}$ 
\While{there exists an edge with $\AF{3} < \Delta$ in $G$}
\If{ there is a unique edge $e_\text{min}$ with minimal $\AF{3}$}
\State remove $e_\text{min}$ from $G$
\Else
\State choose an edge $e_\text{min}$ from the edges with minimal
\State $\AF{3}$ uniformly at random and remove it from $G$
\EndIf
\State re-calculate $\AF{3}$ for all affected edges in $G$
\EndWhile
\State assign the same label $l$ to each vertex $v$ in a connected component of $G$\\
\Return a list of tuples $(v, l)$

\end{algorithmic}
\hrule
\vspace{0.3em}
\caption{AFRC community detection algorithm, here with $\AF{3}$.}\label{alg:cap}
\end{figure}

\noindent Different choices for $\Delta$ are possible and may improve the community detection accuracy. Empirically, we find that when there is less community structure, $\mathcal{N}_1$ and $\mathcal{N}_2$ can overlap significantly, which makes finding a good value for $\Delta$ harder and tends to degrade accuracy. Choosing $\Delta$ more aggressively, i.e. closer to the upper (lower) mode when using $\AF{3}$ ($\AF{4}$) seems to work well in these cases.\\

\begin{figure}[!t]
  \begin{subfigure}[t]{.48\textwidth}
    \centering
    \includegraphics[width=\linewidth]{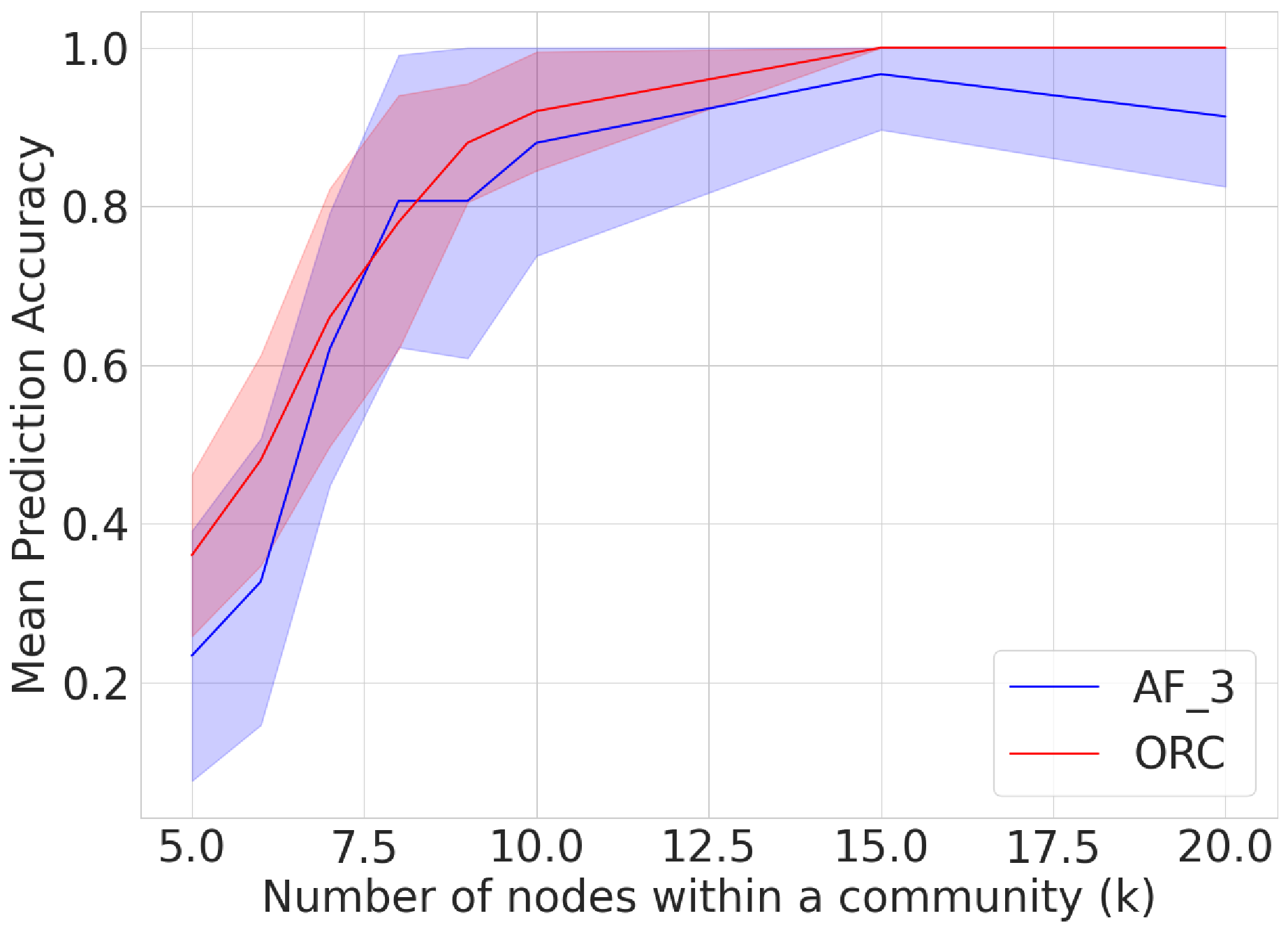}
    \caption{}
  \end{subfigure}
  \hfill
  \begin{subfigure}[t]{.48\textwidth}
    \centering
    \includegraphics[width=\linewidth]{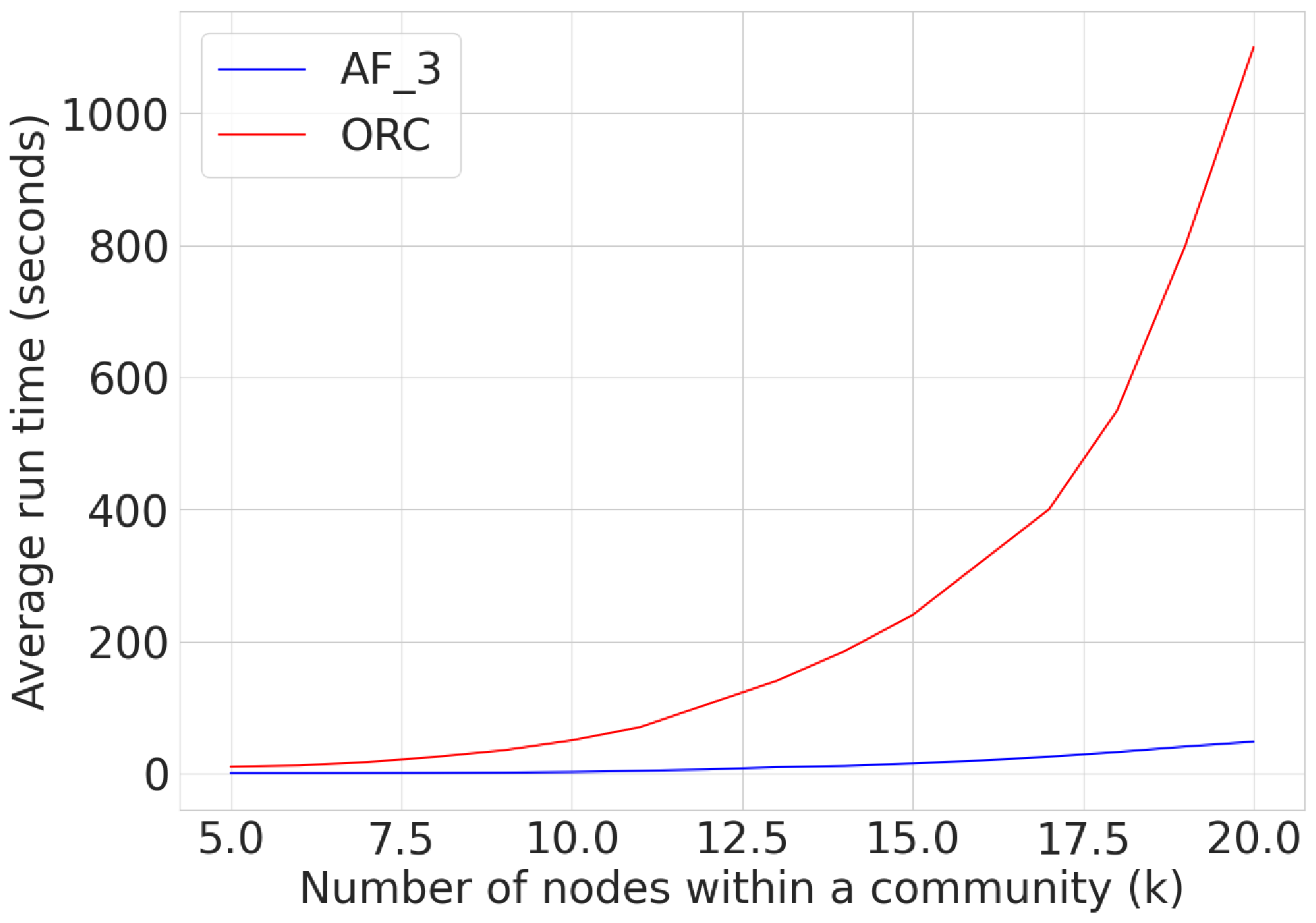}
    \caption{}
  \end{subfigure}

    \medskip

    \begin{subfigure}[t]{.48\textwidth}
    \centering
    \includegraphics[width=\linewidth]{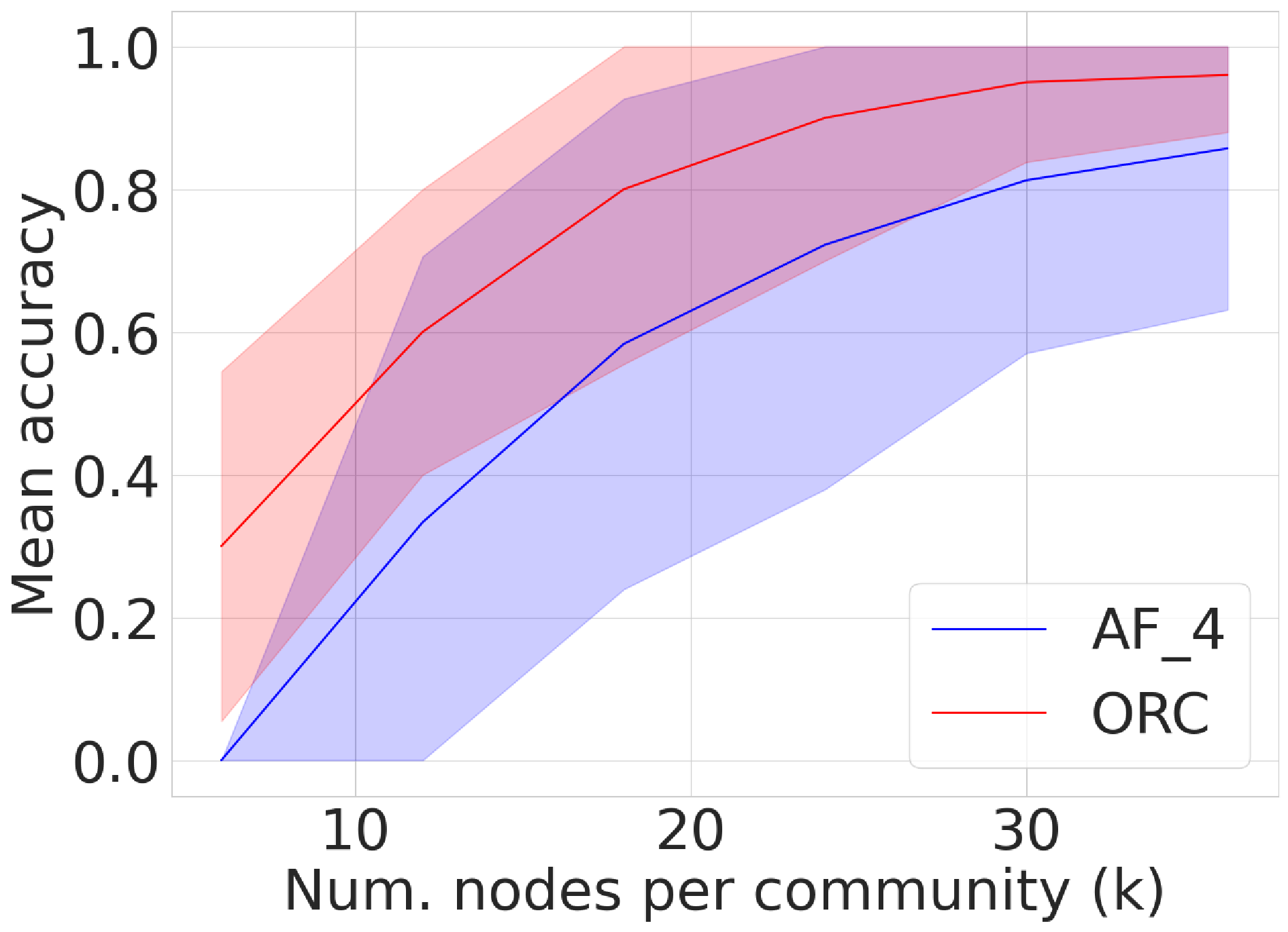}
    \caption{}
  \end{subfigure}
  \hfill
  \begin{subfigure}[t]{.48\textwidth}
    \centering
    \includegraphics[width=\linewidth]{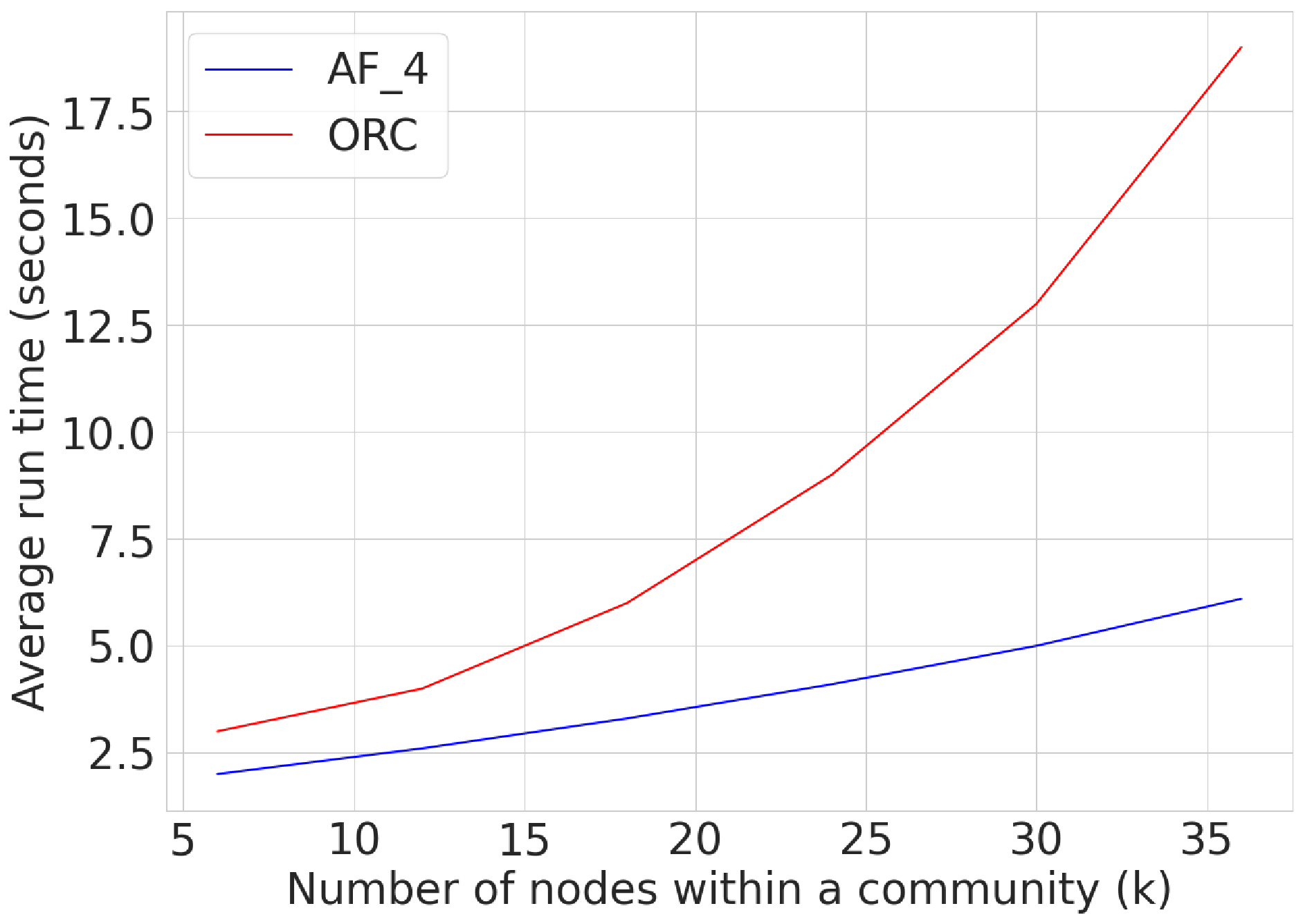}
    \caption{}
  \end{subfigure}
  
    \caption{Mean accuracy (percentage of correctly identified communities) and standard deviation (shaded areas) of ORC-based (red) and $\AF{3}$-based community-detection (blue) in \textbf{(a)}, and the time required for the algorithm to terminate (in seconds) in \textbf{(b)} for an SBM with parameters $(10, k, 0.7, 0.05)$. \textbf{(c)} and \textbf{(d)} depict the same using the $\AF{4}$ for a hierarchical bipartite graph with parameters $(k + k, 0.5, 0.05)$}\label{fig:sbm_accuracies}
\end{figure}

\noindent In general, it is possible to use a priori knowledge of a network to choose an appropriate augmentation of the Forman curvature to detect communities in a computationally efficient way, but with accuracy similar to that achieved using the ORC. For example, as proposed in Section \ref{sec:CurvGaps}, for unipartite, SBM-like networks with relatively clear community structure (i.e. large ratio $\hat{p} / \hat{q}$), we can use $\AF{3}$. As Figure \ref{fig:sbm_accuracies} shows, the accuracy of this method is comparable to that of the ORC-based method (panel \textbf{(a)}), while being orders of magnitude faster to run (panel \textbf{(b)}). \\

\noindent For networks which we suspect to be (approximately) bipartite, we can include the 4-cycles, which makes the algorithm more costly to run. However, as can be seen in Figure \ref{fig:sbm_accuracy}, panels \textbf{(c)} and \textbf{(d)}, this is still significantly faster than using the ORC, while suffering only minor drops in accuracy. We also note that the variance in our HBG results is quite large, which we believe is because we only have two communities. Nonetheless, for a community size of $\sim30$, $\AF{4}(e)$ and ORC are on average competitive. The computational efficiency of our algorithm follows from the fact that for an edge $e = (u,v)$ with $d_\text{max} = \text{max} \{\text{deg}(u), \text{deg}(v) \}$, computing $\AF{3}(e)$ scales as $O(d_\text{max})$, $\AF{4}(e)$ as $O(d_\text{max}^2)$, and $\mathcal{O}(e)$ as $O(d_\text{max}^3)$. The number of edges deleted by the algorithm is approximately equal for all three curvatures, so this difference in computational complexities explains the increased speed observed when using $\AF{3}(e)$ or $\AF{4}(e)$.\\

\noindent In Appendix Figures \ref{fig:football_accuracy} and \ref{fig:karate_accuracy}, we present additional results on real-world networks and also compare our algorithm against two well-established community detection algorithms: the Girvan-Newman Algorithm and the Louvain Method. On the American College Football network, the ORC-based algorithm correctly identifies 6 out of 12 communities, while the $\AF{3}$- and $\AF{4}$-based algorithms each identify 5 communities correctly. Girvan-Newman performs best and correctly detects 7 communities, while Louvain only identifies 3 ground-truth communities. On the Karate Club, all five algorithms fail to correctly identify the two underlying communities, but the Figure \ref{alg:cap} algorithm with $\AF{3}$ and Girvan-Newman arguably perform best as both only mislabel two nodes.\\
    
\noindent We also conduct some additional experiments to compare the performance of our algorithm with the three curvatures around the community detection threshold, and find that for networks with few or no cycles within communities, all three curvatures yield small curvature gaps, which makes community detection difficult or even impossible. We quantify this using the notion of tree-SBMs, i.e. SBMs where the communities are not Erd\H{o}s-R\'{e}nyi random graphs, but randomly chosen trees with a small probability $p$ of having additional edges within communities. For $p$ close to zero, only the ORC picks up on the emerging community structure. As we increase $p$, and our tree-SBM becomes more similar to a standard SBM, the augmentations of the Forman-Ricci curvature begin to also yield increasingly large curvature gaps --- see Appendix Table \ref{table:5}.\\


\noindent Similarly, we find that for standard SBMs where the probability $p$ of an edge within a community is only marginally larger than the probability $q$ of an edge between communities, the ORC yields larger curvature gaps and more accurate community detection results than our augmentations of the Forman-Ricci curvature. As we increase $p$, all three curvatures begin to yield similarly good accuracies -- see a comparison of ORC vs $\af{3}$ in Appendix Figure \ref{fig:sbm_p_in}. This underlines the idea that for networks with a particularly difficult community structure, especially in the very sparse regime, the ORC is generally the best, and sometimes only, curvature for community detection. In most cases however, the augmentations of the Forman-Ricci curvature give similarly good results, which together with their computational efficiency makes them preferable to the ORC.

\section{Conclusion}
\label{sec:Conclusion}
Discrete notions of curvature have recently become the subject of interest in network science, both for theoretical investigations and in applications such as community detection. Much of this interest has been focused on the Ollivier-Ricci Curvature: to the best of our knowledge, all previous curvature-based community-detection algorithms have used the ORC. While this yields state-of-the-art results for many networks, it is also computationally expensive and quickly becomes inapplicable for larger graphs.\\

\noindent In this work, we chose a different approach and provided a detailed study of augmentations of the Ricci-curvature discretization proposed by Forman (AFRC). We discussed the relation of these augmentations to the FRC and ORC both empirically and theoretically -- in particular, we provided evidence that in many artificial and real-world networks, Augmented Forman-Ricci curvature can be used to attain similar or even better community detection results than the traditionally used Ollivier-Ricci curvature. Established Ollivier-Ricci curvature-based algorithms for community detection are therefore not without well-performing and computationally cheaper alternatives.\\

\noindent The advantages of AFRC over ORC in community detection suggest that the former could represent important improvements in other areas where different discrete curvature measures have already found an application -- for instance network embedding or graph learning. Moreover, the fact that AFRC and ORC behave differently in the studied networks indicates that the former is not just a computational improvement over the latter, but it also reflects different geometric properties of networks. The augmentations proposed in this paper may thus provide further insight into curvature-related problems, and create new possibilities in contexts where other curvature measures are ineffective. For instance, we found the augmentation based on 4-cycles to significantly outperform ORC in community detection in bipartite networks, which opens opportunities for the design of curvature-based methods for this kind of networks.

\newpage

\appendix
\section{Studied networks}\label{appendix: studied networks}

\noindent\textit{\textbf{Model networks}}
\begin{itemize}
    \item \textbf{Erd\H{o}s-Rényi graph:} ER$(m,p)$. Network of $m$ nodes where each edge has an independent probability $p$ of being present. The bipartite version BG$(m,p)$ consists on two groups of $m$ nodes, with the inter-group edges having a probability $p$ of appearing.

    \item \textbf{Stochastic Block Model:} SBM$(l,k,p,q)$. Network of $l\cdot k$ nodes divided into $l$ communities of size $k$. Each intra-group (resp. inter-group) edge has an independent probability $p$ (resp. $q$) of being present.

    \item \textbf{Hierarchical Bipartite Graph:} HBG$(m,p,q)$. Bipartite network of $4m$ nodes divided into $2$ communities of size $2m$, where each community contains $m$ nodes of each partition.  Each intra-group (resp. inter-group) edge has an independent probability $p$ (resp. $q$) of being present.\\
\end{itemize}\vspace{1em}

\noindent\textit{\textbf{Real-world networks}}

\begin{table}[h]
\centering
\begin{tabular}{|l|c|c|c|c|c|c|c|}
\hline
\textbf{Network}                                               & Nodes       & Edges       & Triangles & Squares   & $\hat{p} / \hat{q}$ \\ \hline
Dolphin \cite{lusseau2003bottlenose}          & 62          & 159         & 95        & 278                        & -          \\ \hline
Power Grid \cite{watts1998collective}         & 4941        & 6594        & 651       & 979                        & -          \\ \hline
Word Adjacency \cite{newman2006finding}       & 112         & 425         & 284       & 2579                     & -          \\ \hhline{|=|=|=|=|=|=|}
Southern Women \cite{davis1941deep}           & 32          & 89          & 0         & 341                       & -          \\ \hline
Corporate Interlocks \cite{scott1980scottish} & 244 & 358 & 0         & 245                 & -          \\ \hhline{|=|=|=|=|=|=|}
Karate Club \cite{zachary1977information}     & 34          & 78          & 45        & 154                      & $6.57$    \\ \hline
College Football \cite{girvan2002community}   & 115         & 613         & 810       & 3915                     & $19.51$    \\ \hline
School 1 \cite{stehle2011high}                & 236         & 5899        & 48401     & 1678496                  & $5.89$    \\ \hline
School 2 \cite{stehle2011high}                & 238         & 5539        & 46687     & 1560642                   & $6.76$   \\ \hline
\end{tabular}
\caption{Summary of the studied real-world networks. The three sections correspond respectively to networks without additional structure, bipartite networks, and networks with community metadata.}
\label{table:realnetworks}
\end{table}

\section{Example: computing $\af(e)$ in a small graph}\label{appendix: detailed example}
\begin{figure}[h!]
    \centering
    \includegraphics[width=0.4\textwidth]{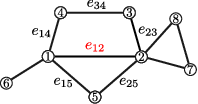}
    \caption{}
    \label{fig: example detailed}
\end{figure}

We revisit the example from Section \ref{sec:Definitions} and calculate the FRC and AFRC of edge $e=e_{12}=(1,2)$ in Figure \ref{fig: example detailed}.\\

\noindent For the Forman-Ricci curvature \eqref{eq:FRC}, we only have to count the degrees of the end points of edge $e_{12}$. These are $\deg(1)=4$ and $\deg(2)=5$ and thus we obtain 
$$
\mathcal{F}(e) = 4 - \deg(1)-\deg(2) = -5.
$$
For the augmented Forman-Ricci curvature, we need to compute the count matrix $\Gamma$. Since we focus on edge $e_{12}$ it suffices to calculate the ${e_{12}}^{\text{th}}$ row and column of $\Gamma$. We fix the vertex order $1<2<\dots<7<8$ which induces the order $e_{12}<e_{14}<e_{15}<e_{23}<e_{25}<e_{34}$ on the edges. Since $e_{12}$ is contained in two cycles, $\gamma_1=1234$ and $\gamma_2=125$, we have
$$
\Gamma_{e_{12}e_{12}}=2
$$
Cycle $\gamma_1=1234$ traverses $e_{12},e_{23},e_{34}$ from small to large, and $e_{14}$ from large to small. Applying definition \eqref{eq:definition counts}, we find 
$$
\begin{cases}
\Gamma_{e_{12}e_{23}}=1 \text{~~(since $e_{12},e_{23}$ are aligned in $\gamma_1$, and $e_{12}<e_{23}$)}\\
\Gamma_{e_{12}e_{34}}=1 \text{~~(since $e_{12},e_{34}$ are aligned $\gamma_1$, and $e_{12}<e_{34}$)}\\
\Gamma_{e_{14}e_{12}}=1 \text{~~(since $e_{12},e_{14}$ are not aligned in $\gamma_1$, and $e_{14}>e_{12}$)}
\end{cases}
$$
Cycle $\gamma_2=125$ traverses $e_{12}, e_{25}$ from small to large and $e_{15}$ from large to small. Applying definition \eqref{eq:definition counts}, we find
$$
\begin{cases}
\Gamma_{e_{12}e_{25}}=1 \text{~~(since $e_{12},e_{25}$ are aligned in $\gamma_2$, and $e_{12}<e_{25}$)}\\
\Gamma_{e_{15}e_{12}}=1 \text{~~(since $e_{12},e_{15}$ are not aligned in $\gamma_2$ and $e_{15}>e_{12}$)}
\end{cases}
$$
All other entries in the ${e_{12}}^{\text{th}}$ row and column of $\Gamma$ are zero since they correspond to edges that do not share a cycle with $e_{12}$ (e.g. $\Gamma_{e_{12}e_{16}}=0$). Finally, the neighbors of $e_{12}$ are $\lbrace e_{14},e_{15},e_{16},e_{23},e_{25},e_{28},e_{27}\rbrace \sim e_{12}$ and the non-neighbors are $\lbrace e_{34},e_{78}\rbrace \not\sim e_{12}$. The AFRC is calculated from expression \eqref{eq:AFRC} as
\begin{align*}
\af(e_{12}) &= 2 + \Gamma_{e_{12}e_{12}} - \sum_{e'\sim e_{12}}\vert \Gamma_{e_{12}e'} + \Gamma_{e'e_{12}}-1\vert - \sum_{e'\not\sim e_{12}} \vert \Gamma_{e_{12}e'} - \Gamma_{e'e_{12}}\vert
\\
&= 2 + 2 - \left(4\cdot\vert 1-1\vert + 3\cdot \vert0-1\vert\right) - (\vert 1\vert + \vert 0\vert)
\\
&= 0.
\end{align*}
The code to calculate the AFRC in a graph will be made available along with the publication of the paper. This includes a calculation of $\mathcal{AF}(e_{12})$ in the example above using the codebase.
\clearpage

\section{Additional Numerical Experiments}
\label{appendix: additional numerical experiments}

\begin{figure}[h!]
  \begin{subfigure}[t]{.24\textwidth}
    \centering
    \includegraphics[width=\linewidth]{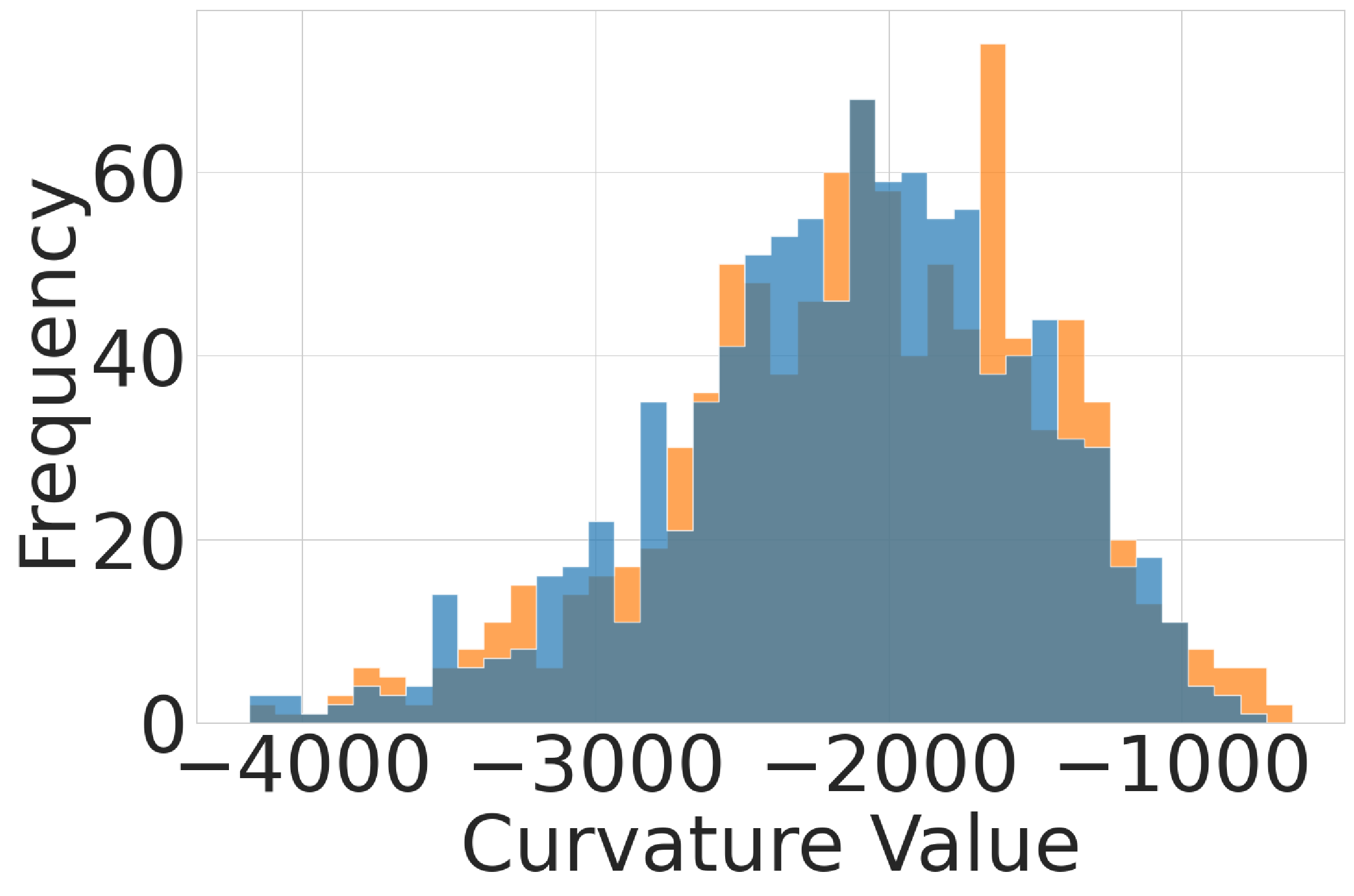}
  \end{subfigure}
  \hfill
  \begin{subfigure}[t]{.24\textwidth}
    \centering
    \includegraphics[width=\linewidth]{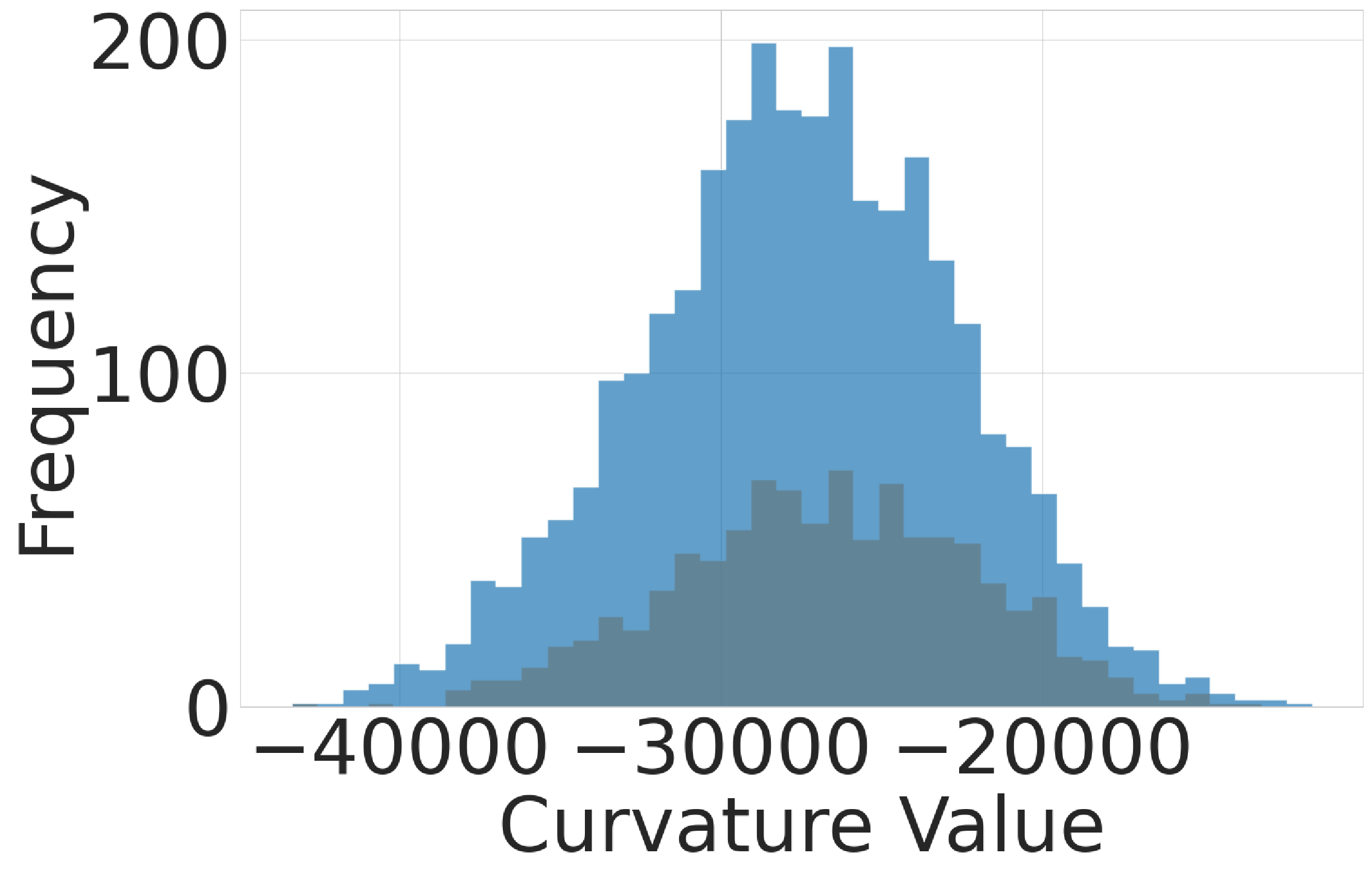}
  \end{subfigure}
    \begin{subfigure}[t]{.24\textwidth}
    \centering
    \includegraphics[width=\linewidth]{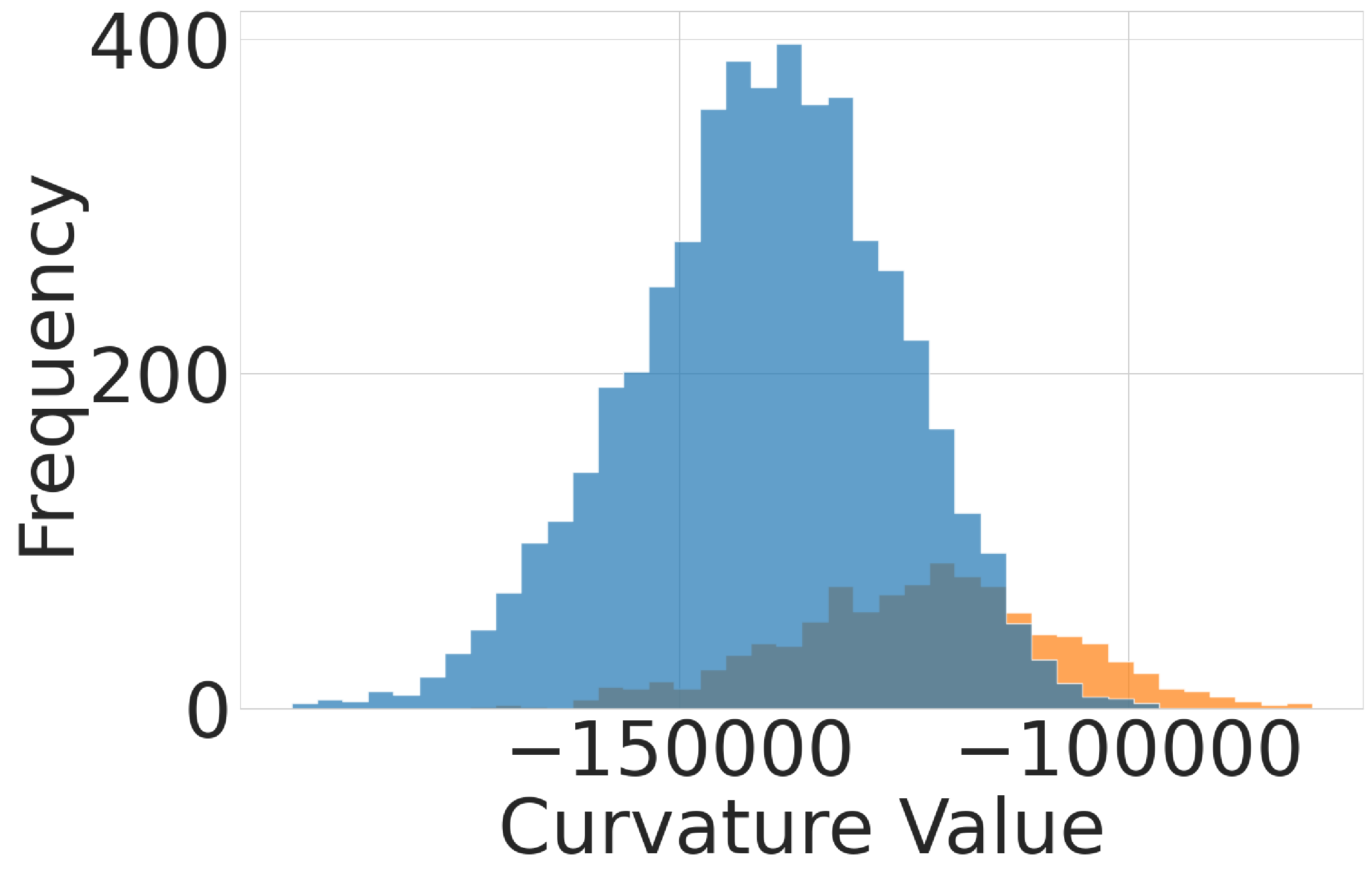}
  \end{subfigure}
  \hfill
  \begin{subfigure}[t]{.24\textwidth}
    \centering
    \includegraphics[width=\linewidth]{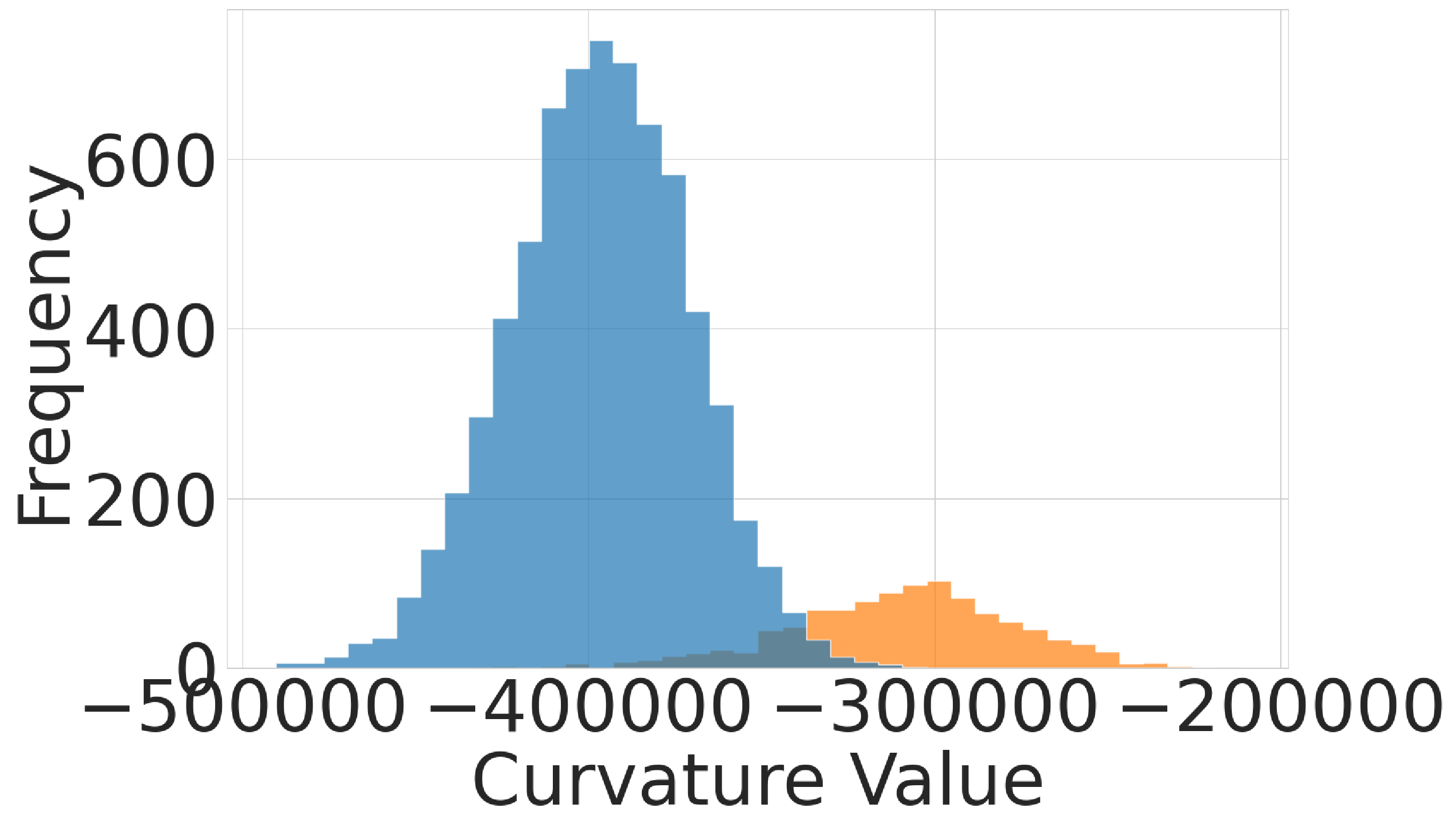}
  \end{subfigure}

  \medskip

  \begin{subfigure}[t]{.24\textwidth}
    \centering
    \includegraphics[width=\linewidth]{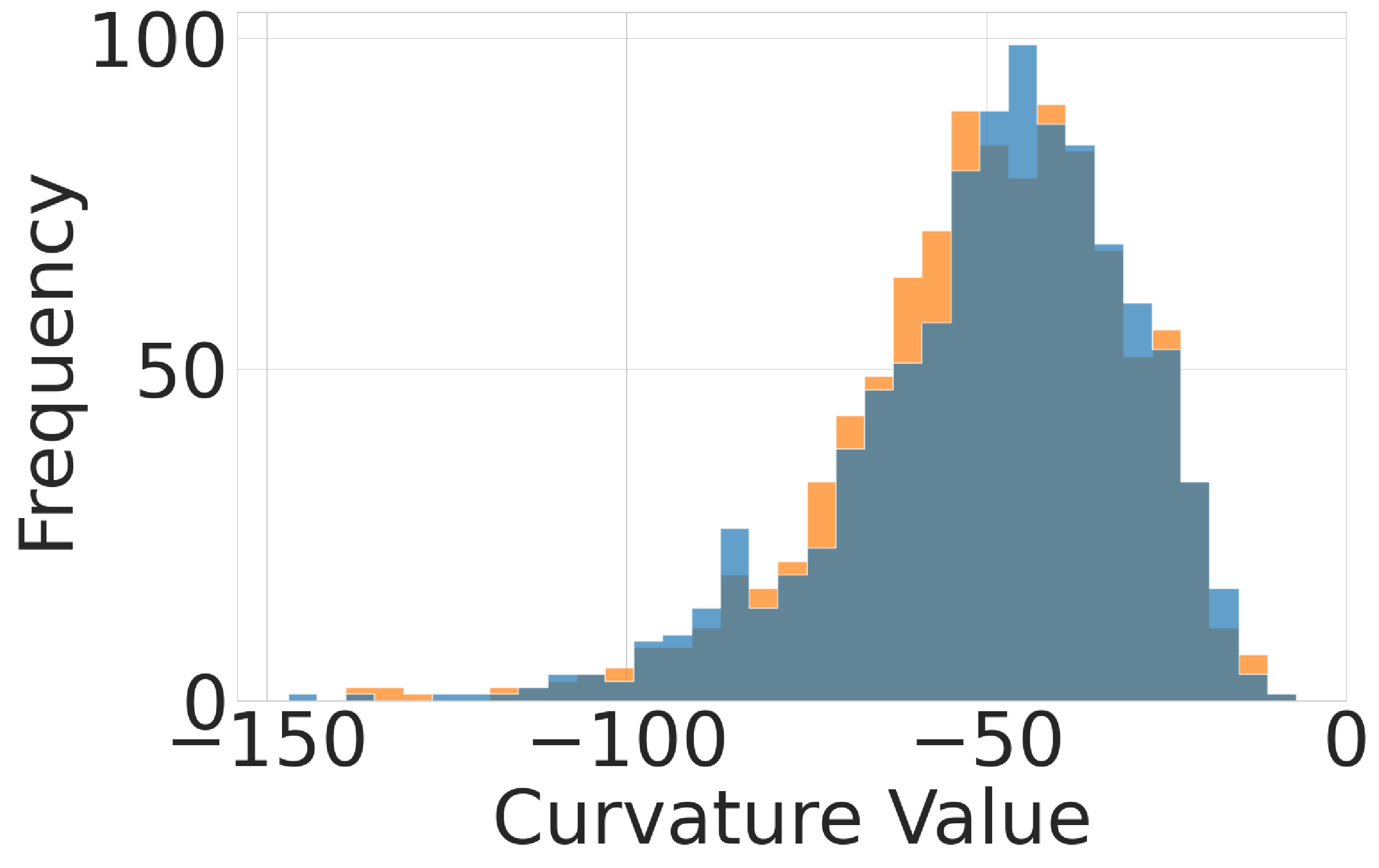}
  \end{subfigure}
  \hfill
  \begin{subfigure}[t]{.24\textwidth}
    \centering
    \includegraphics[width=\linewidth]{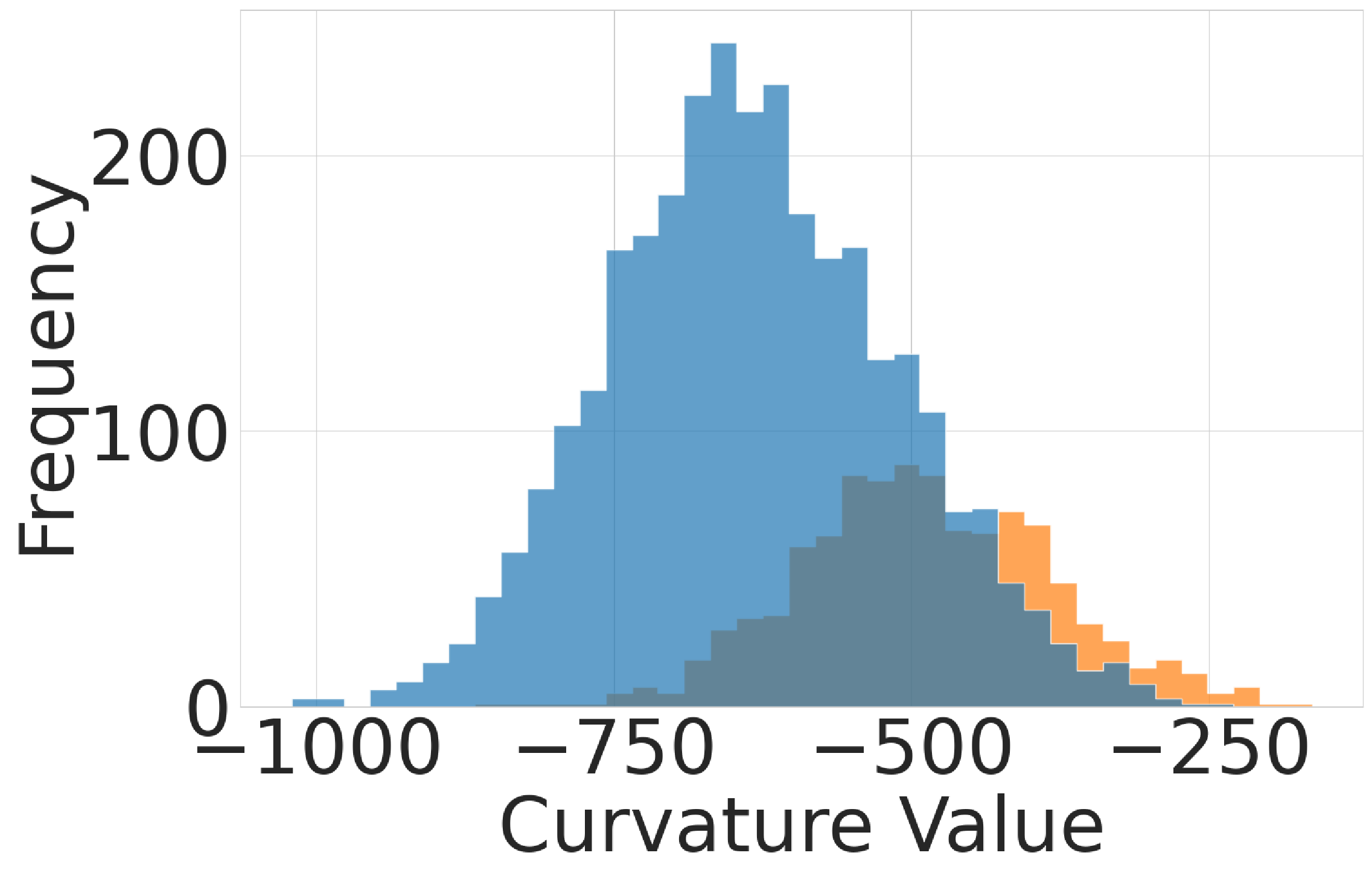}
  \end{subfigure}
    \begin{subfigure}[t]{.24\textwidth}
    \centering
    \includegraphics[width=\linewidth]{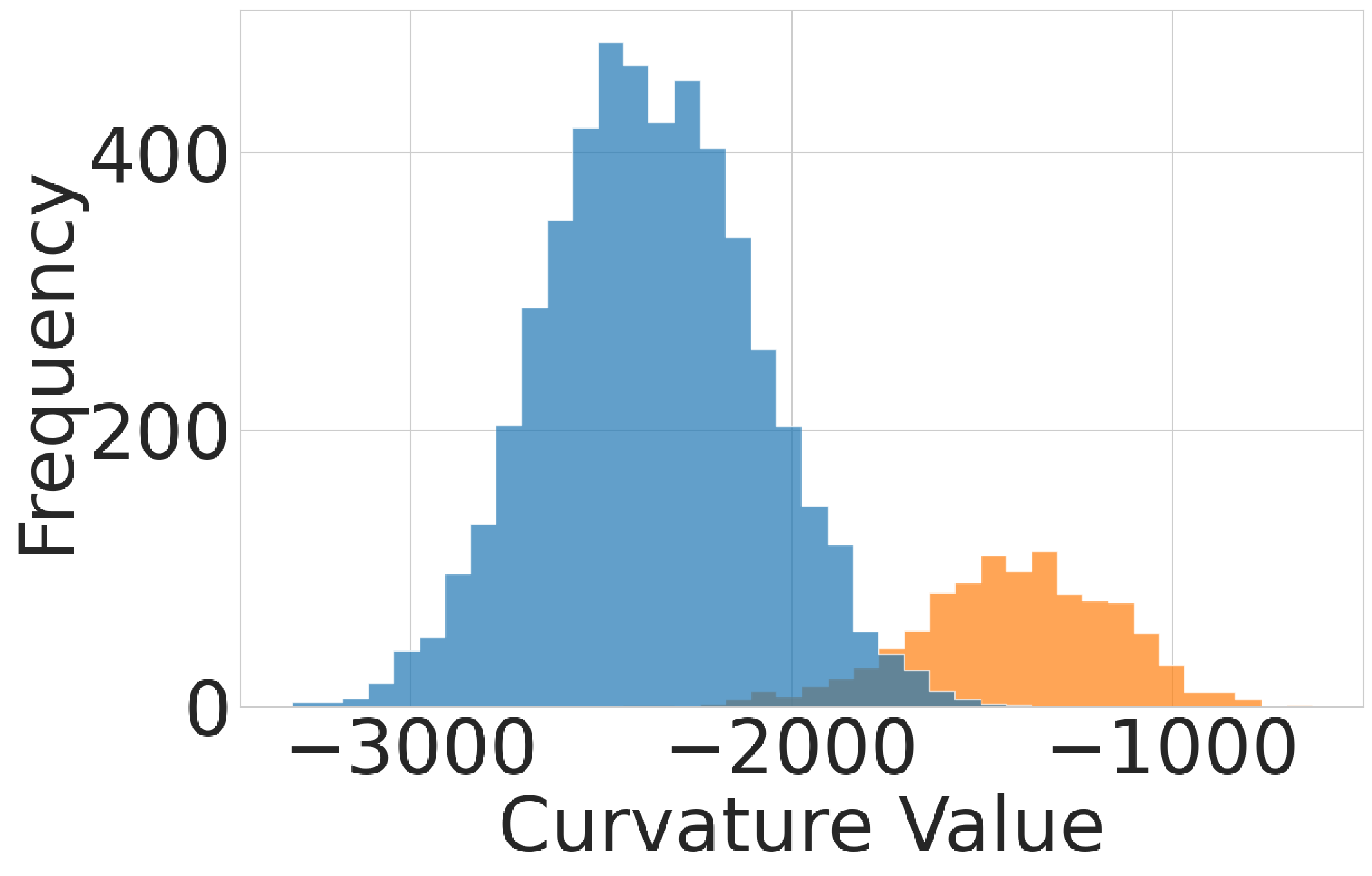}
  \end{subfigure}
  \hfill
  \begin{subfigure}[t]{.24\textwidth}
    \centering
    \includegraphics[width=\linewidth]{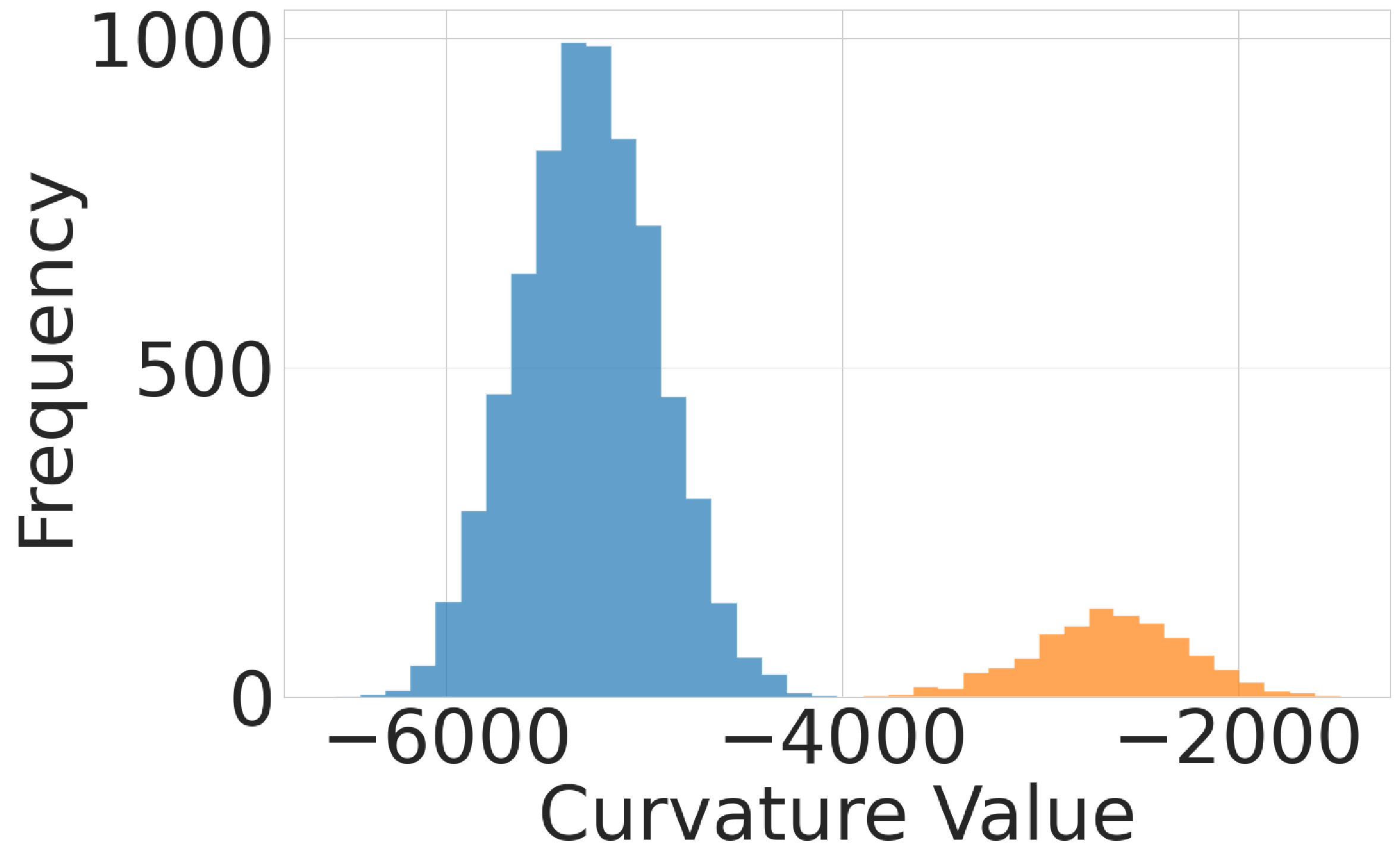}
  \end{subfigure}

  \medskip

  \begin{subfigure}[t]{.24\textwidth}
    \centering
    \includegraphics[width=\linewidth]{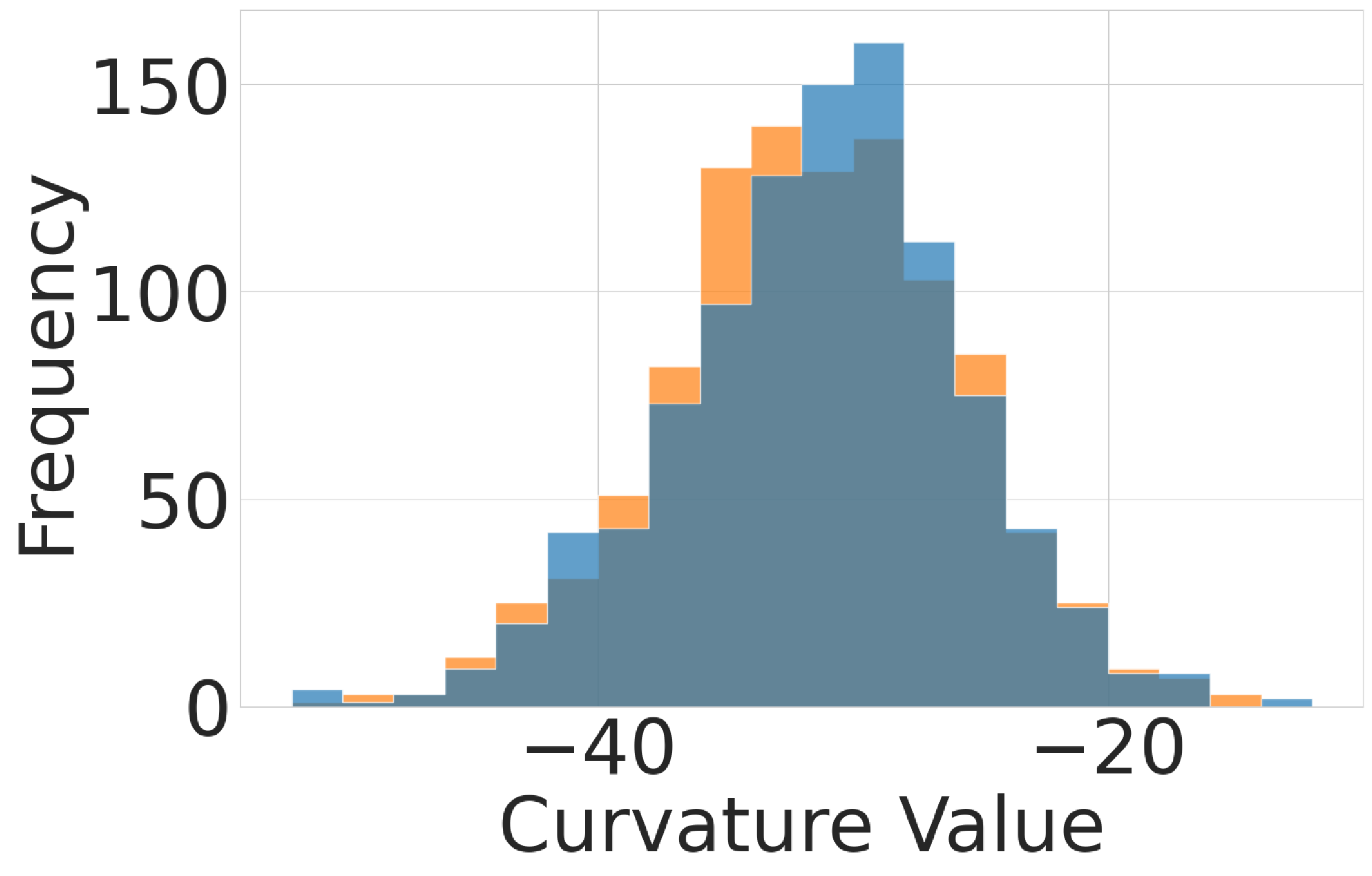}
  \end{subfigure}
  \hfill
  \begin{subfigure}[t]{.24\textwidth}
    \centering
    \includegraphics[width=\linewidth]{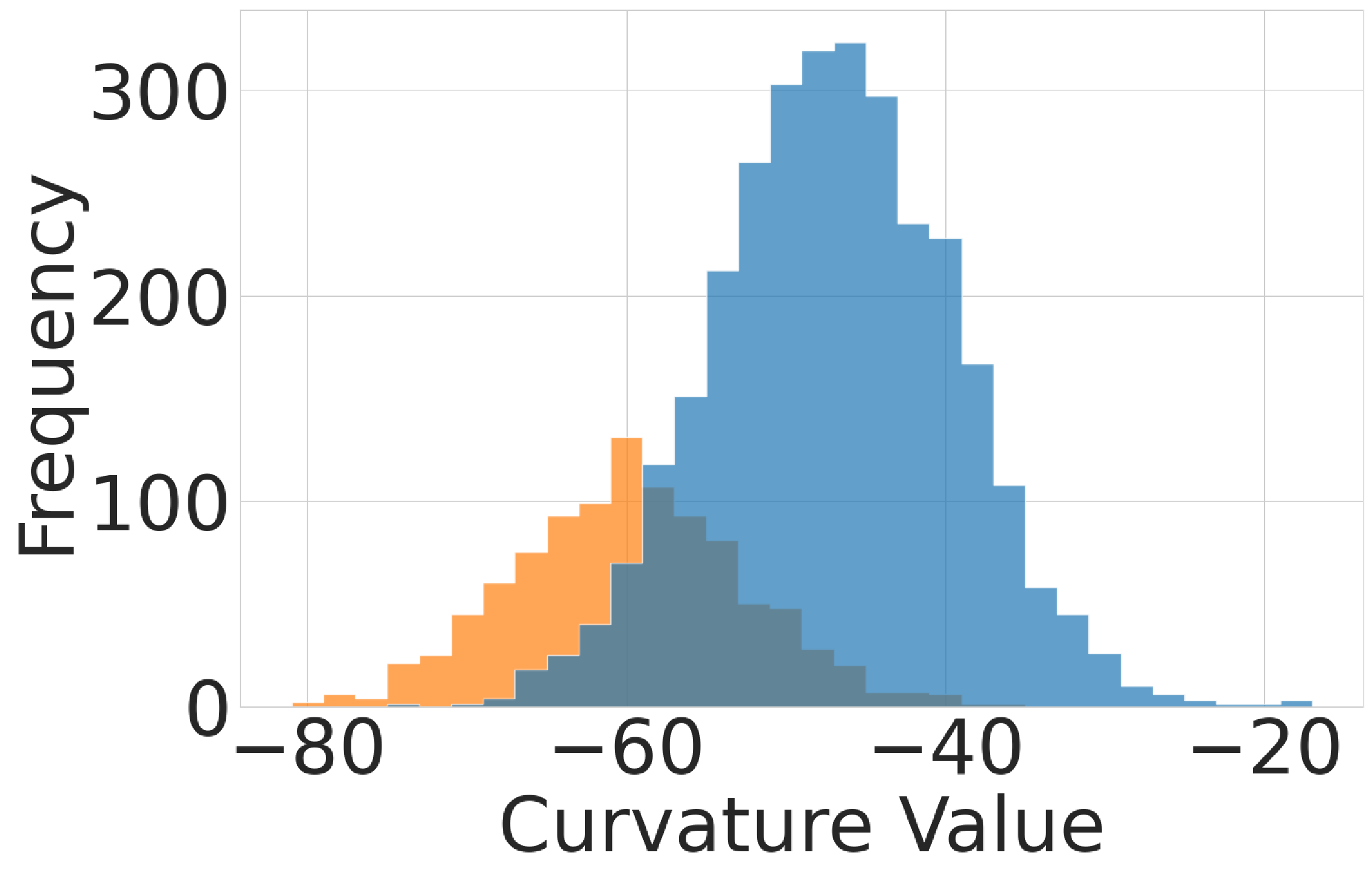}
  \end{subfigure}
    \begin{subfigure}[t]{.24\textwidth}
    \centering
    \includegraphics[width=\linewidth]{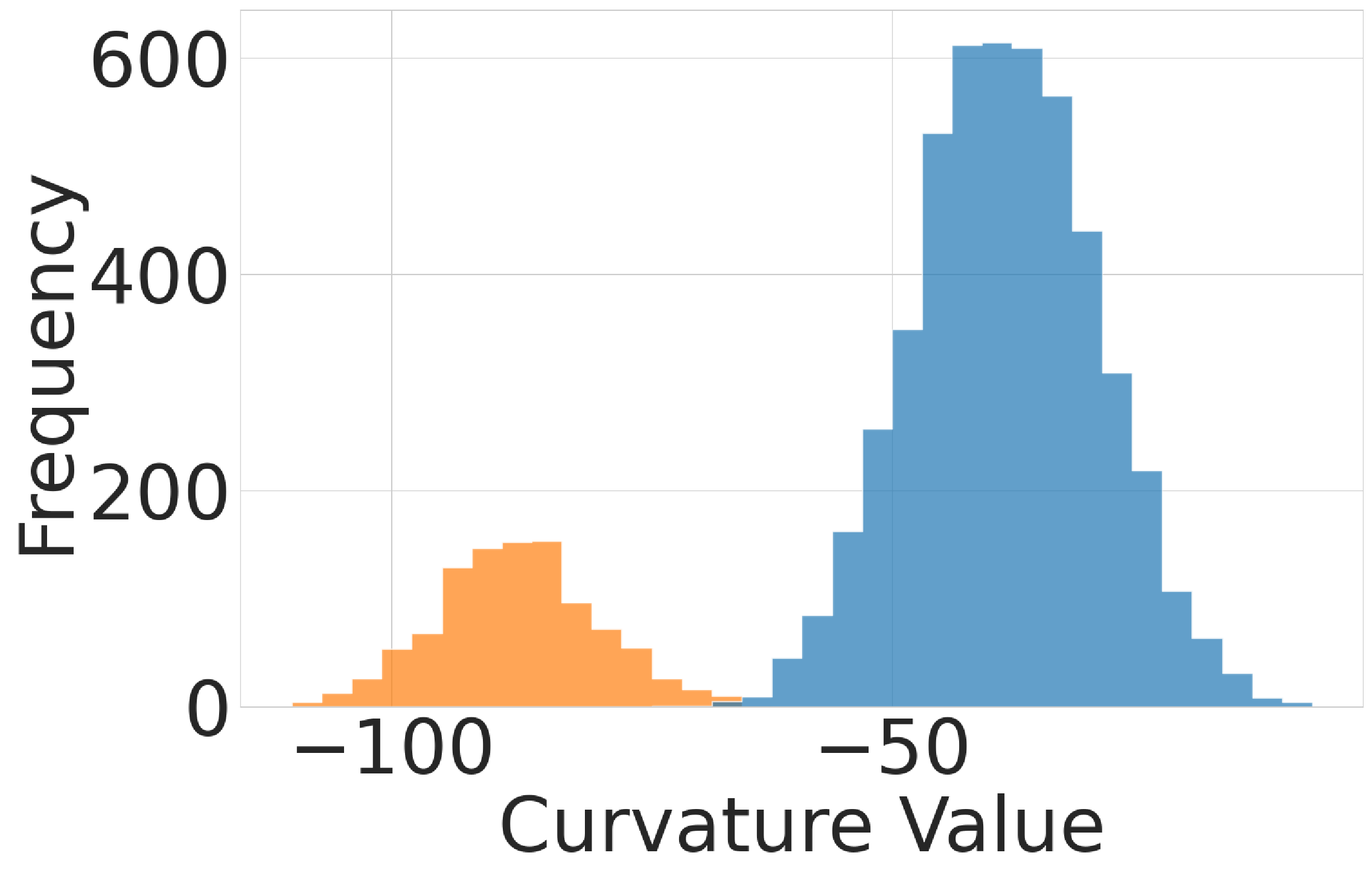}
  \end{subfigure}
  \hfill
  \begin{subfigure}[t]{.24\textwidth}
    \centering
    \includegraphics[width=\linewidth]{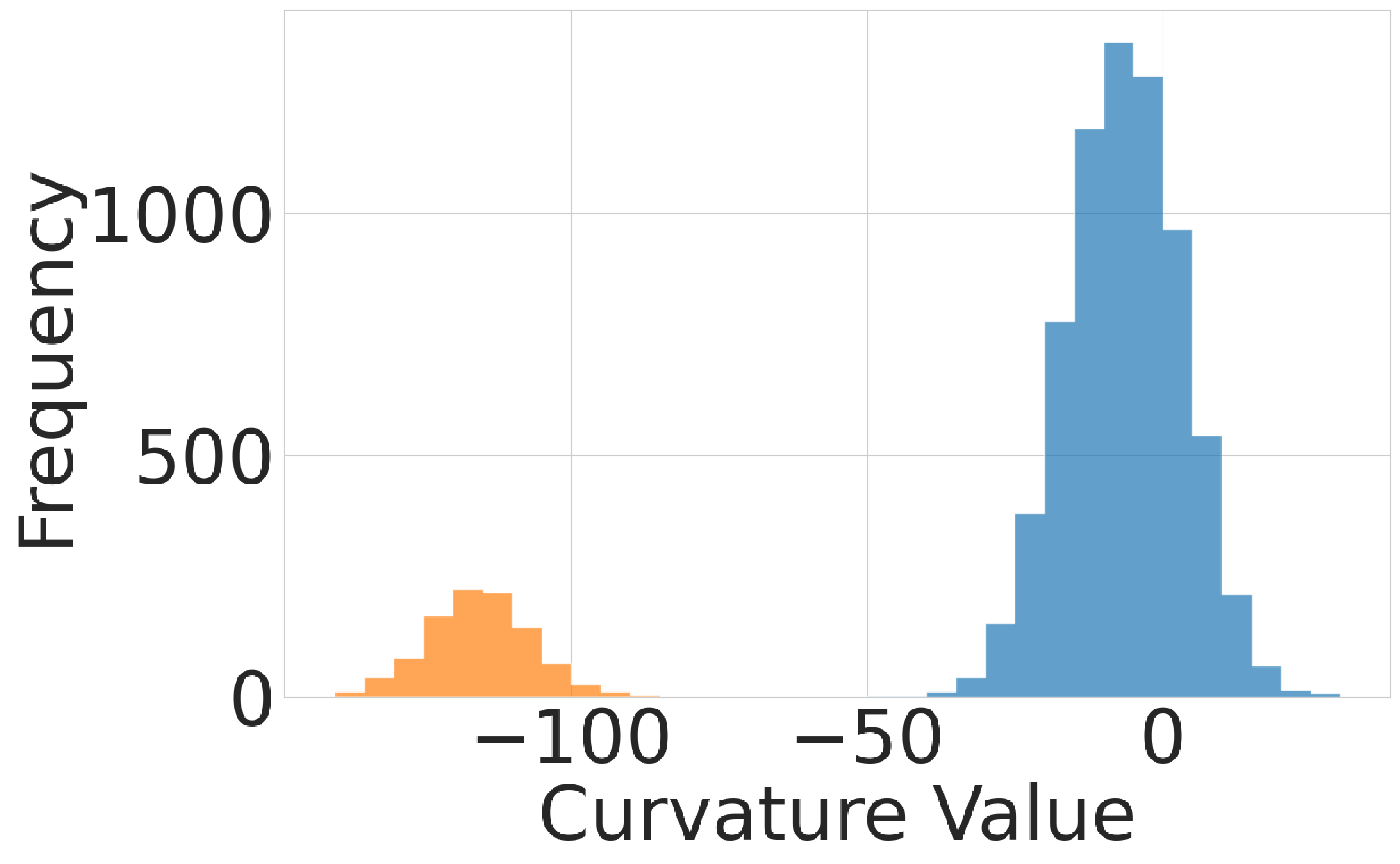}
  \end{subfigure}

  \medskip

  \begin{subfigure}[t]{.24\textwidth}
    \centering
    \includegraphics[width=\linewidth]{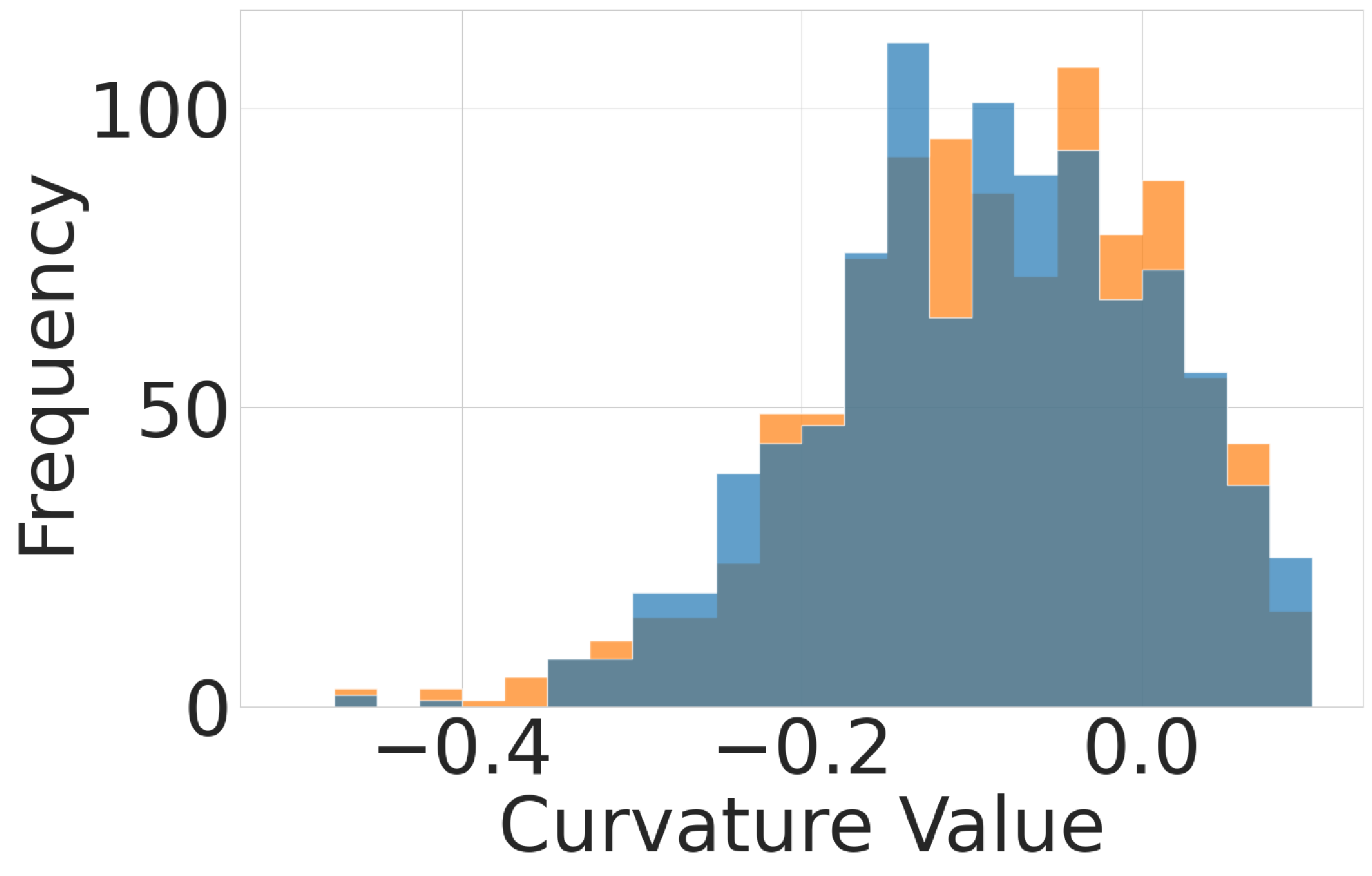}
  \end{subfigure}
  \hfill
  \begin{subfigure}[t]{.24\textwidth}
    \centering
    \includegraphics[width=\linewidth]{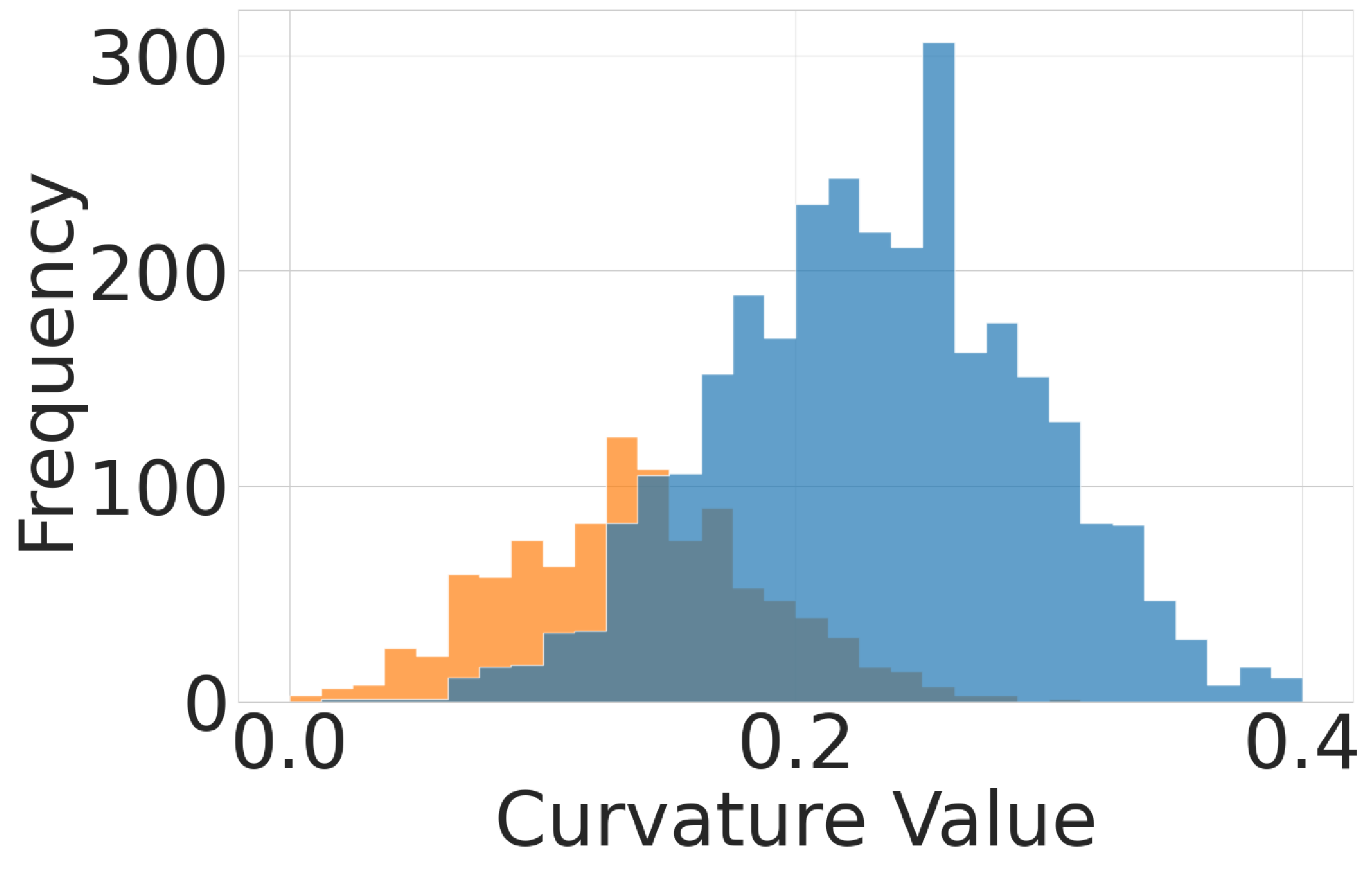}
  \end{subfigure}
    \begin{subfigure}[t]{.24\textwidth}
    \centering
    \includegraphics[width=\linewidth]{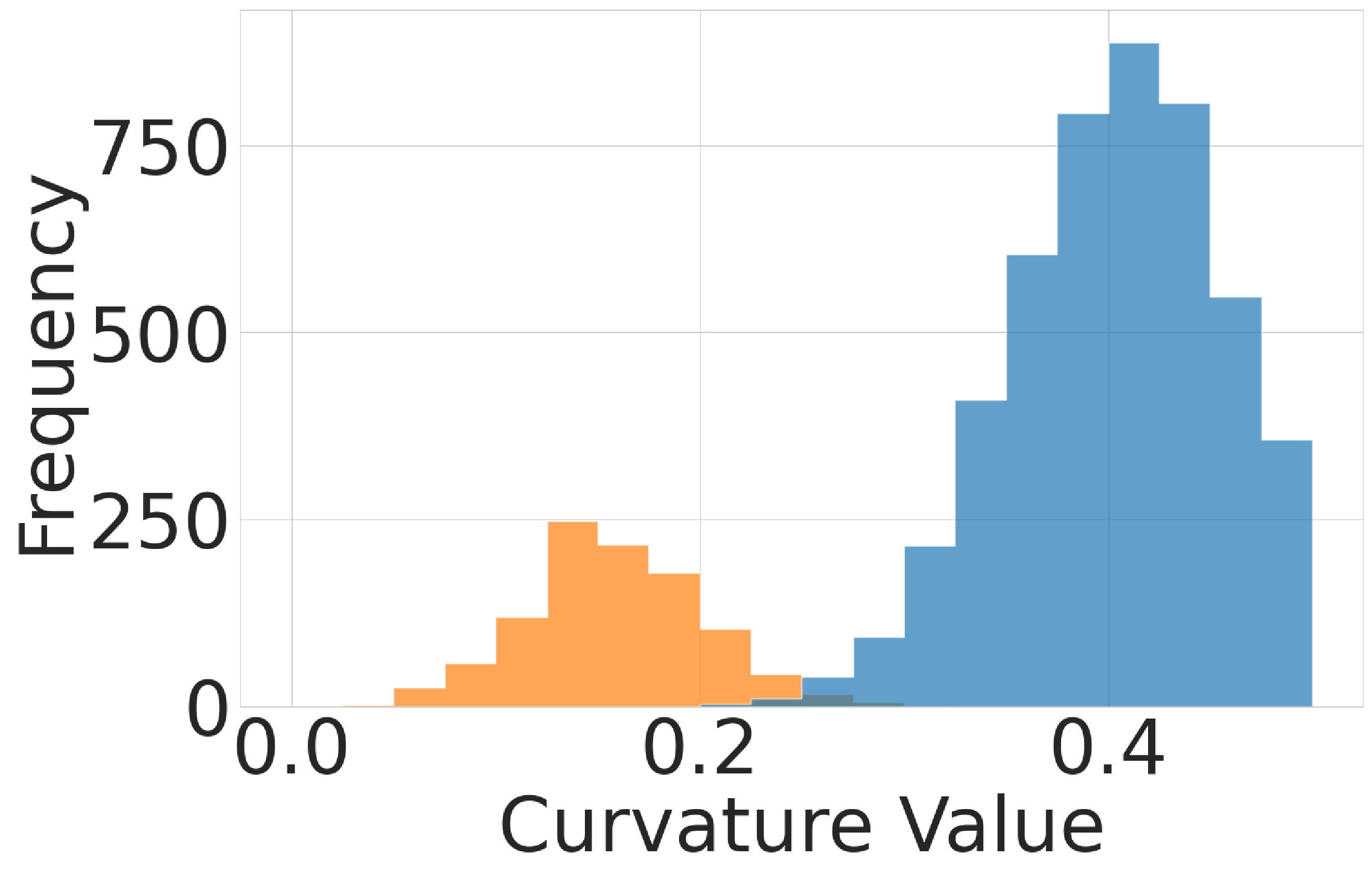}
  \end{subfigure}
  \hfill
  \begin{subfigure}[t]{.24\textwidth}
    \centering
    \includegraphics[width=\linewidth]{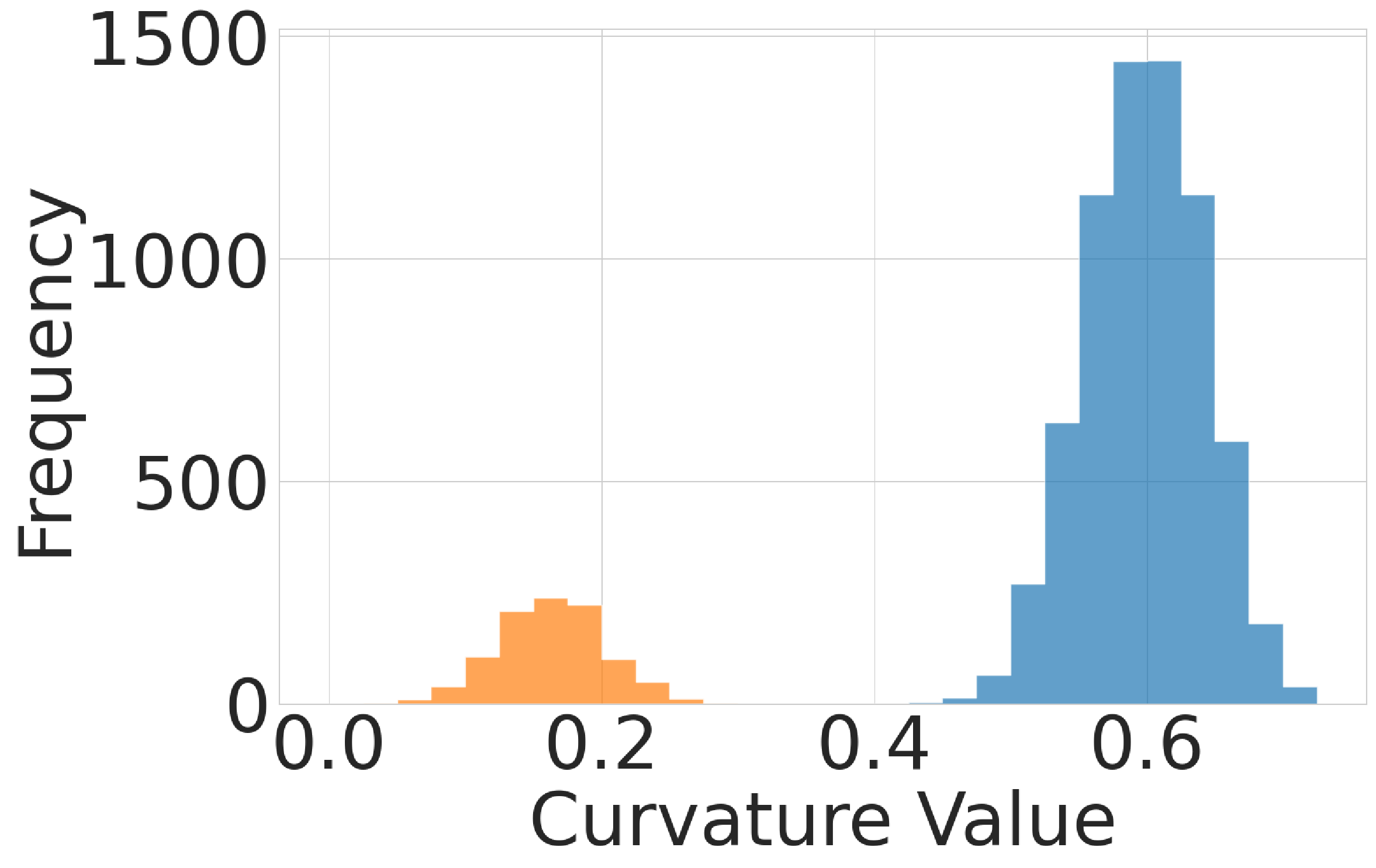}
  \end{subfigure}
  
  \caption{Top to bottom: Distributions of $\AF{5}$, $\AF{4}$, $\AF{3}$, and ORC for SBMs($2, 100, p, 0.1$) for $p \in \{0.1, 0.3, 0.5, 0.7\}$ (left to right). Edges within (between) communities in blue (orange).}
  \label{fig:AppC_1}
\end{figure}

\begin{figure}
  \begin{subfigure}[t]{.24\textwidth}
    \centering
    \includegraphics[width=\linewidth]{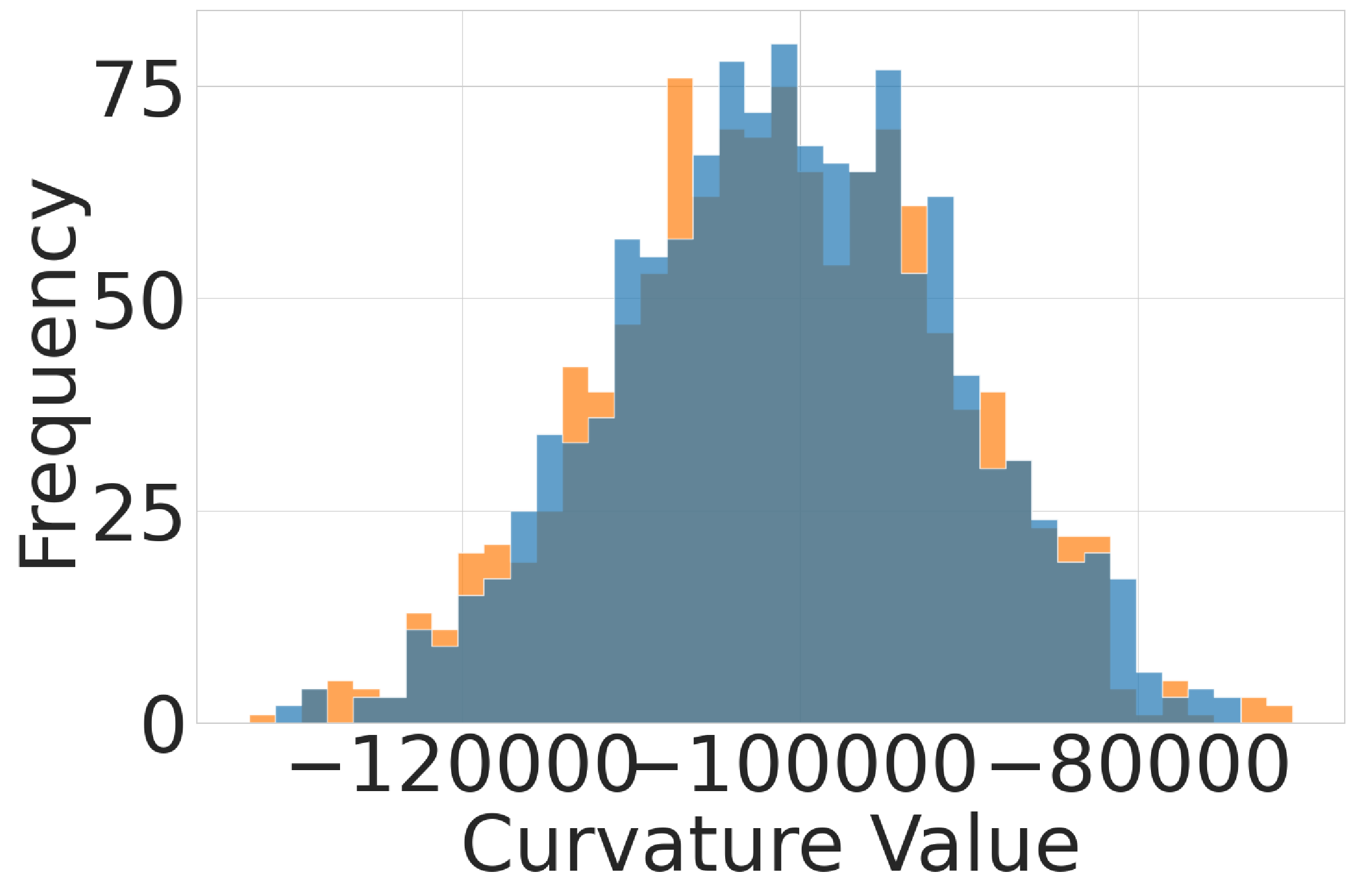}
  \end{subfigure}
  \hfill
  \begin{subfigure}[t]{.24\textwidth}
    \centering
    \includegraphics[width=\linewidth]{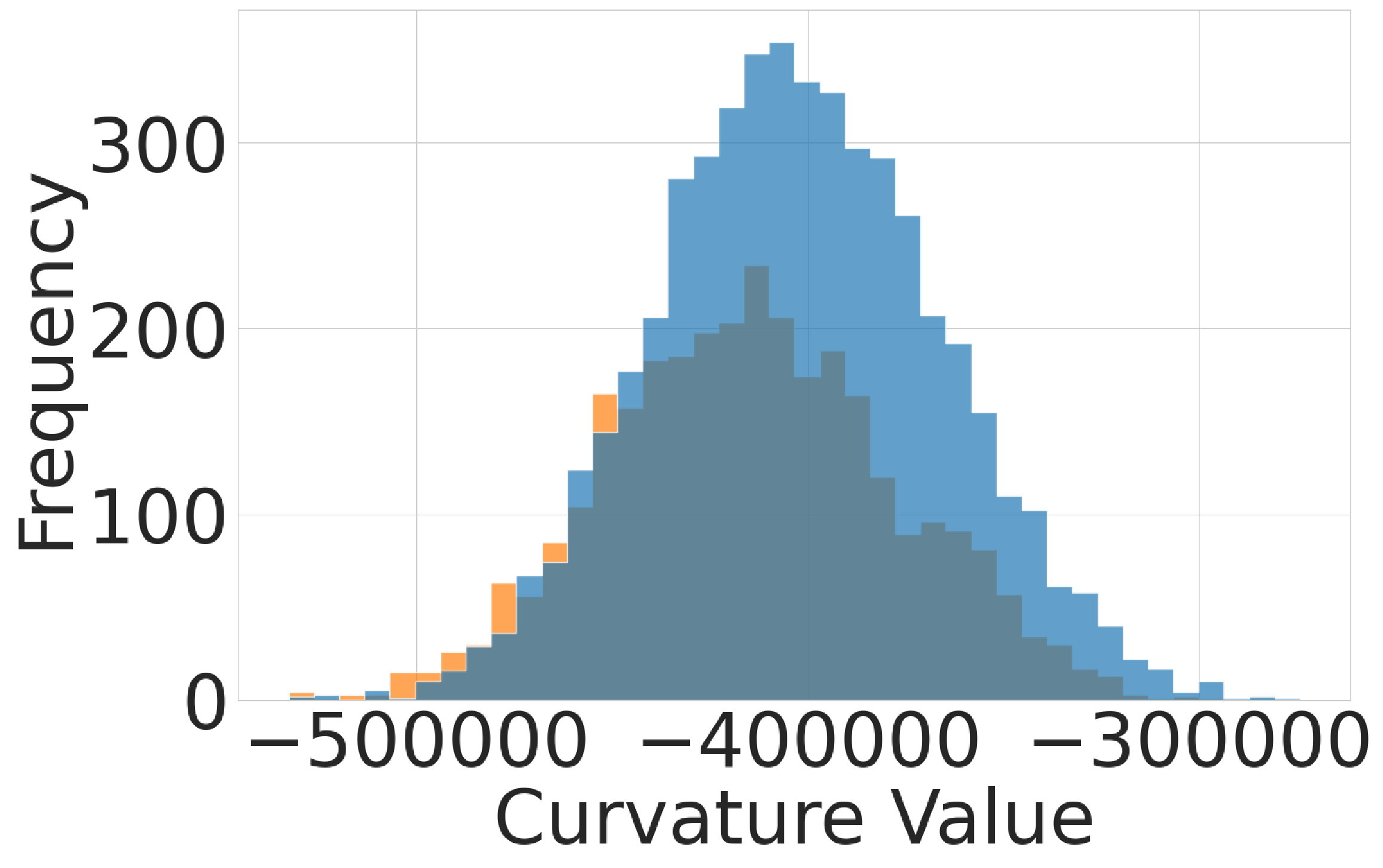}
  \end{subfigure}
    \begin{subfigure}[t]{.24\textwidth}
    \centering
    \includegraphics[width=\linewidth]{figures/appendix_figures/sbm_2_100_0.5_0.1_afrc_5.eps}
  \end{subfigure}
  \hfill
  \begin{subfigure}[t]{.24\textwidth}
    \centering
    \includegraphics[width=\linewidth]{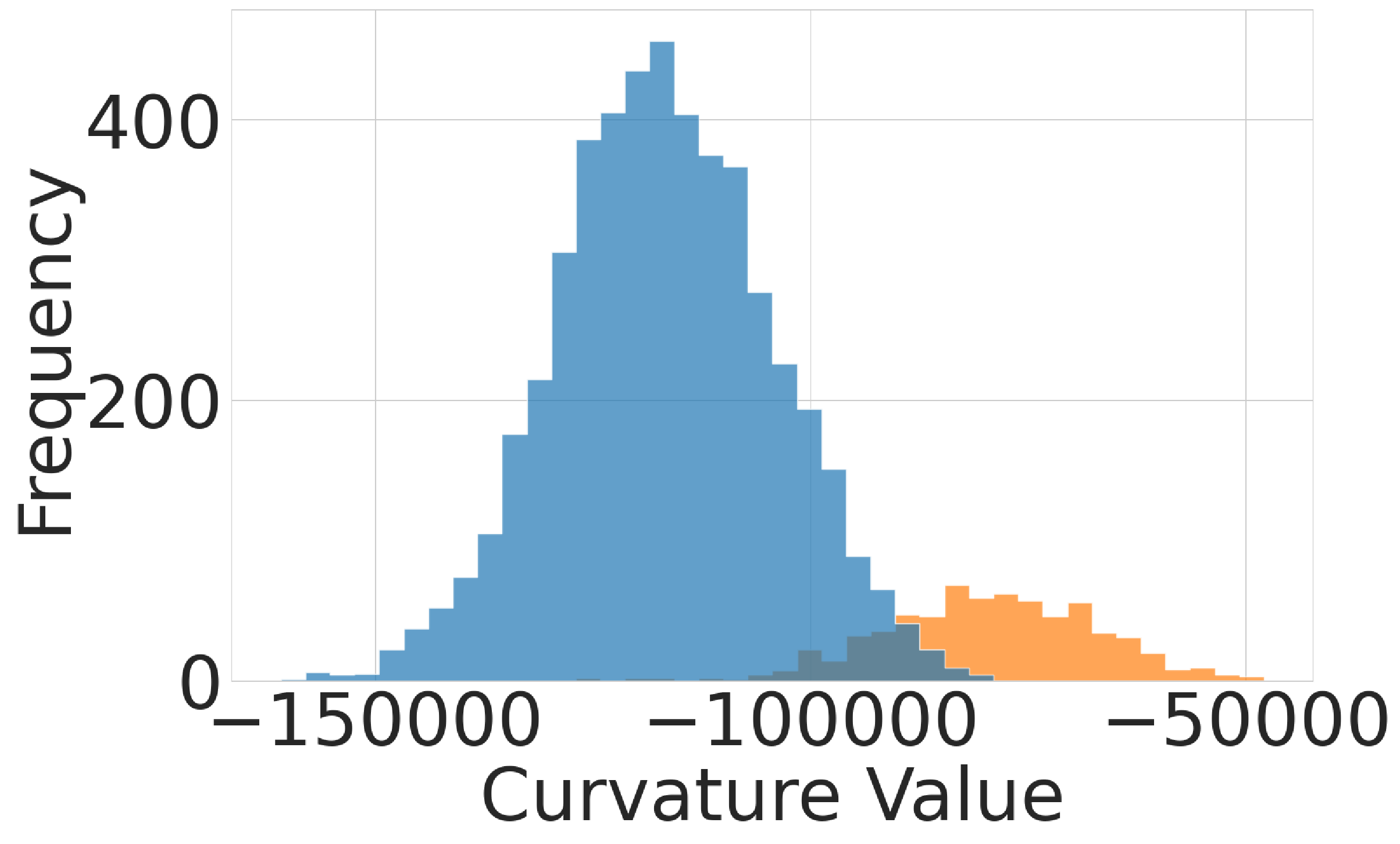}
  \end{subfigure}

  \medskip
  
  \begin{subfigure}[t]{.24\textwidth}
    \centering
    \includegraphics[width=\linewidth]{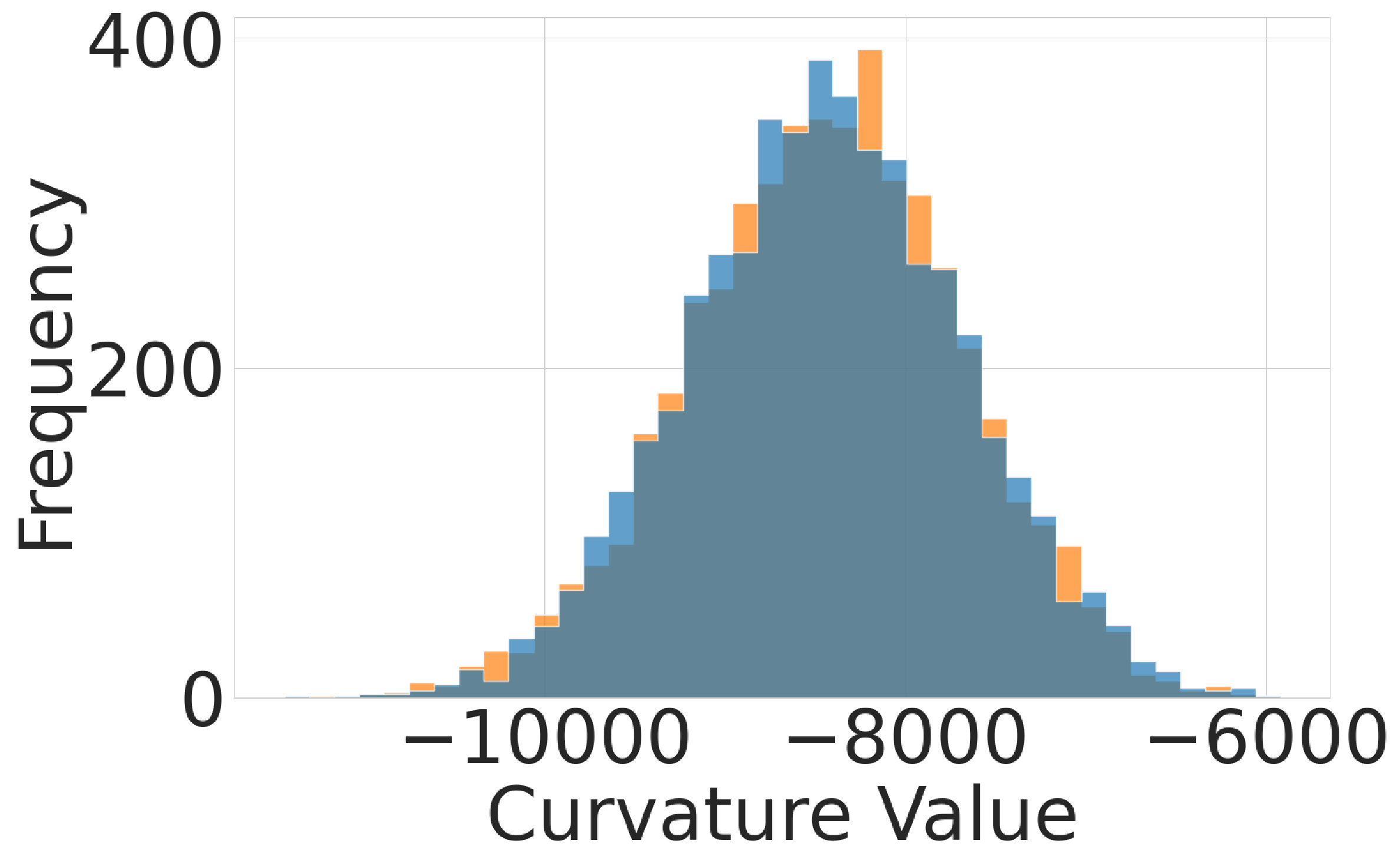}
  \end{subfigure}
  \hfill
  \begin{subfigure}[t]{.24\textwidth}
    \centering
    \includegraphics[width=\linewidth]{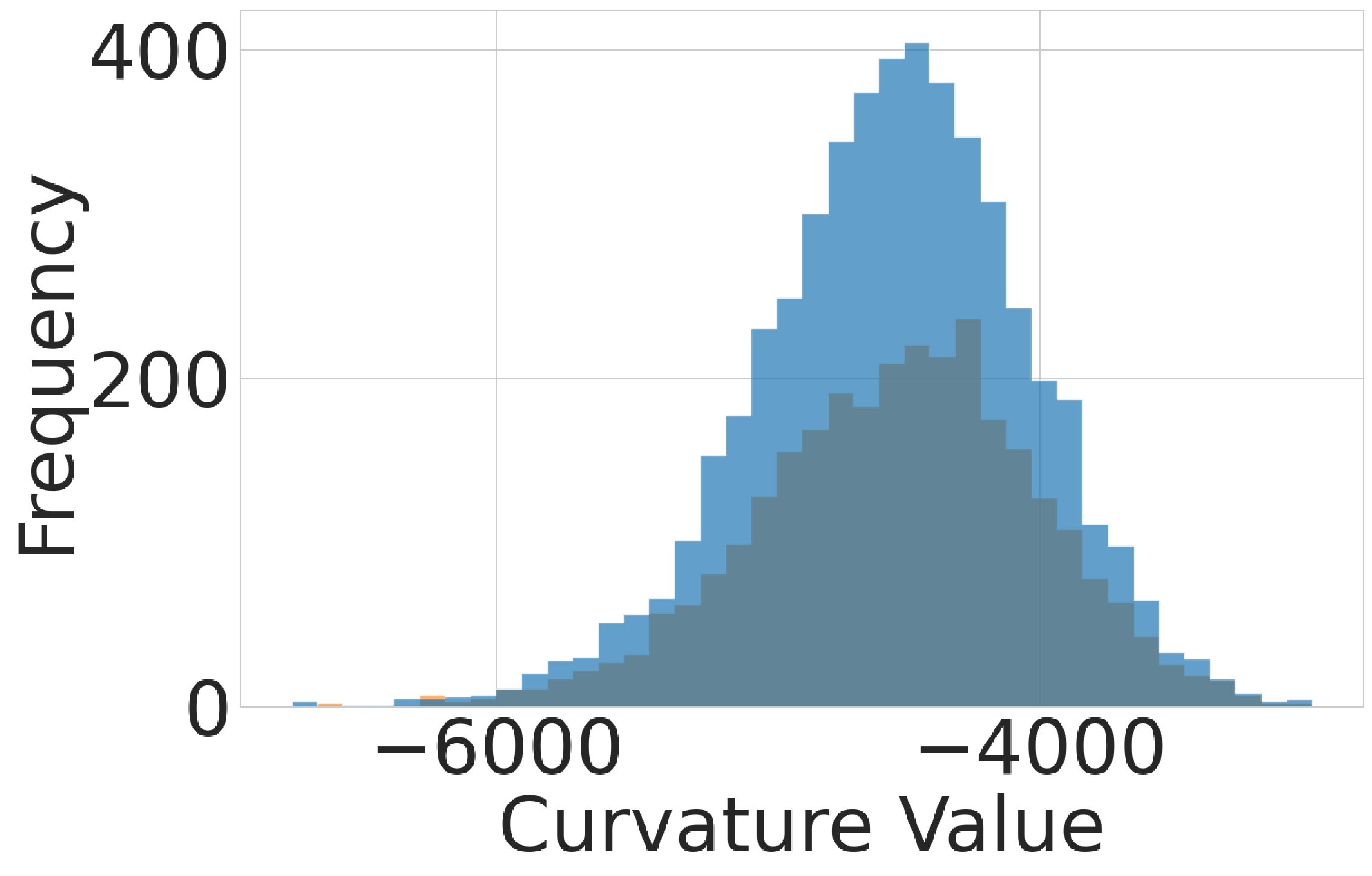}
  \end{subfigure}
    \begin{subfigure}[t]{.24\textwidth}
    \centering
    \includegraphics[width=\linewidth]{figures/appendix_figures/sbm_2_100_0.5_0.1_afrc.eps}
  \end{subfigure}
  \hfill
  \begin{subfigure}[t]{.24\textwidth}
    \centering
    \includegraphics[width=\linewidth]{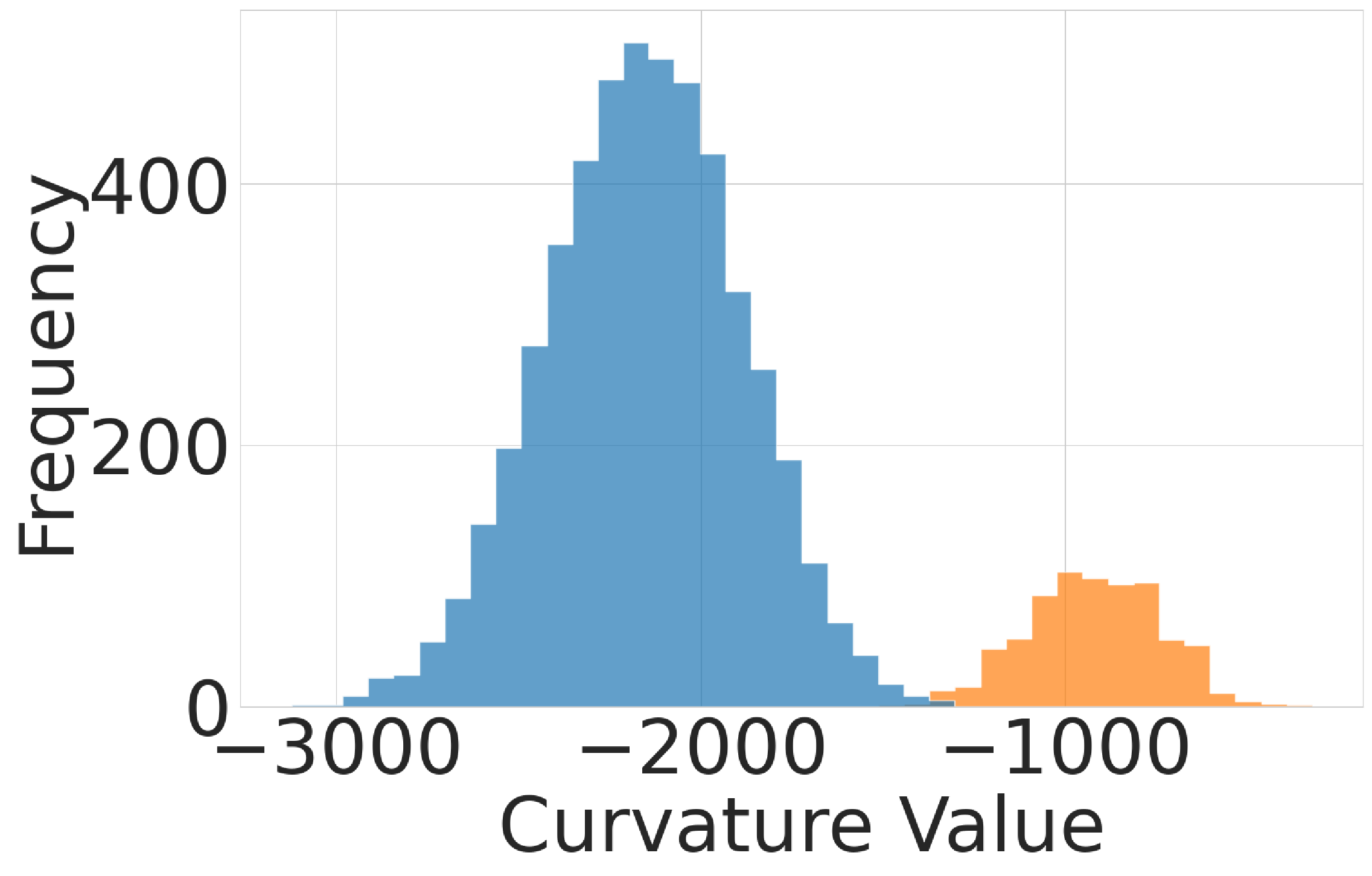}
  \end{subfigure}

  \medskip

    \begin{subfigure}[t]{.24\textwidth}
    \centering
    \includegraphics[width=\linewidth]{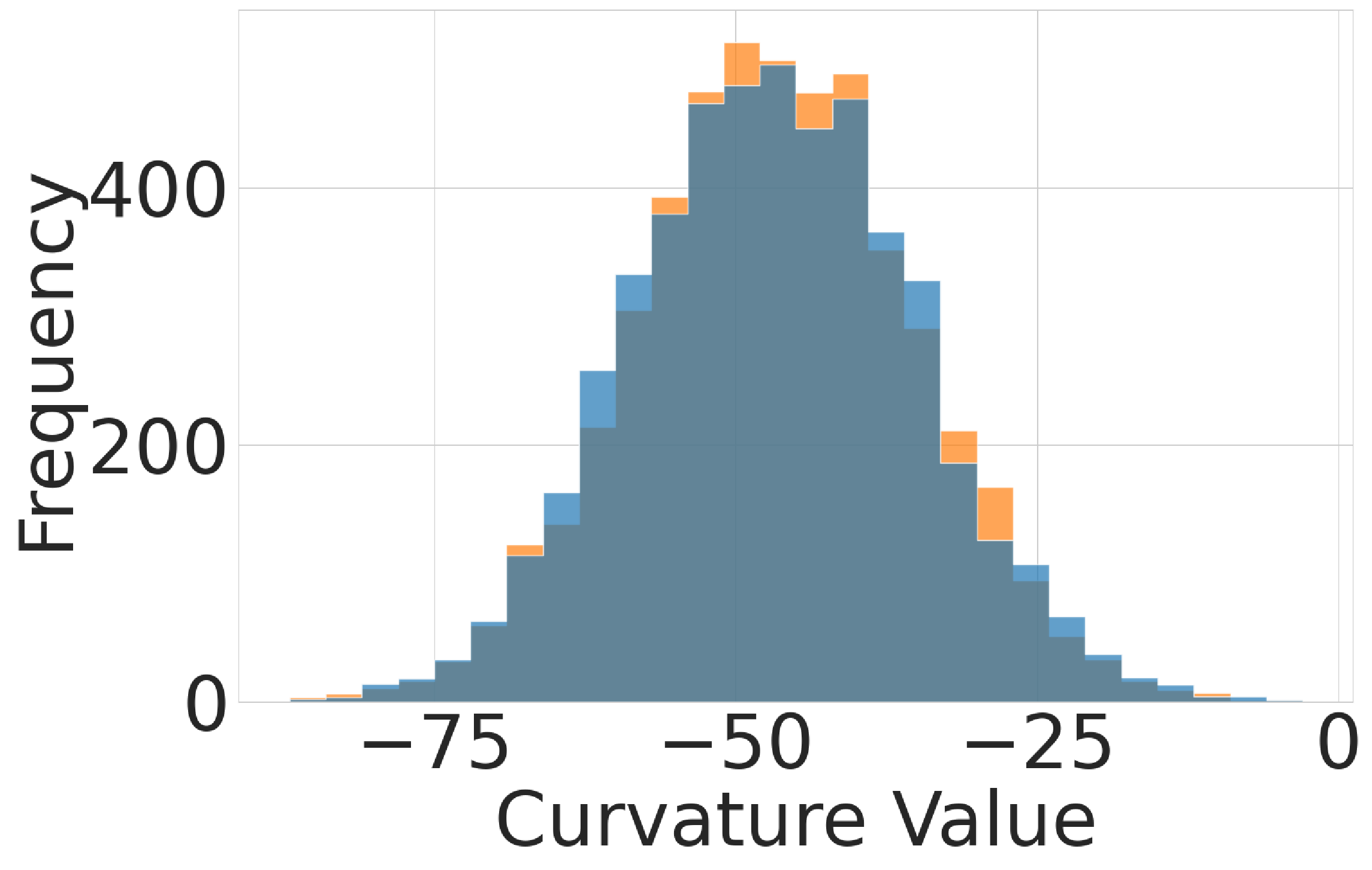}
  \end{subfigure}
  \hfill
  \begin{subfigure}[t]{.24\textwidth}
    \centering
    \includegraphics[width=\linewidth]{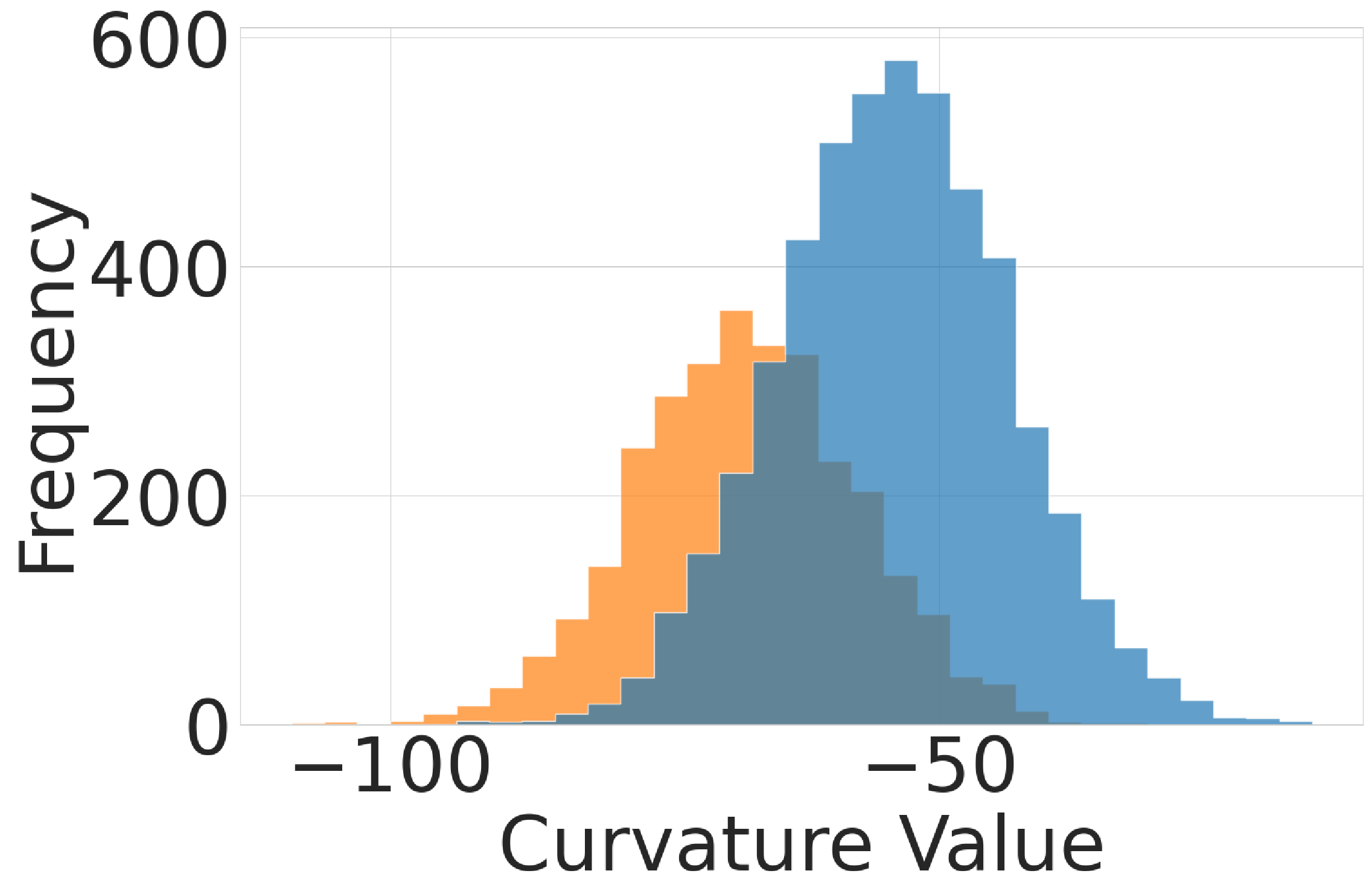}
  \end{subfigure}
    \begin{subfigure}[t]{.24\textwidth}
    \centering
    \includegraphics[width=\linewidth]{figures/appendix_figures/sbm_2_100_0.5_0.1_afrc_3.eps}
  \end{subfigure}
  \hfill
  \begin{subfigure}[t]{.24\textwidth}
    \centering
    \includegraphics[width=\linewidth]{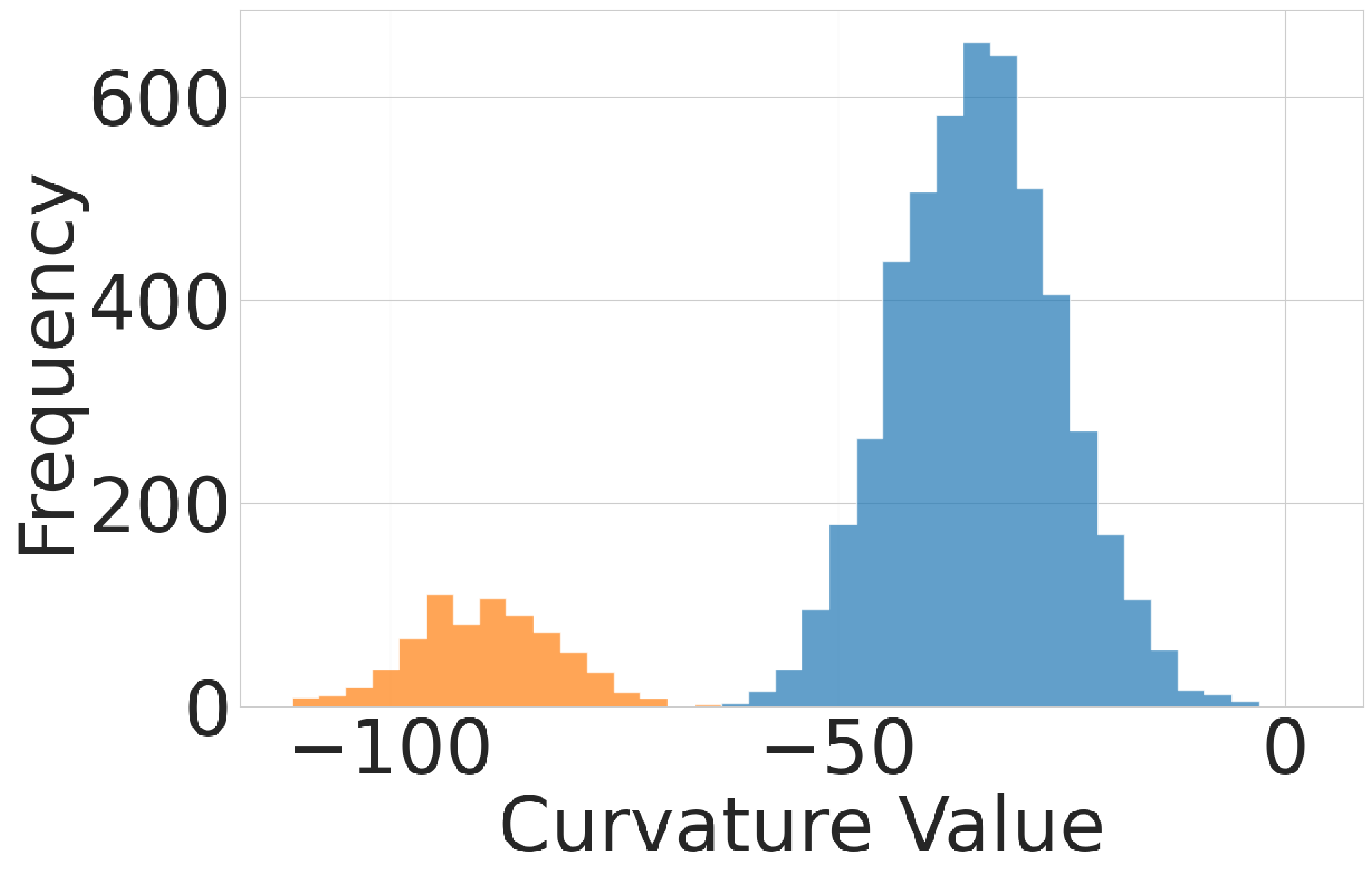}
  \end{subfigure}

  \medskip

  \begin{subfigure}[t]{.24\textwidth}
    \centering
    \includegraphics[width=\linewidth]{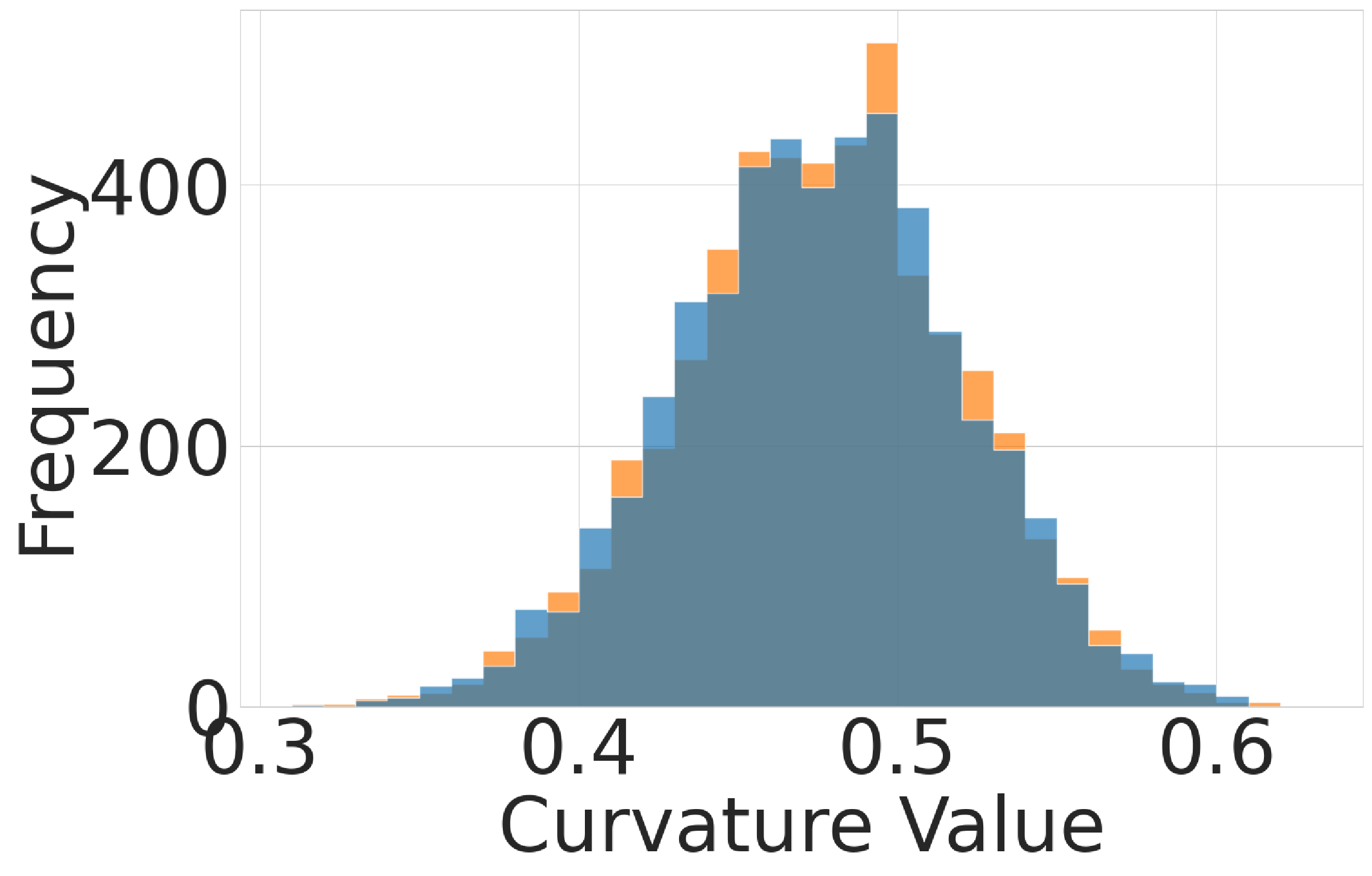}
  \end{subfigure}
  \hfill
  \begin{subfigure}[t]{.24\textwidth}
    \centering
    \includegraphics[width=\linewidth]{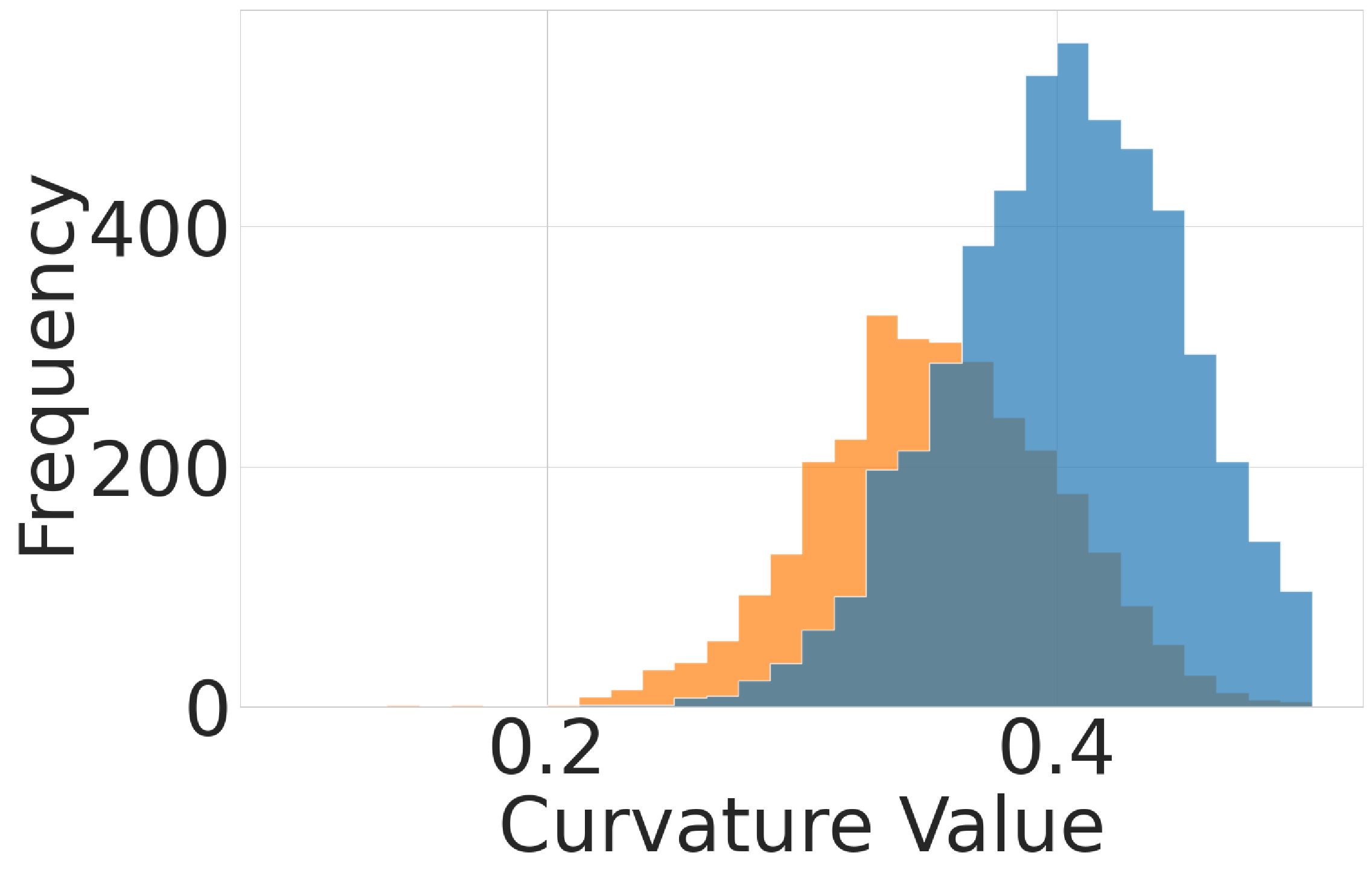}
  \end{subfigure}
    \begin{subfigure}[t]{.24\textwidth}
    \centering
    \includegraphics[width=\linewidth]{figures/appendix_figures/sbm_2_100_0.5_0.1_orc.eps}
  \end{subfigure}
  \hfill
  \begin{subfigure}[t]{.24\textwidth}
    \centering
    \includegraphics[width=\linewidth]{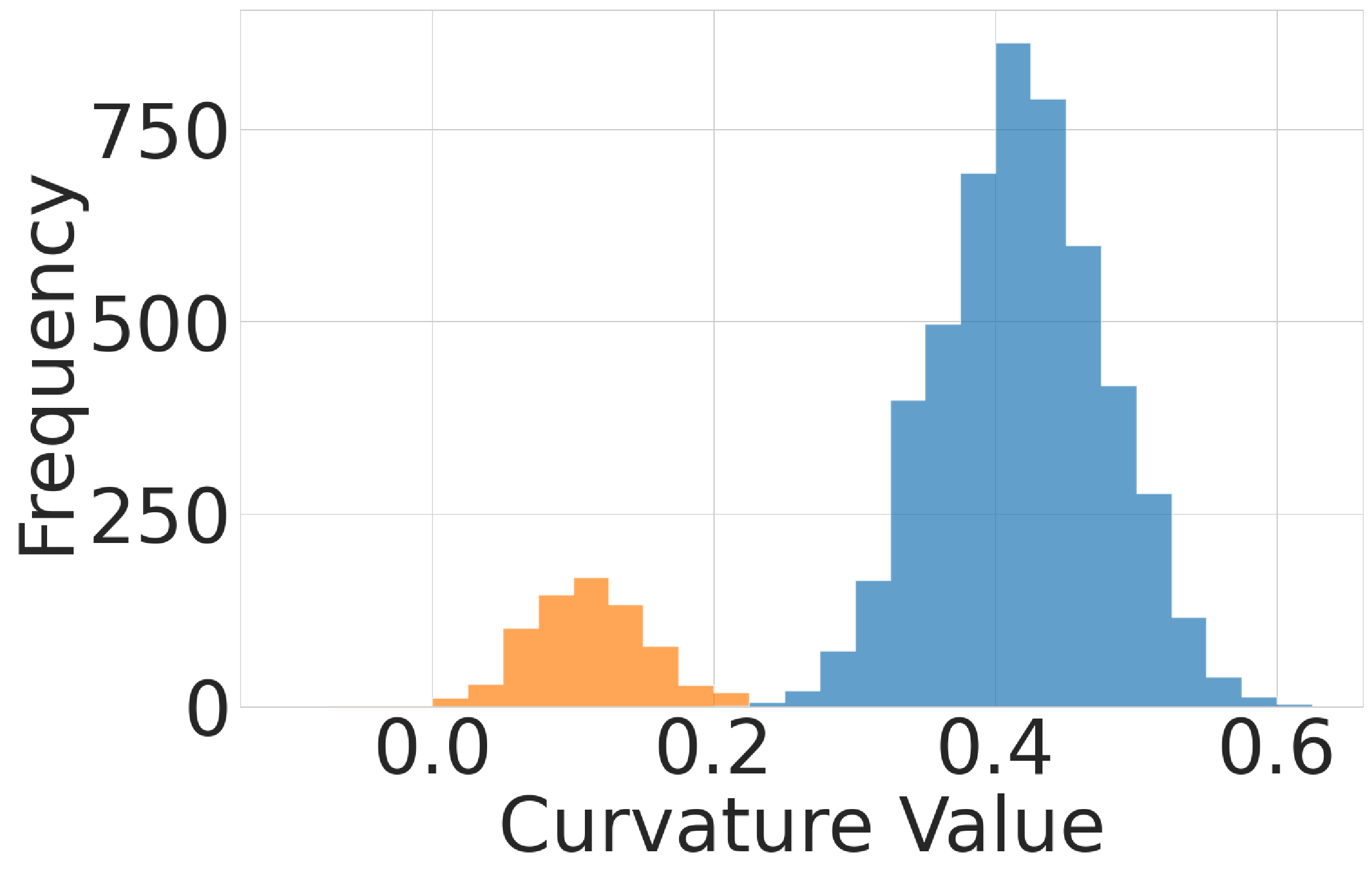}
  \end{subfigure}
  
  \caption{Top to bottom: Distributions of $\AF{5}$, $\AF{4}$, $\AF{3}$, and ORC for SBMs($2, 100, 0.5, q$) for $q \in \{0.5, 0.3, 0.1, 0.07\}$ (left to right). Edges within (between) communities in blue (orange).}
  \label{fig:AppC_2}
\end{figure}

\begin{figure}
  \begin{subfigure}[t]{.24\textwidth}
    \centering
    \includegraphics[width=\linewidth]{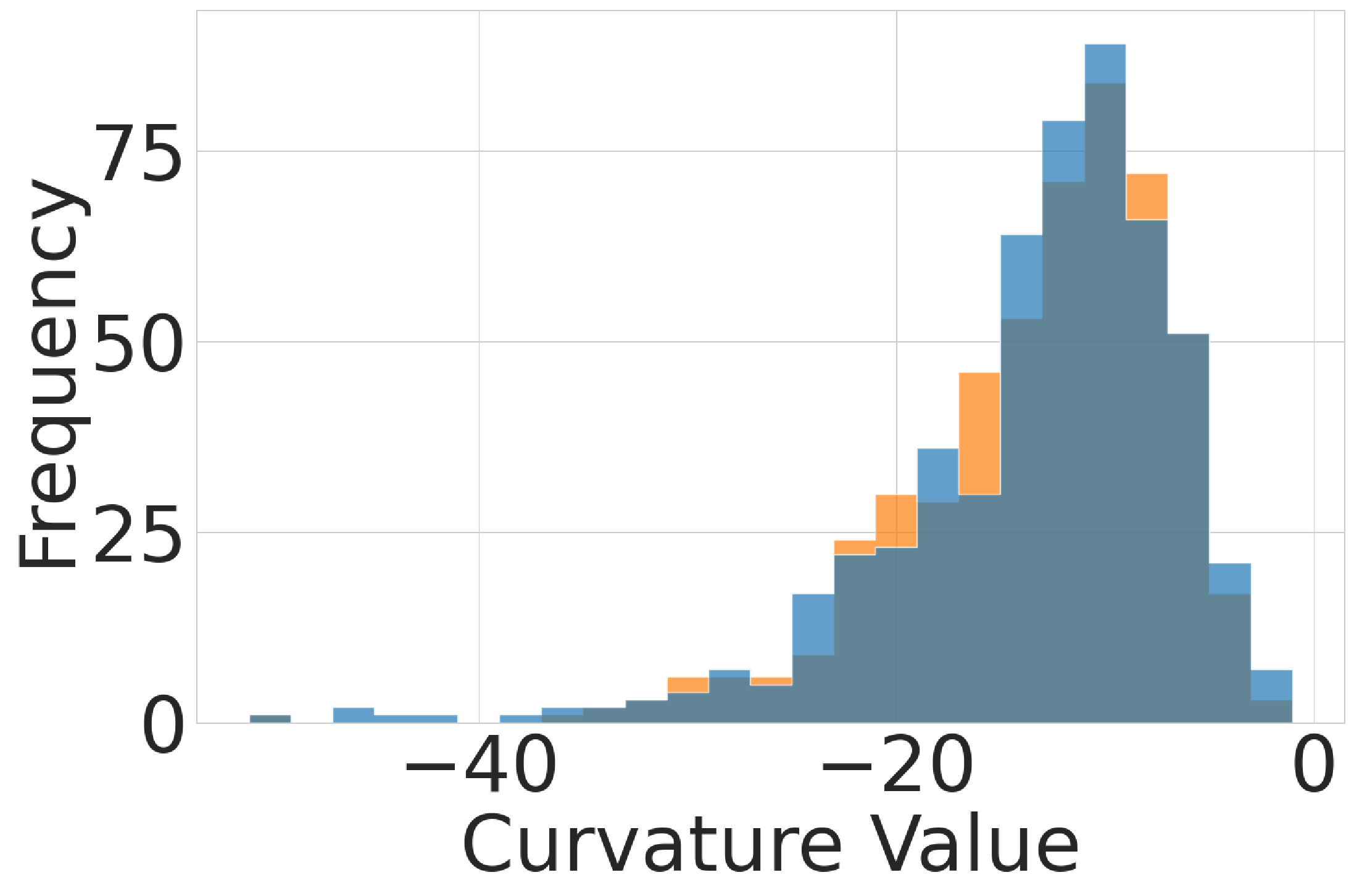}
  \end{subfigure}
  \hfill
  \begin{subfigure}[t]{.24\textwidth}
    \centering
    \includegraphics[width=\linewidth]{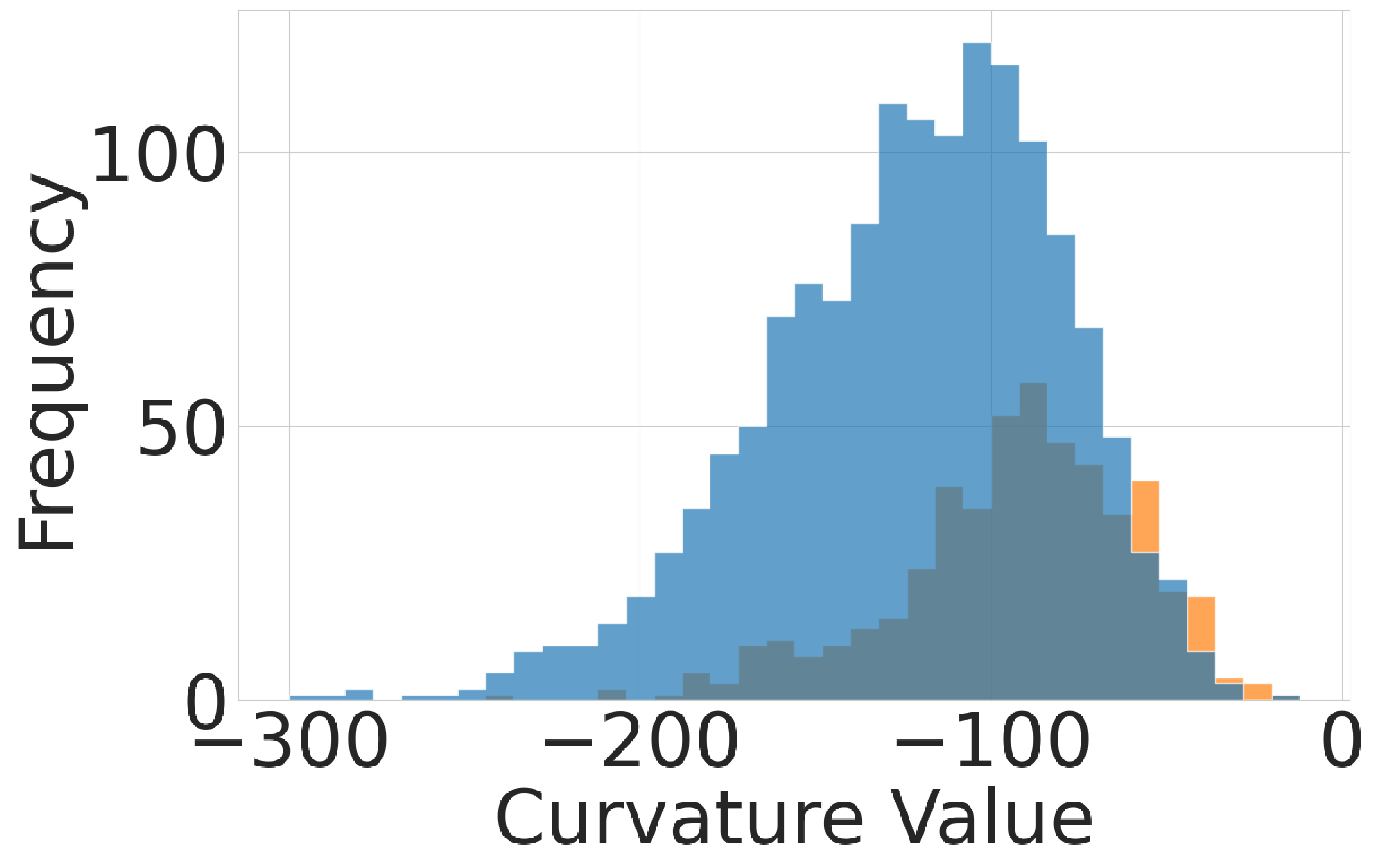}
  \end{subfigure}
    \begin{subfigure}[t]{.24\textwidth}
    \centering
    \includegraphics[width=\linewidth]{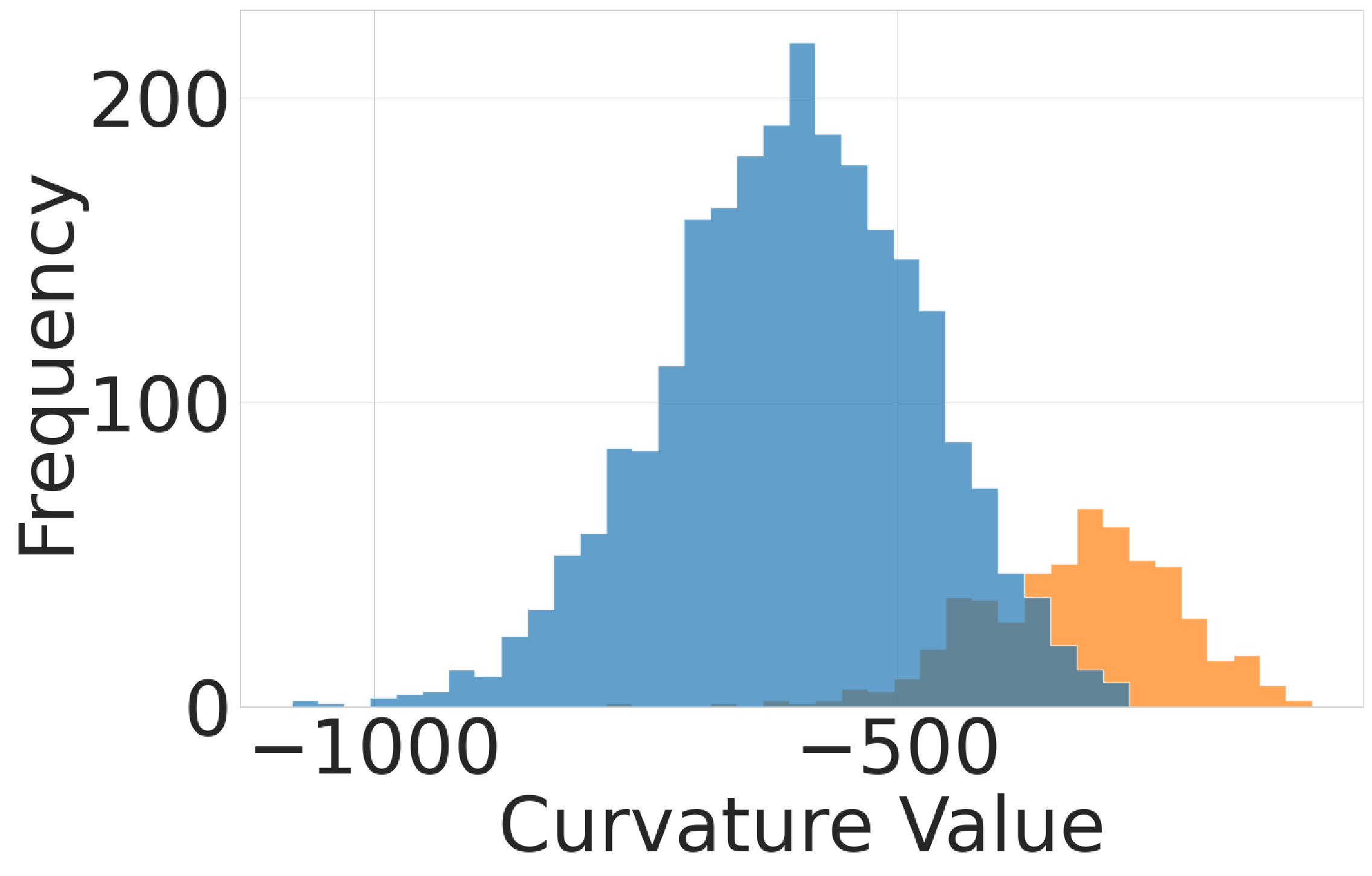}
  \end{subfigure}
  \hfill
  \begin{subfigure}[t]{.24\textwidth}
    \centering
    \includegraphics[width=\linewidth]{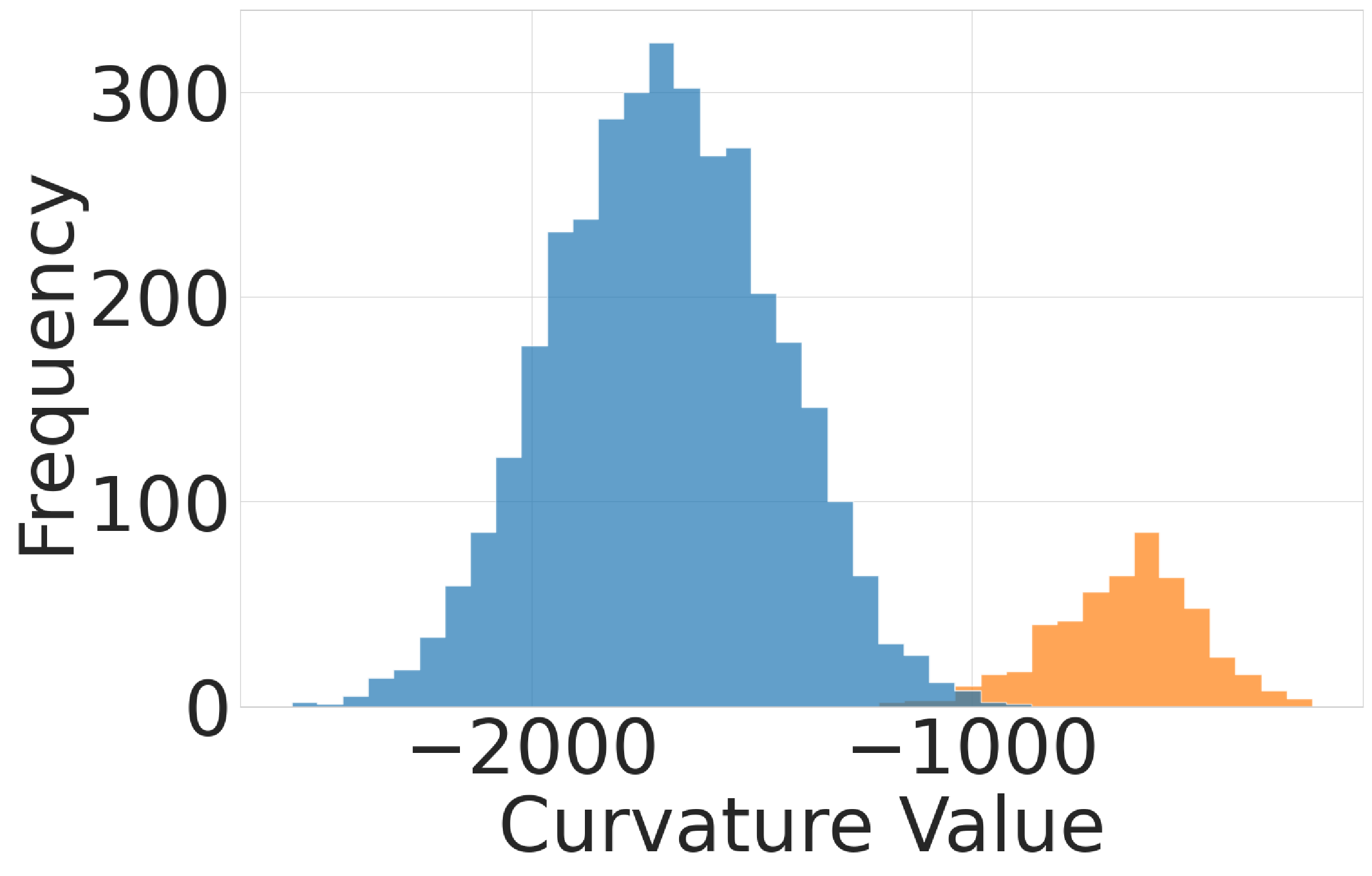}
  \end{subfigure}

  \medskip

  \begin{subfigure}[t]{.24\textwidth}
    \centering
    \includegraphics[width=\linewidth]{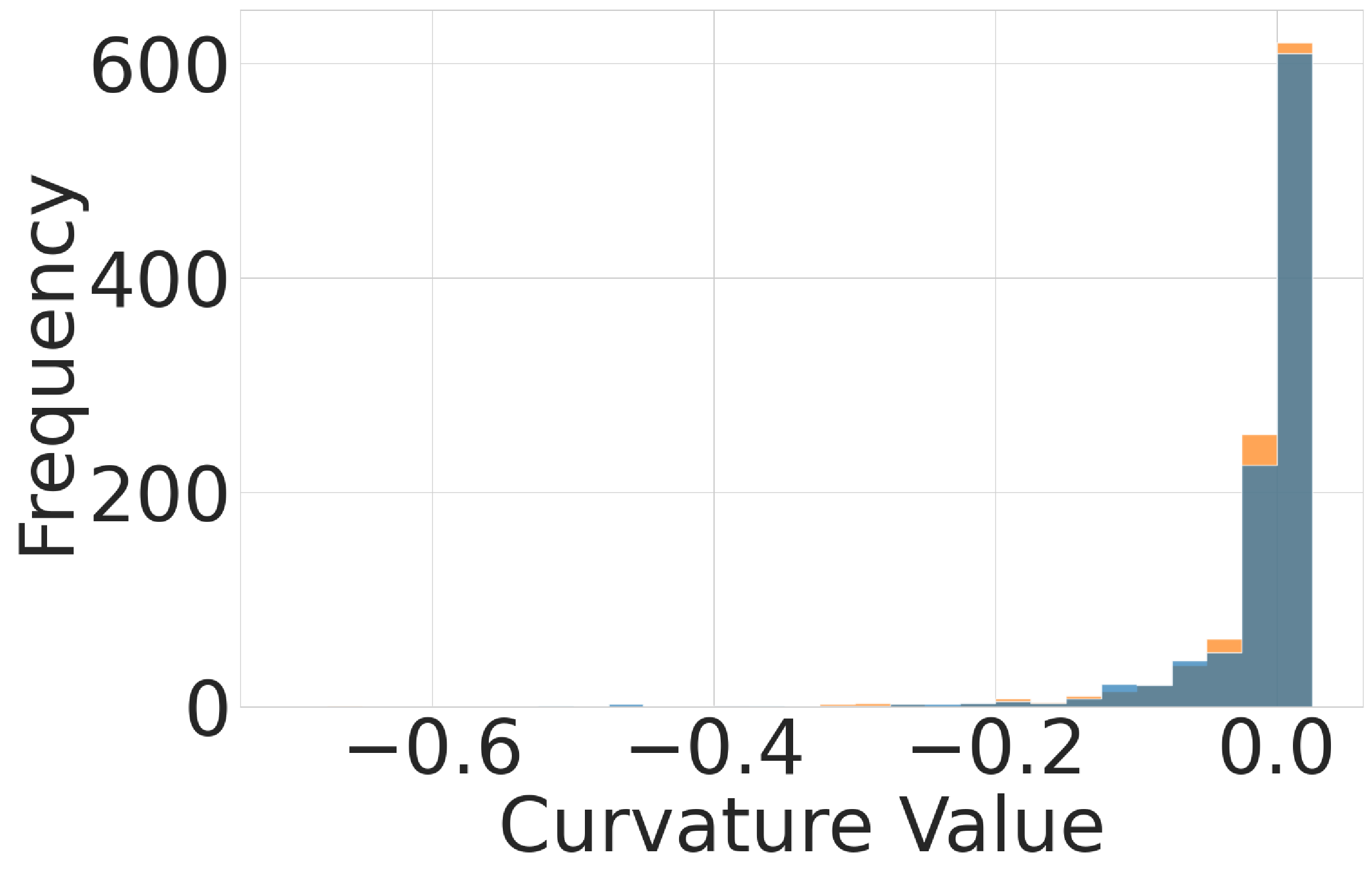}
  \end{subfigure}
  \hfill
  \begin{subfigure}[t]{.24\textwidth}
    \centering
    \includegraphics[width=\linewidth]{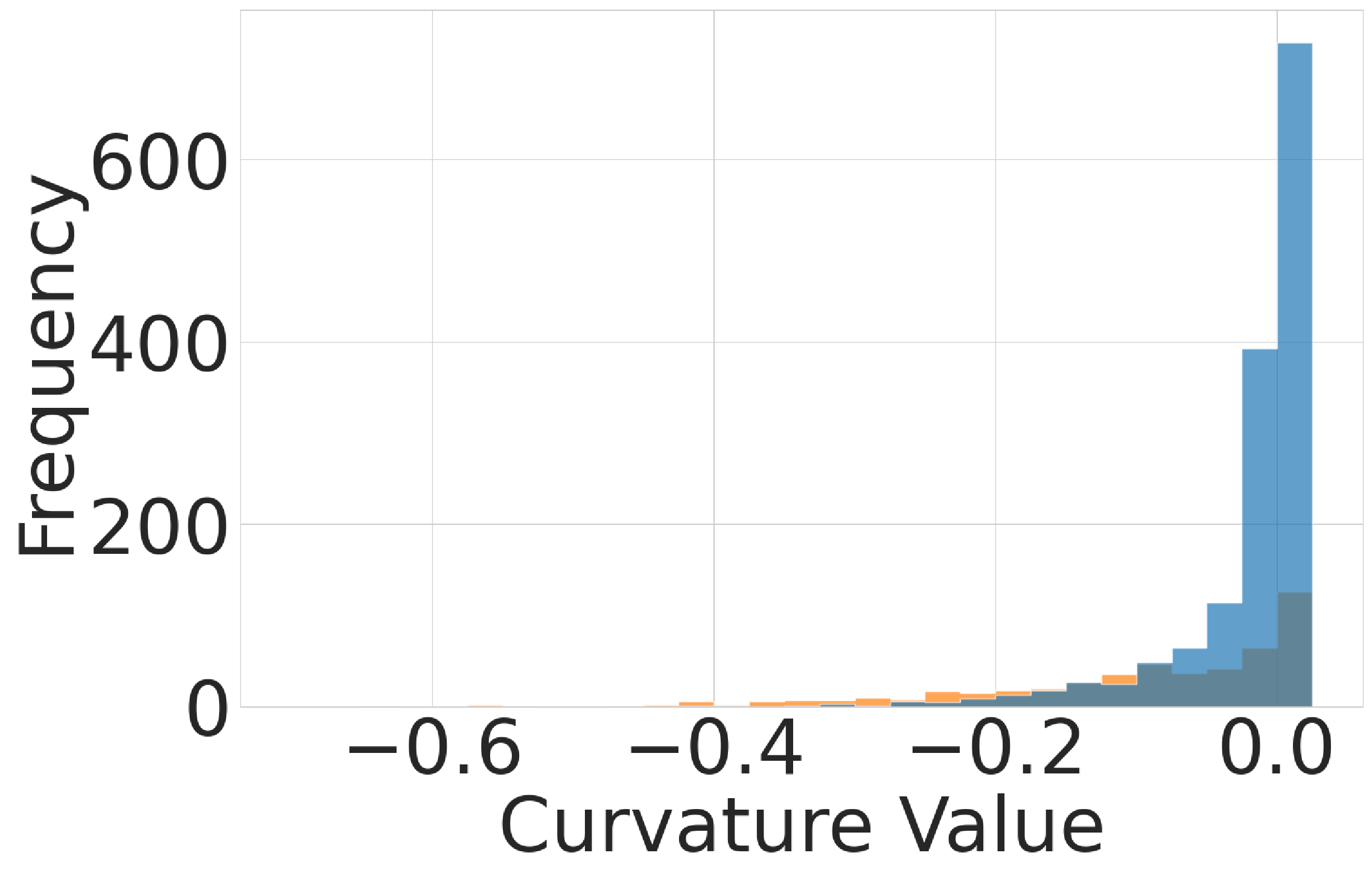}
  \end{subfigure}
    \begin{subfigure}[t]{.24\textwidth}
    \centering
    \includegraphics[width=\linewidth]{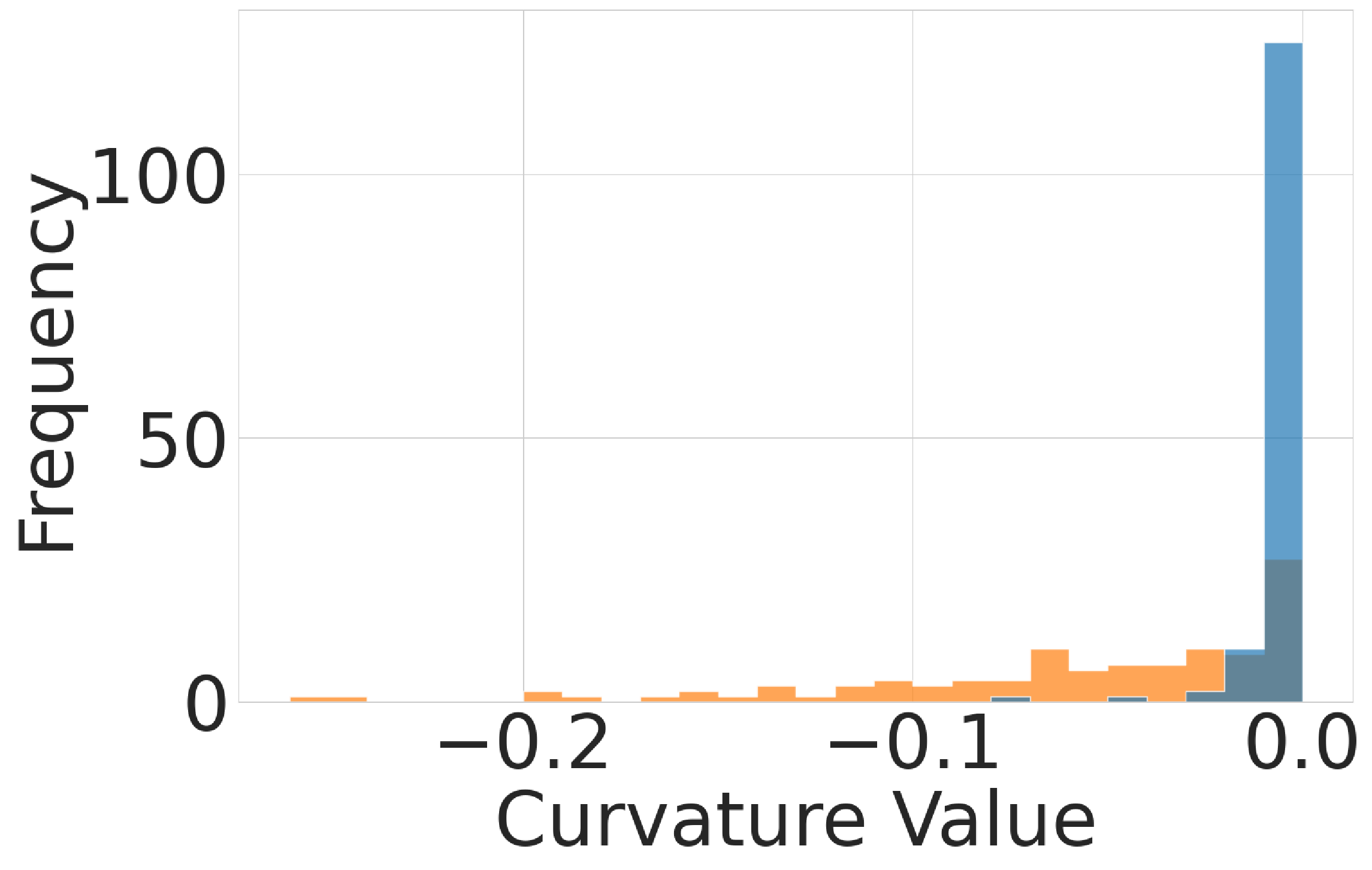}
  \end{subfigure}
  \hfill
  \begin{subfigure}[t]{.24\textwidth}
    \centering
    \includegraphics[width=\linewidth]{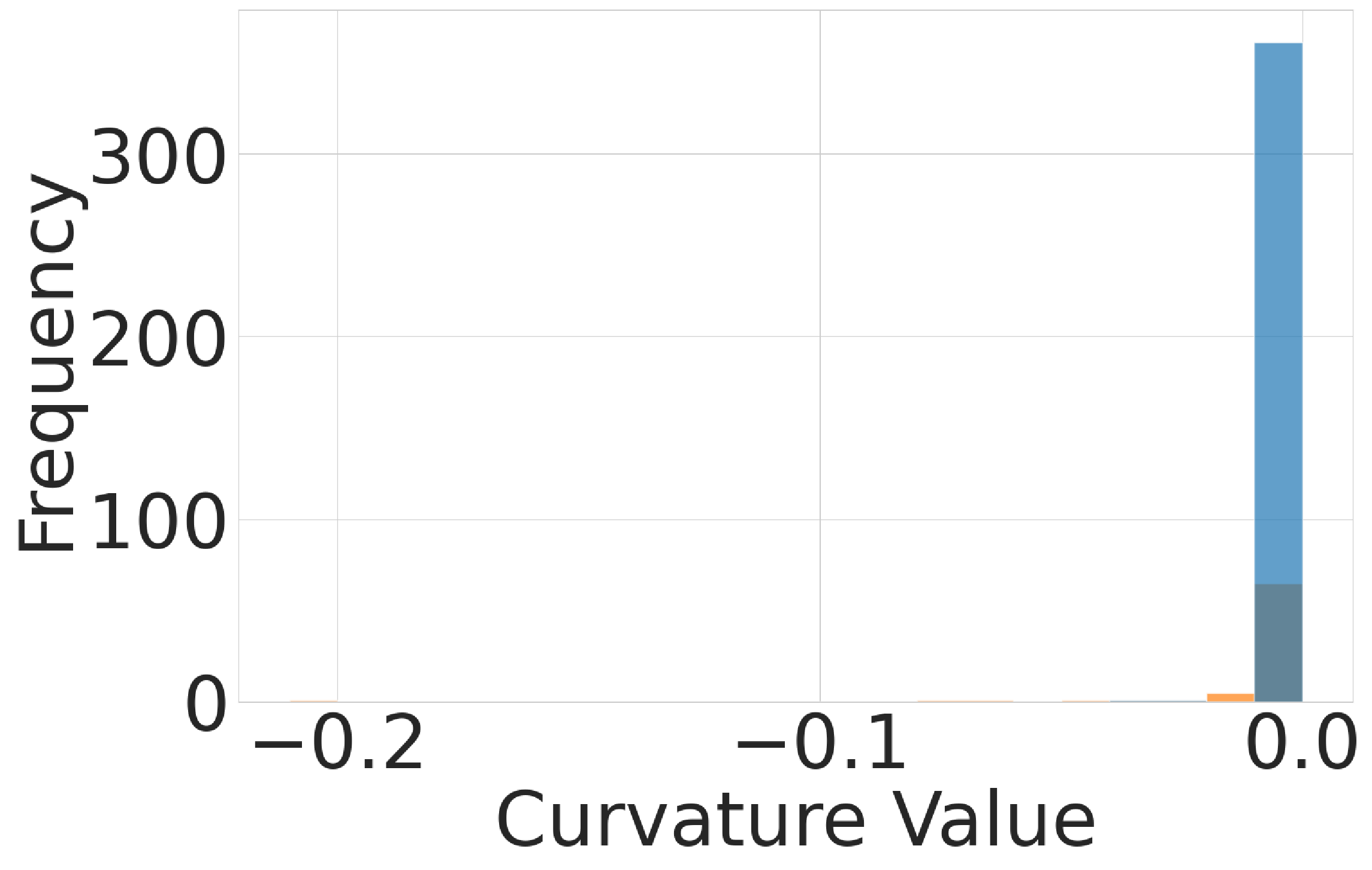}
  \end{subfigure}
  
  \caption{Distributions of the Augmented Forman-Ricci Curvature \textbf{(top)} and the Ollivier-Ricci Curvature \textbf{(bottom)} for HBGs($50, p, 0.1$) for $p \in \{0.1, 0.3, 0.5, 0.7\}$, from left to right.}
  \label{fig:AppC_3}
\end{figure}

\begin{figure}

  \begin{subfigure}[t]{.24\textwidth}
    \centering
    \includegraphics[width=\linewidth]{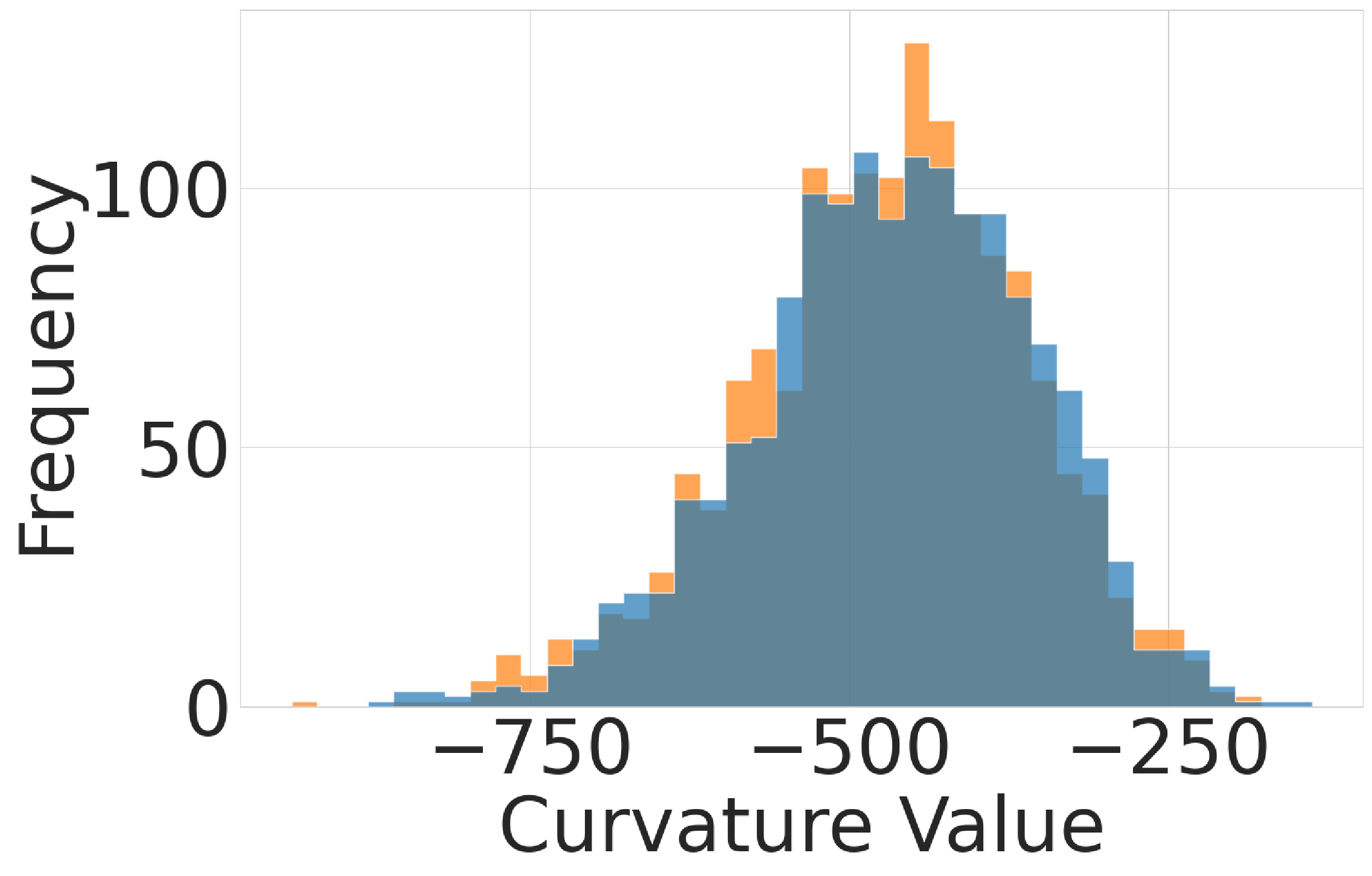}
  \end{subfigure}
  \hfill
  \begin{subfigure}[t]{.24\textwidth}
    \centering
    \includegraphics[width=\linewidth]{figures/appendix_figures/hbg_100_100_0.3_0.1_afrc.eps}
  \end{subfigure}
    \begin{subfigure}[t]{.24\textwidth}
    \centering
    \includegraphics[width=\linewidth]{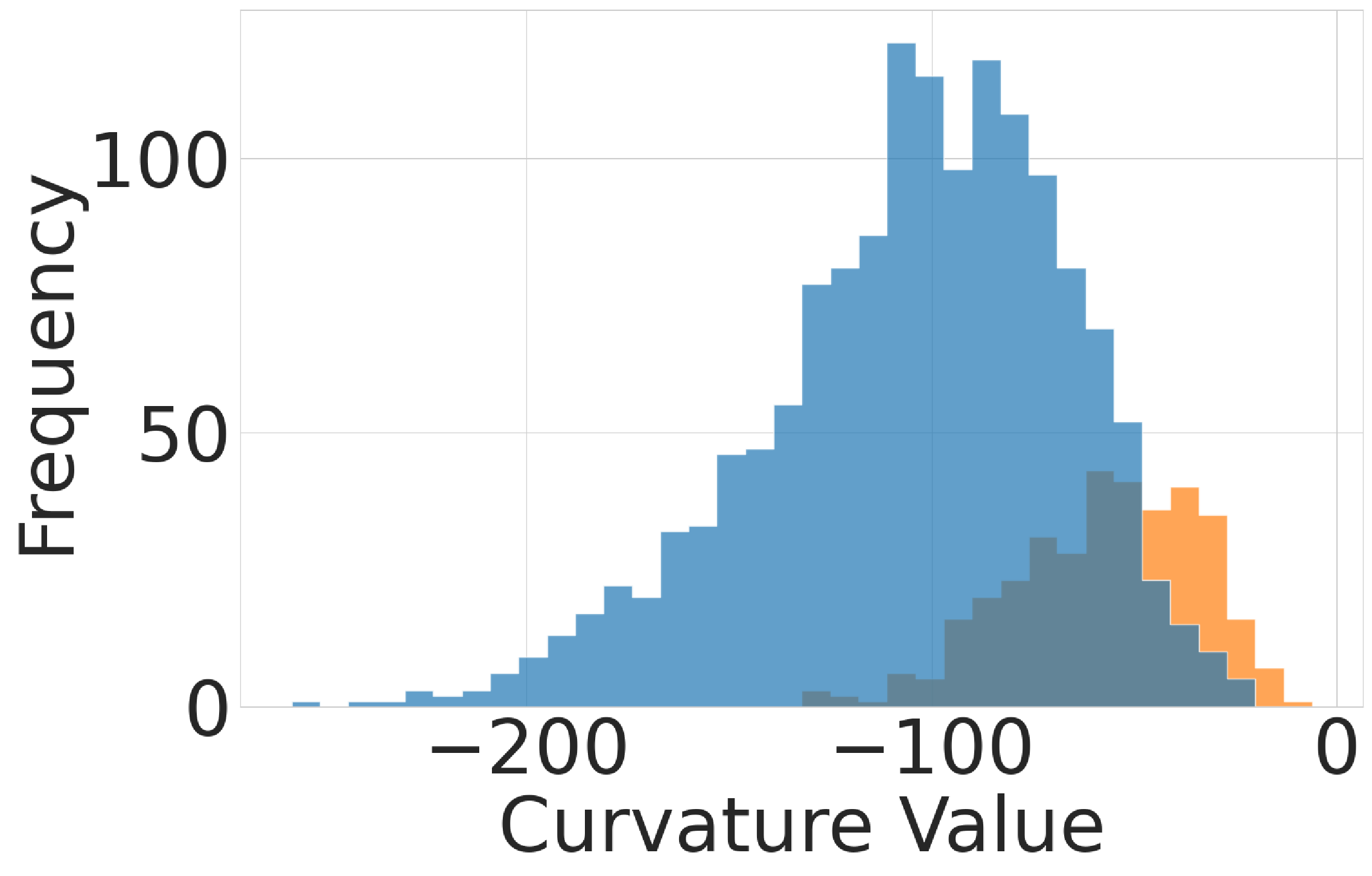}
  \end{subfigure}
  \hfill
  \begin{subfigure}[t]{.24\textwidth}
    \centering
    \includegraphics[width=\linewidth]{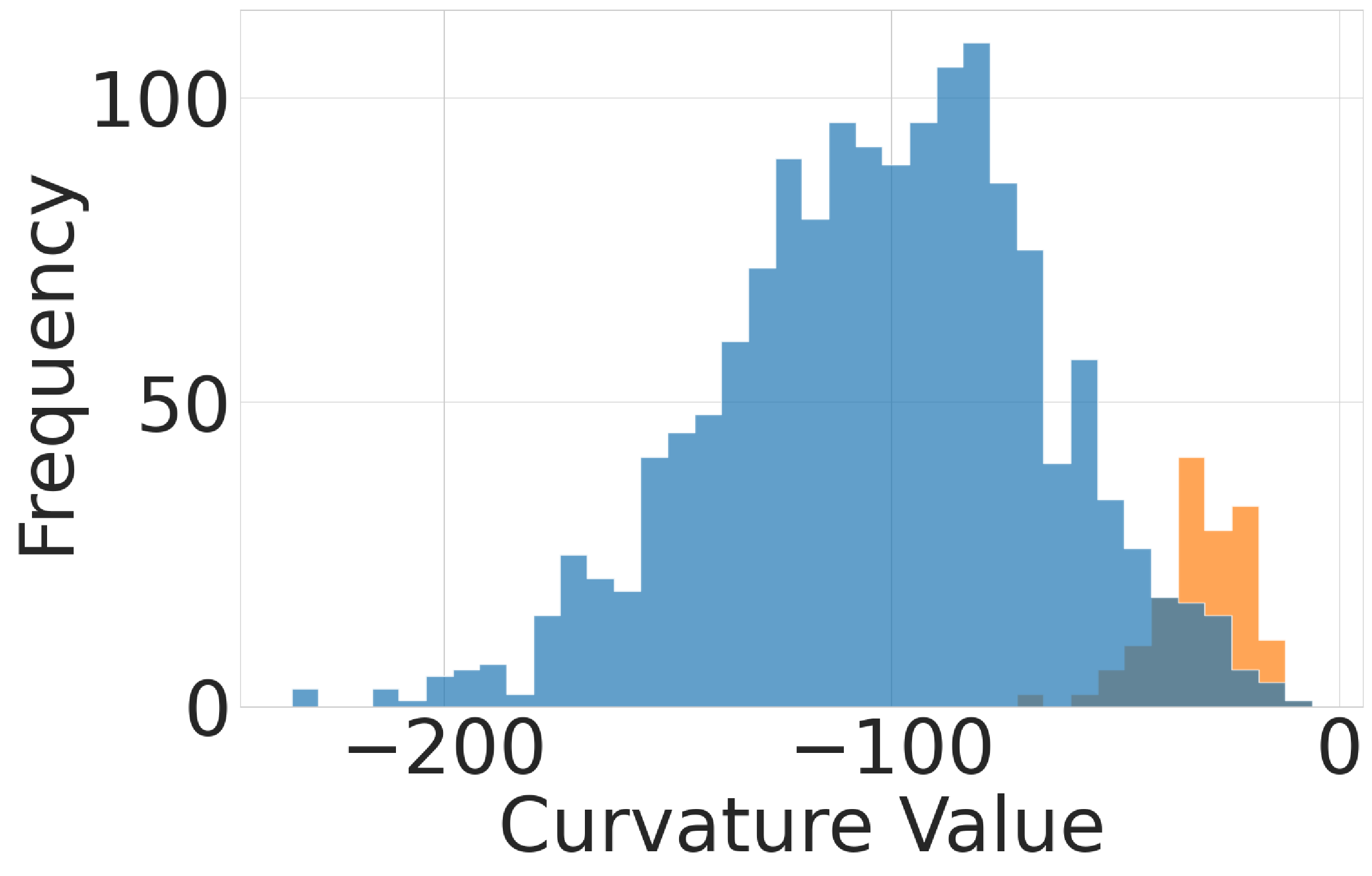}
  \end{subfigure}

  \medskip

  \begin{subfigure}[t]{.24\textwidth}
    \centering
    \includegraphics[width=\linewidth]{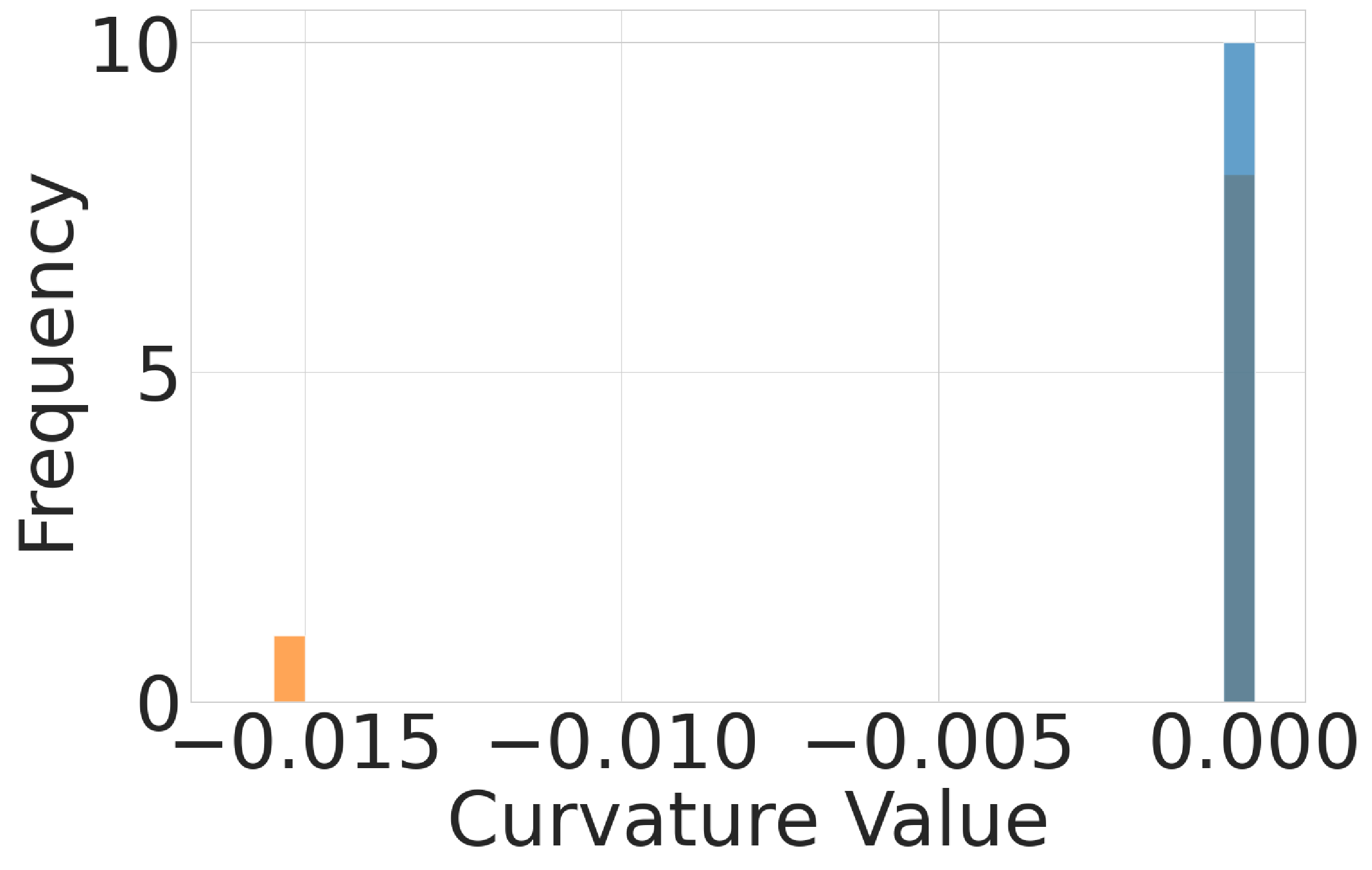}
  \end{subfigure}
  \hfill
  \begin{subfigure}[t]{.24\textwidth}
    \centering
    \includegraphics[width=\linewidth]{figures/appendix_figures/hbg_100_100_0.3_0.1_orc.eps}
  \end{subfigure}
    \begin{subfigure}[t]{.24\textwidth}
    \centering
    \includegraphics[width=\linewidth]{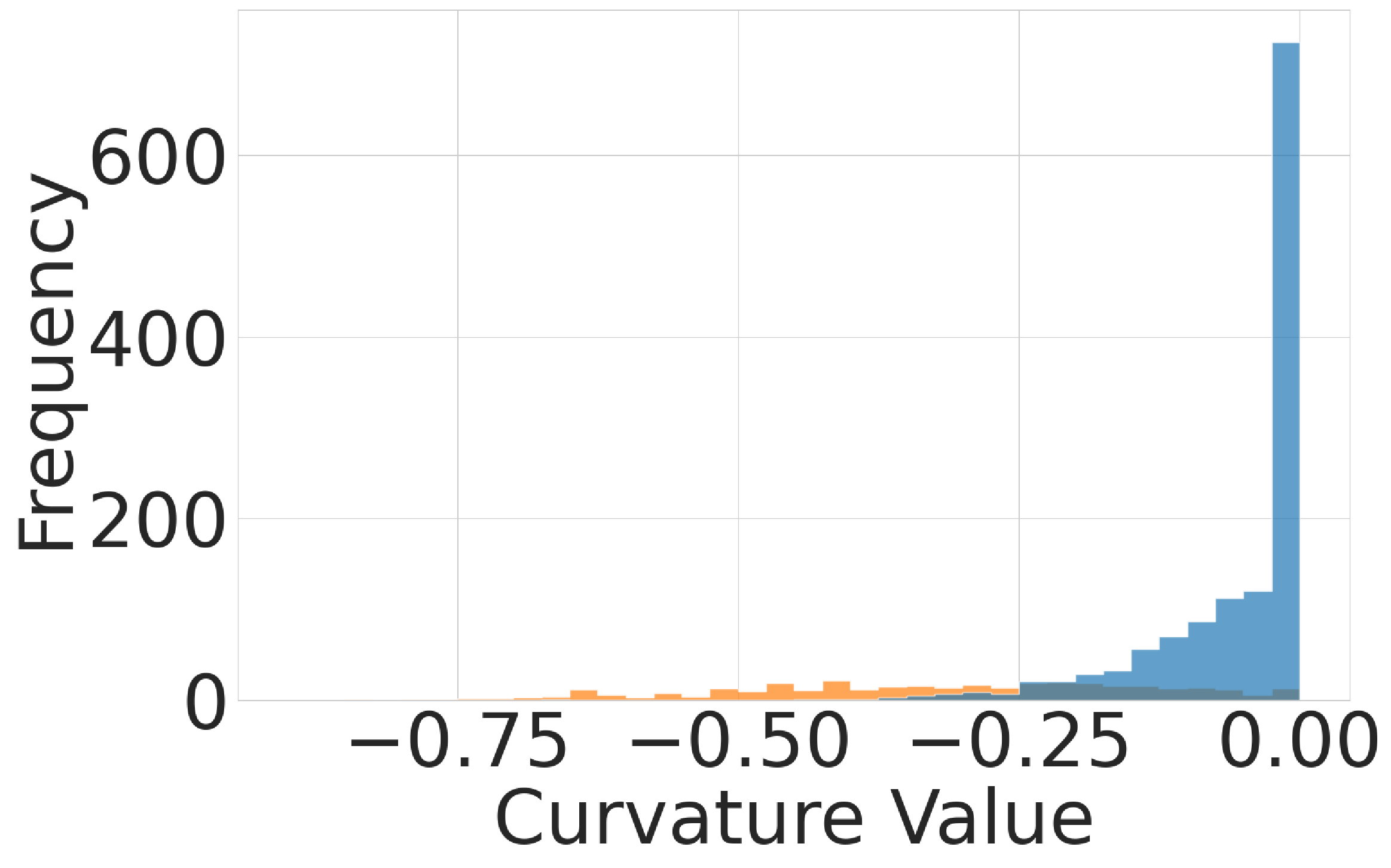}
  \end{subfigure}
  \hfill
  \begin{subfigure}[t]{.24\textwidth}
    \centering
    \includegraphics[width=\linewidth]{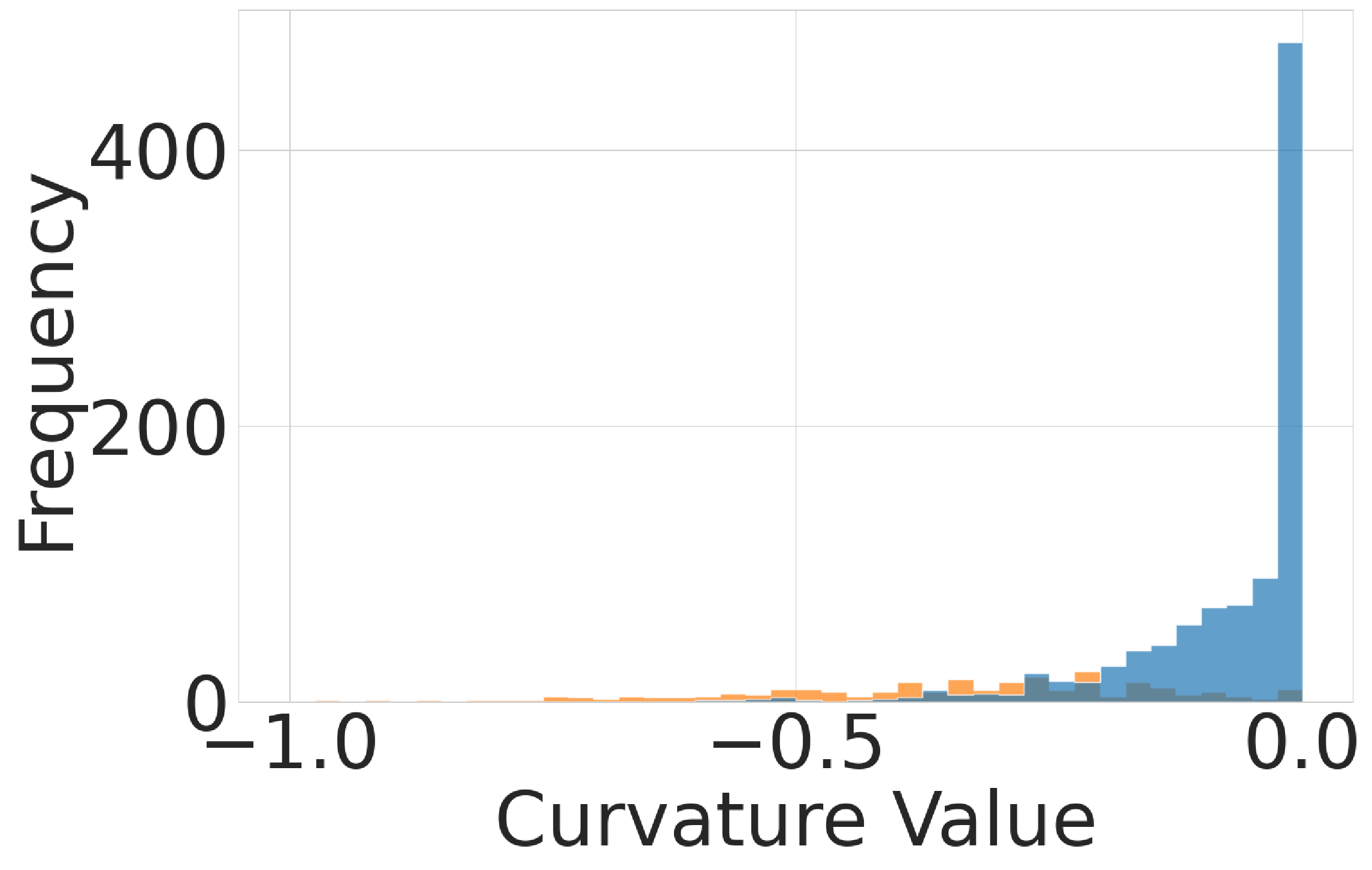}
  \end{subfigure}
  
  \caption{Distributions of the Augmented Forman-Ricci Curvature \textbf{(top)} and the Ollivier-Ricci Curvature \textbf{(bottom)} for HBGs($50, 0.3, q$) for $q \in \{0.3, 0.1, 0.07, 0.03\}$, from left to right.}
  \label{fig:AppC_4}
\end{figure}

\begin{table}[h!]
\centering
\begin{tabular}{ |l|c|c|c| } 
\hline
\textbf{Network} & \textbf{ORC} & \textbf{$\AF{3}$} & \textbf{$\AF{4}$} \\
\hline
t-SBM($2, 50, 0, 0.03$) & 0.52 & 0.46 & 0.40\\ 
\hline
t-SBM($2, 50, 0.1, 0.03$) & 2.14 & 0.93 & 0.53\\ 
\hline
t-SBM($2, 50, 0.2, 0.03$) & 3.98 & 1.72 & 1.41\\ 
\hline
t-SBM($2, 50, 0.3, 0.03$) & 6.18 & 3.04 & 3.29\\ 
\hline
\end{tabular}
\caption{Comparison of curvature gaps attained by the Ollivier-Ricci curvature and augmentations of the Forman-Ricci curvature for tree-SBMs with increasing probability of additional edges within communities.}
\label{table:5}
\end{table}

\begin{figure}[h!]
    \centering
    \includegraphics[width= 0.7\linewidth]{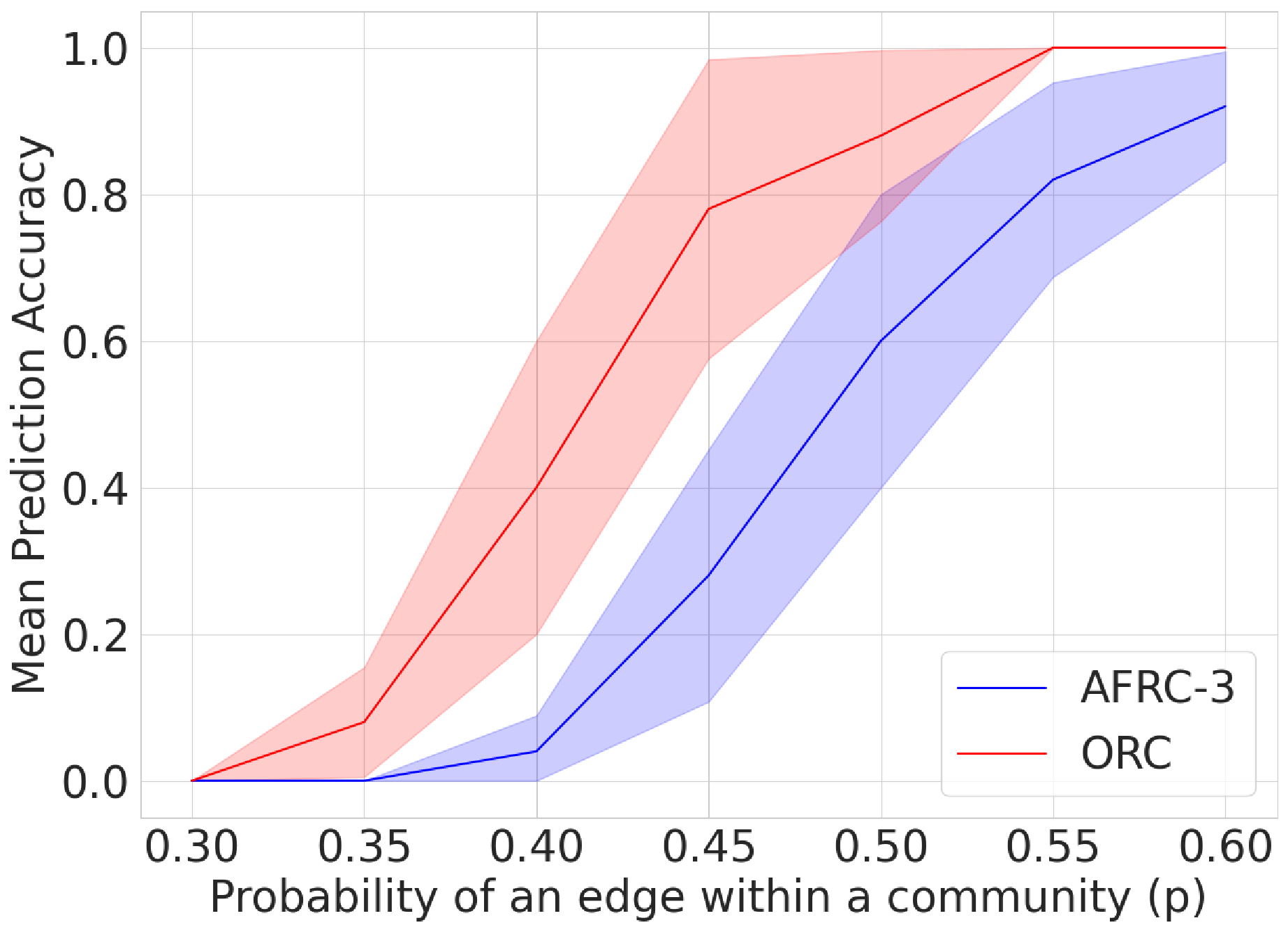}
    \caption{Mean accuracy (percentage of correctly identified communities) and standard deviation (shaded areas) of 15 ORC-based (red) and $\AF{3}$ community-detection (blue) for an SBM with parameters $(5, 20, p, 0.05)$ for $p \in \{0.3, 0.35, 0.4, 0.45, 0.5, 0.55, 0.6 \}$.}
    \label{fig:sbm_p_in}
\end{figure}

\begin{figure}[!t]
  \begin{subfigure}[t]{.39\textwidth}
    \centering
    \includegraphics[width=\linewidth]{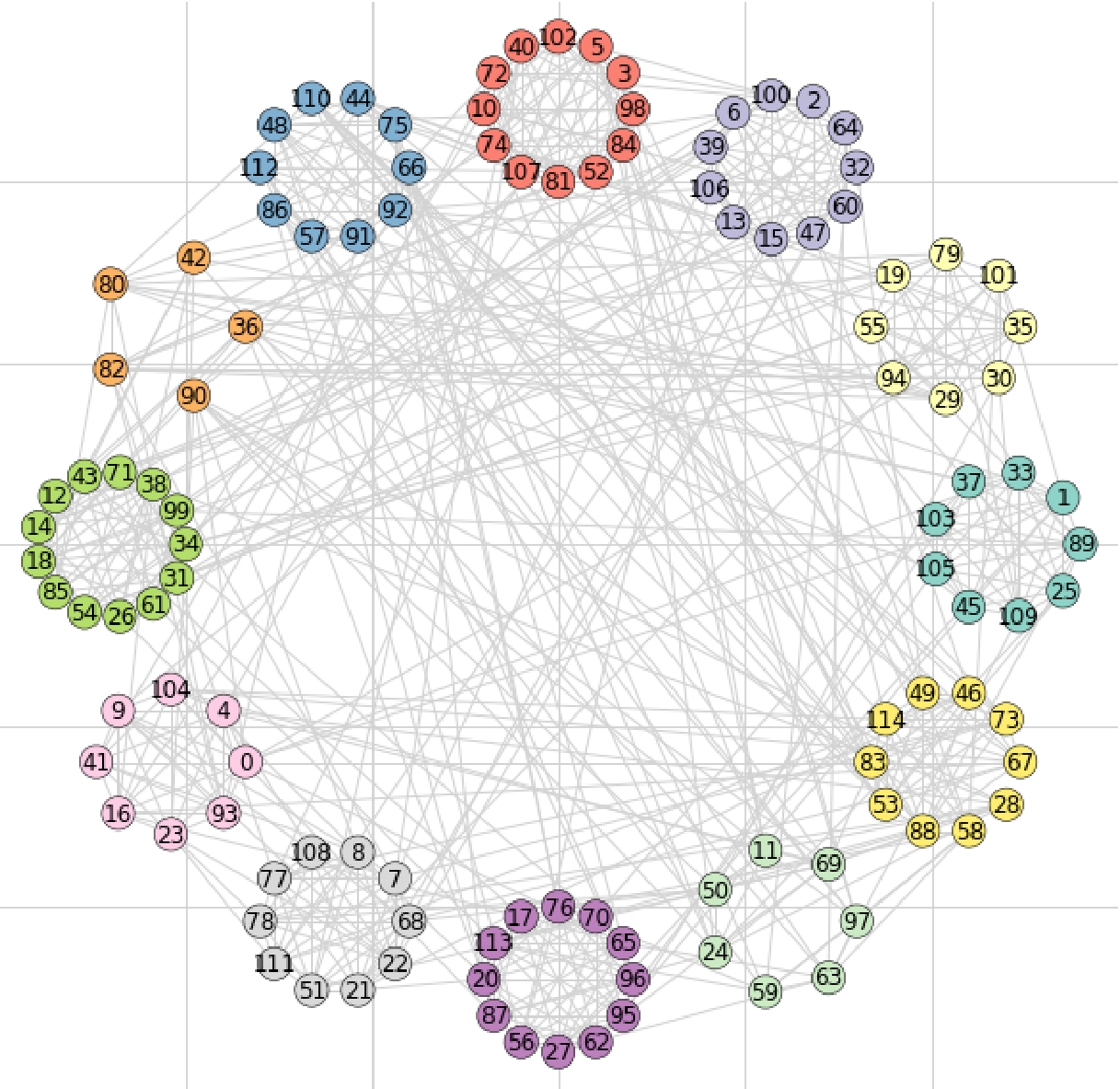}
    \caption{}
  \end{subfigure}
  \hfill
  \begin{subfigure}[t]{.39\textwidth}
    \centering
    \includegraphics[width=\linewidth]{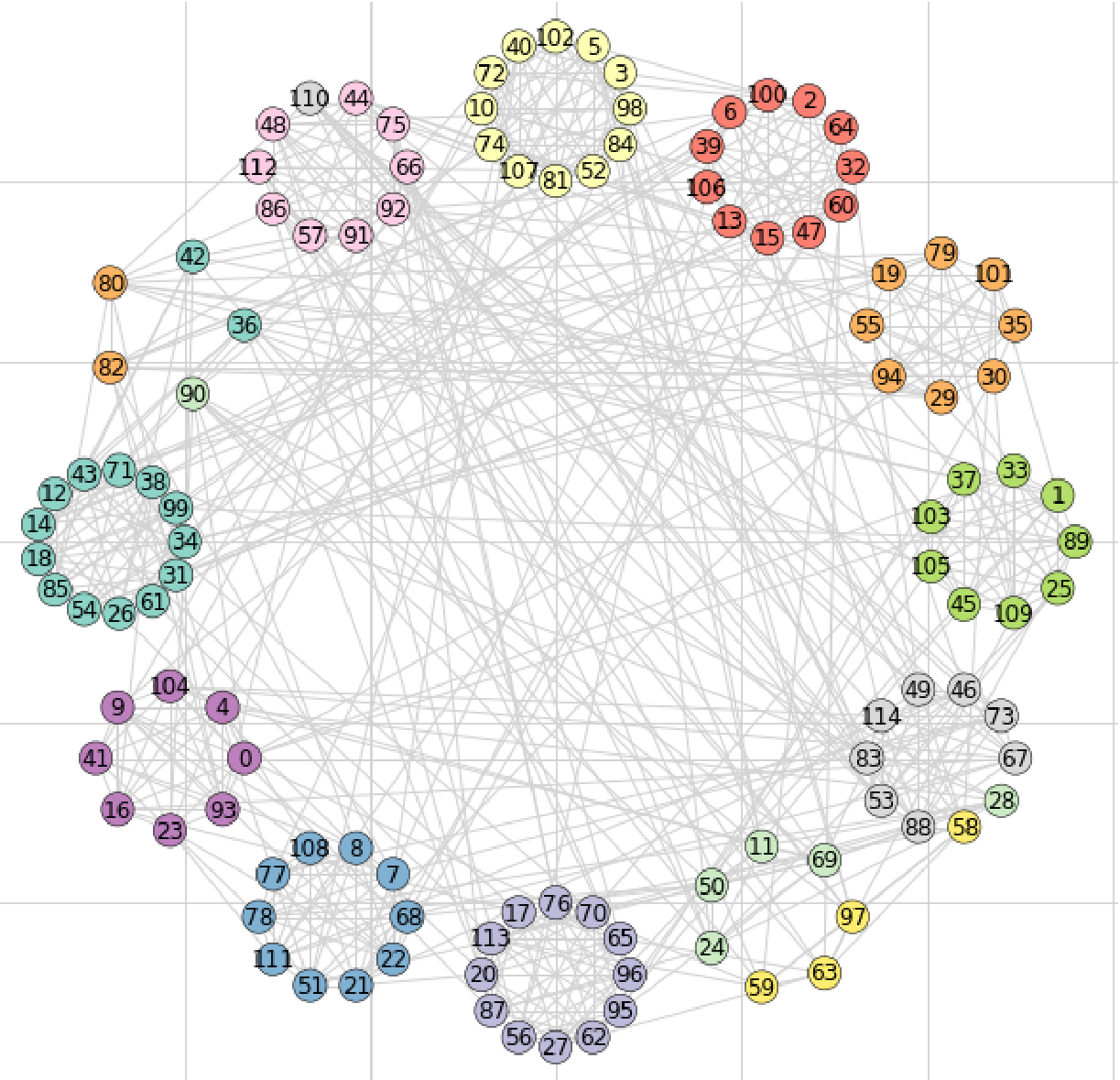}
    \caption{}
  \end{subfigure}

    \medskip

    \begin{subfigure}[t]{.39\textwidth}
    \centering
    \includegraphics[width=\linewidth]{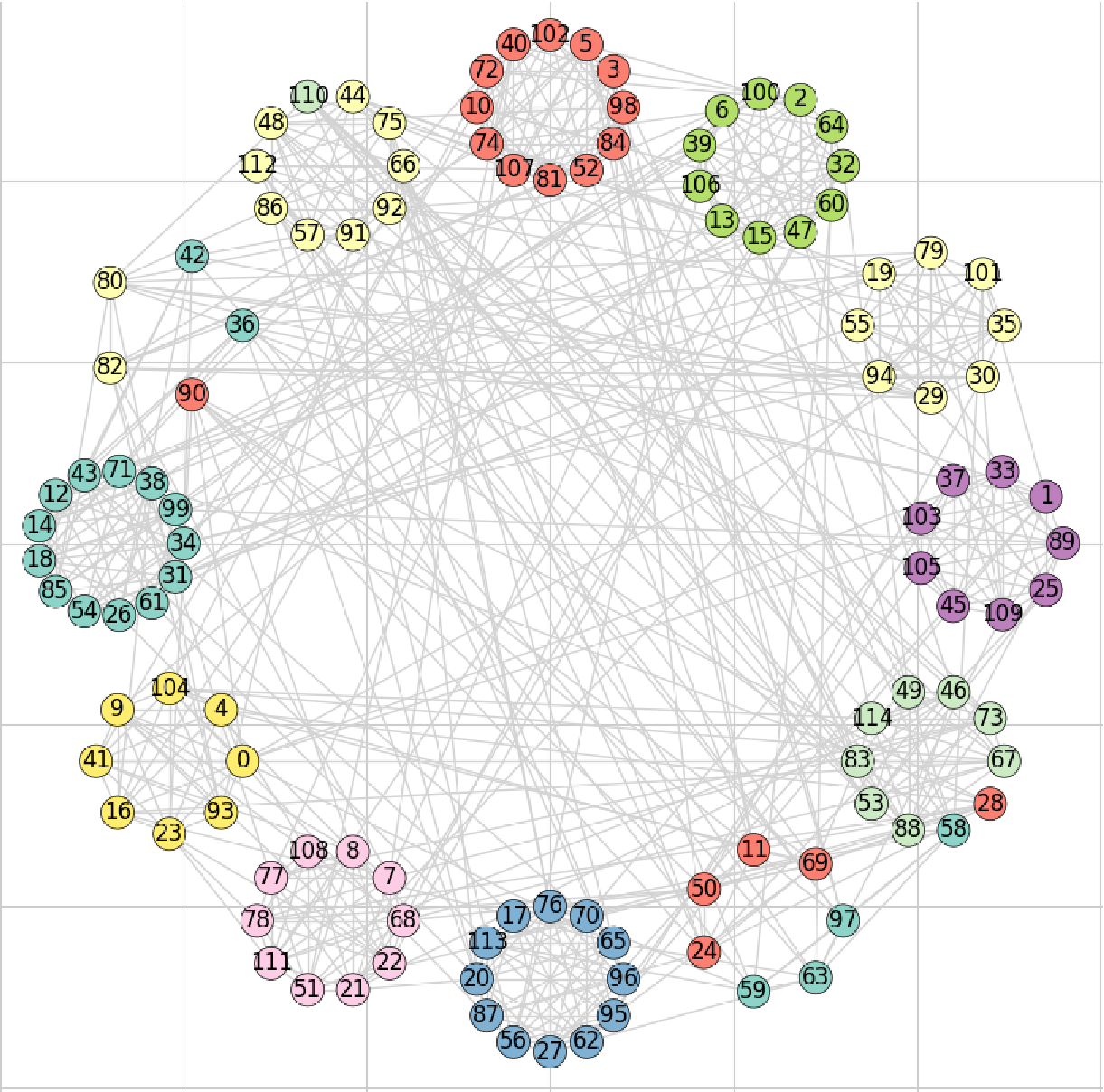}
    \caption{}
  \end{subfigure}
  \hfill
  \begin{subfigure}[t]{.39\textwidth}
    \centering
    \includegraphics[width=\linewidth]{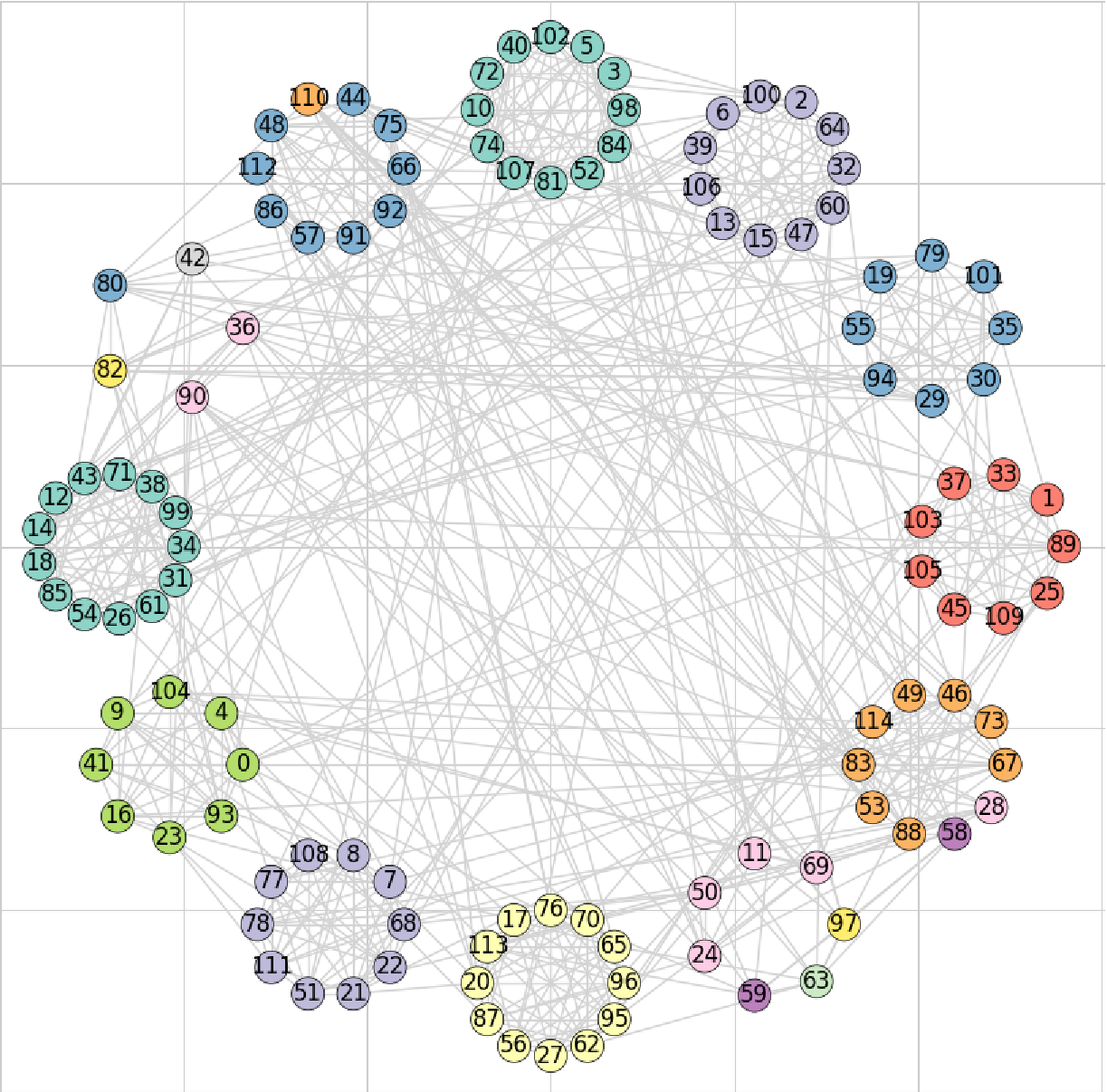}
    \caption{}
  \end{subfigure}

      \medskip

    \begin{subfigure}[t]{.39\textwidth}
    \centering
    \includegraphics[width=\linewidth]{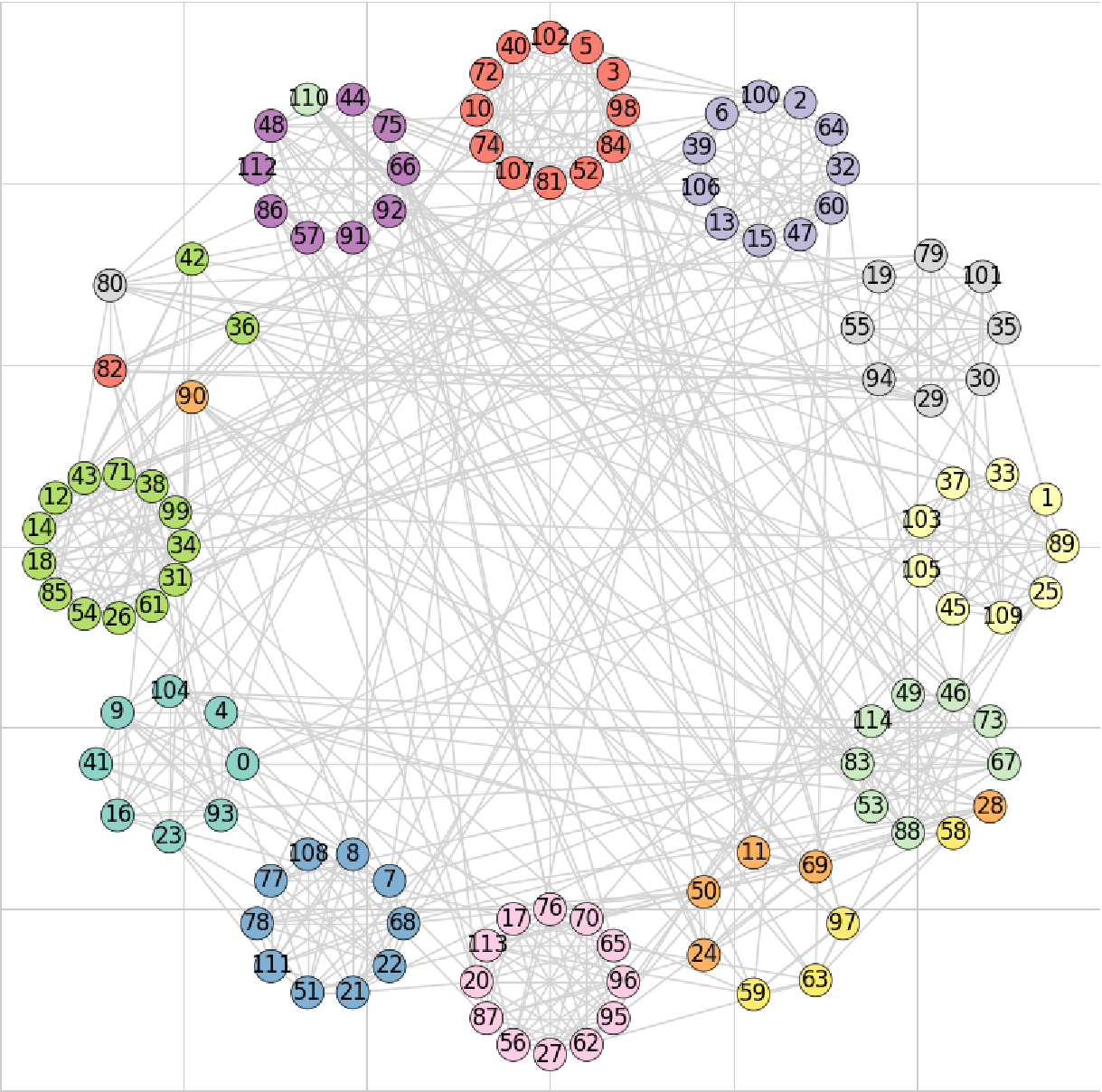}
    \caption{}
  \end{subfigure}
  \hfill
  \begin{subfigure}[t]{.39\textwidth}
    \centering
    \includegraphics[width=\linewidth]{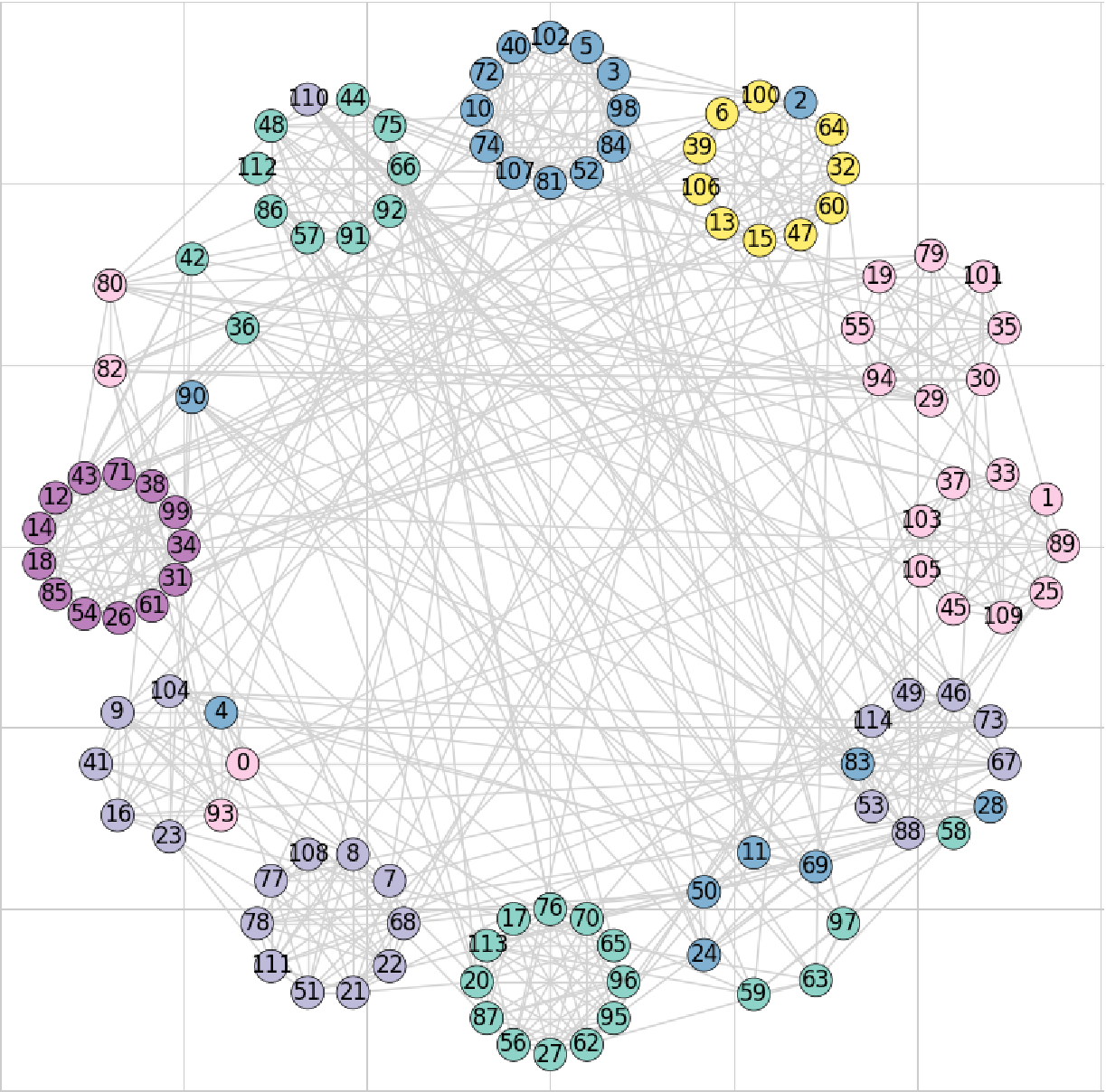}
    \caption{}
   \end{subfigure}
    \caption{Community detection in the American College Football network: ground truth communities \textbf{(a)}, and communities detected by sequential edge-deletion using the ORC \textbf{(b)}, $\AF{3}$ \textbf{(c)}, and $\AF{4}$ \textbf{(d)}. We also display the communities detected by the Girvan-Newman \textbf{(e)} and Louvain \textbf{(f)} algorithms.}\label{fig:football_accuracy}
\end{figure}

\begin{figure}[!t]
    \centering
    \subfloat[]{\includegraphics[width=.39\linewidth]{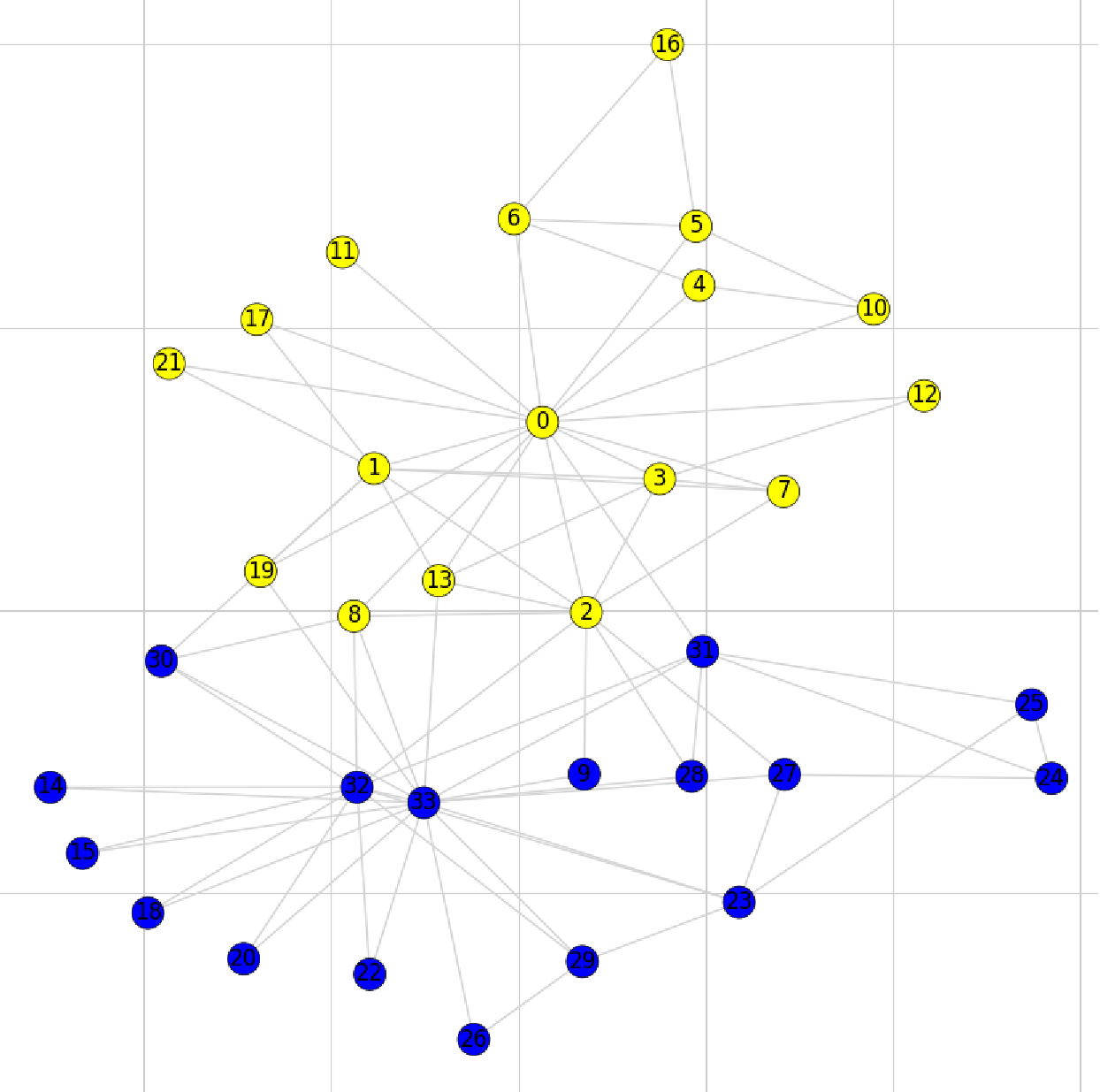}}
    \hfill
    \subfloat[]{\includegraphics[width=.39\linewidth]{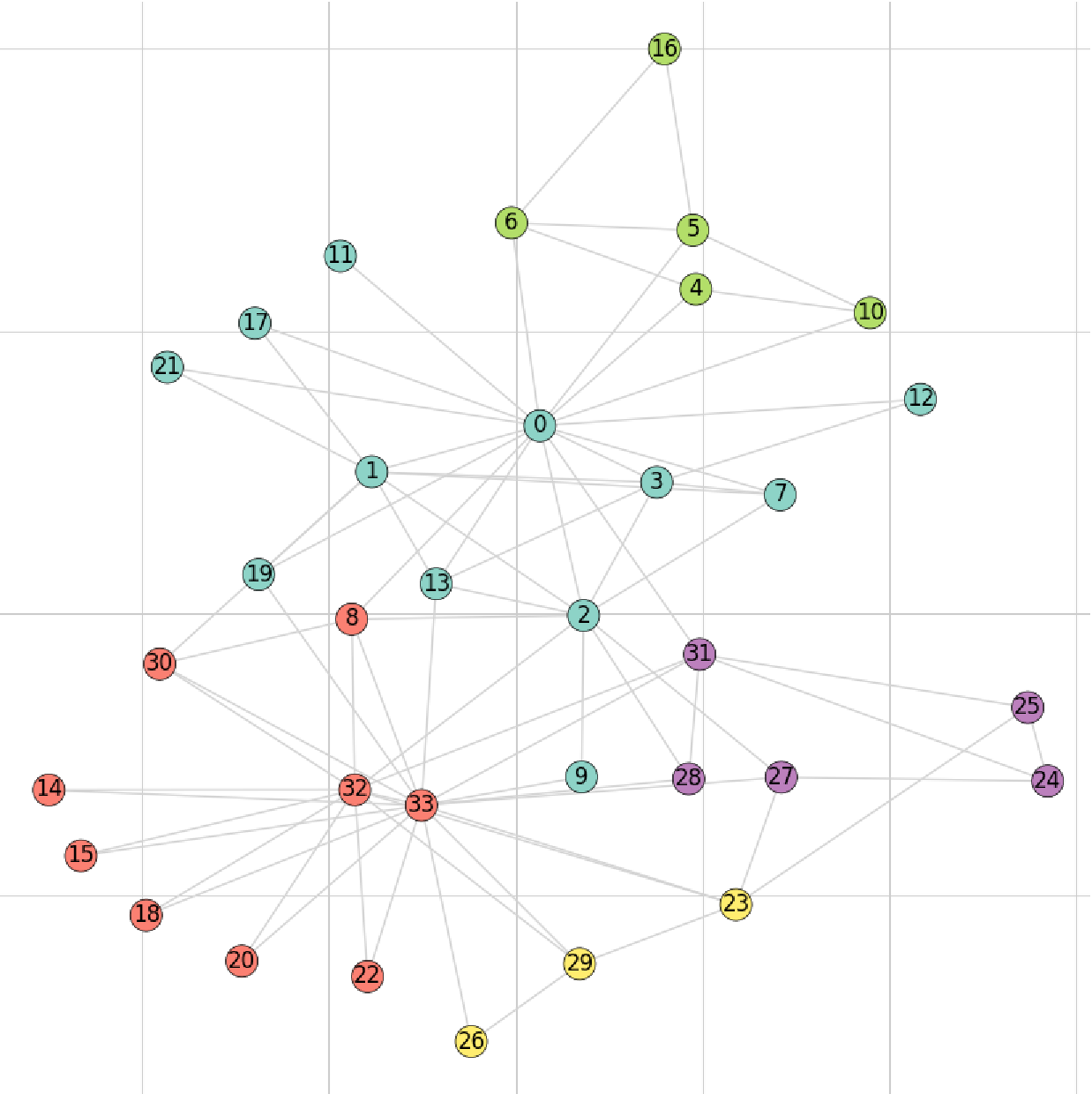}}\\
    \medskip
    \subfloat[]{\includegraphics[width=.39\linewidth]{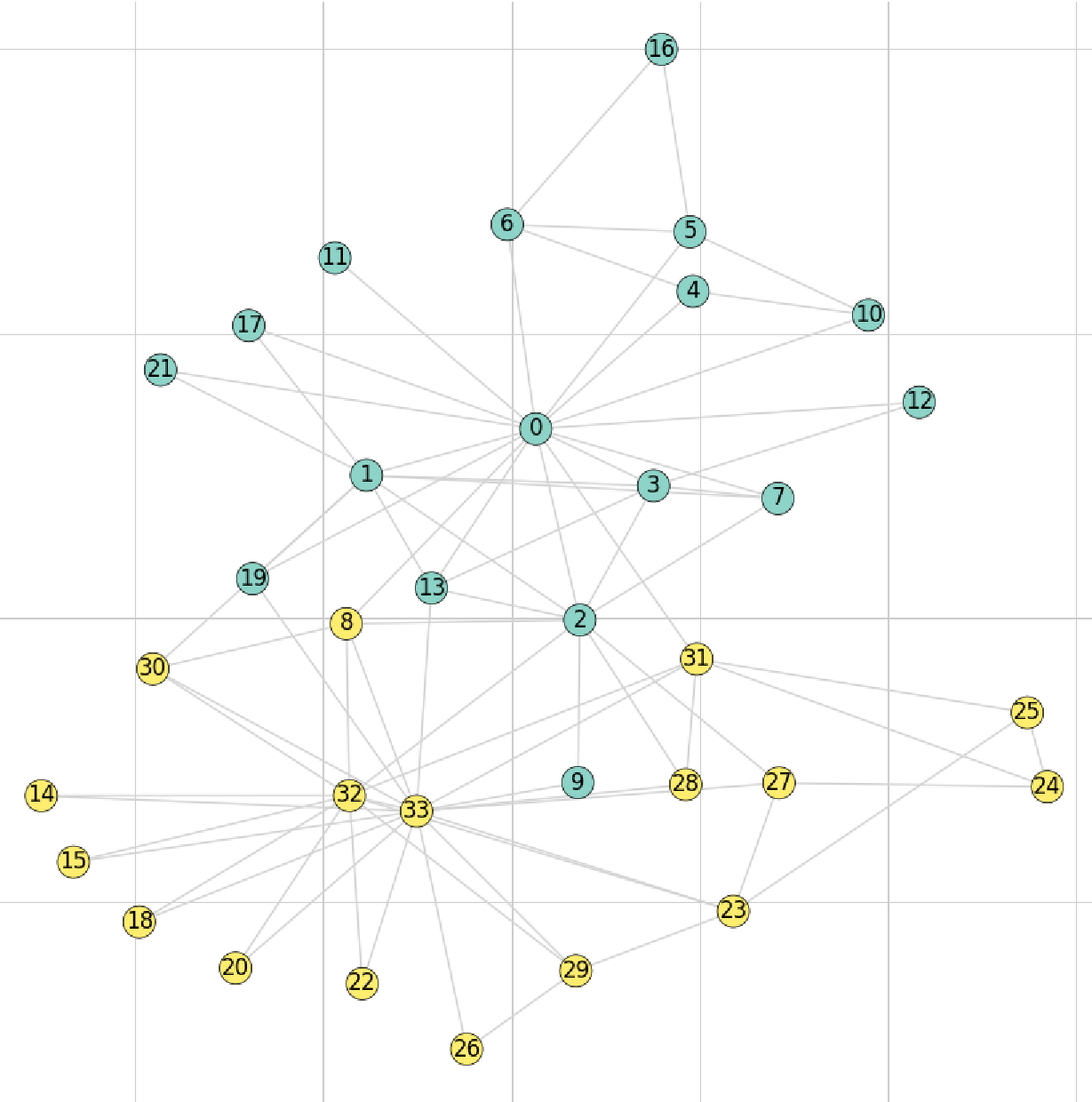}}
    \hfill
    \subfloat[]{\includegraphics[width=.39\linewidth]{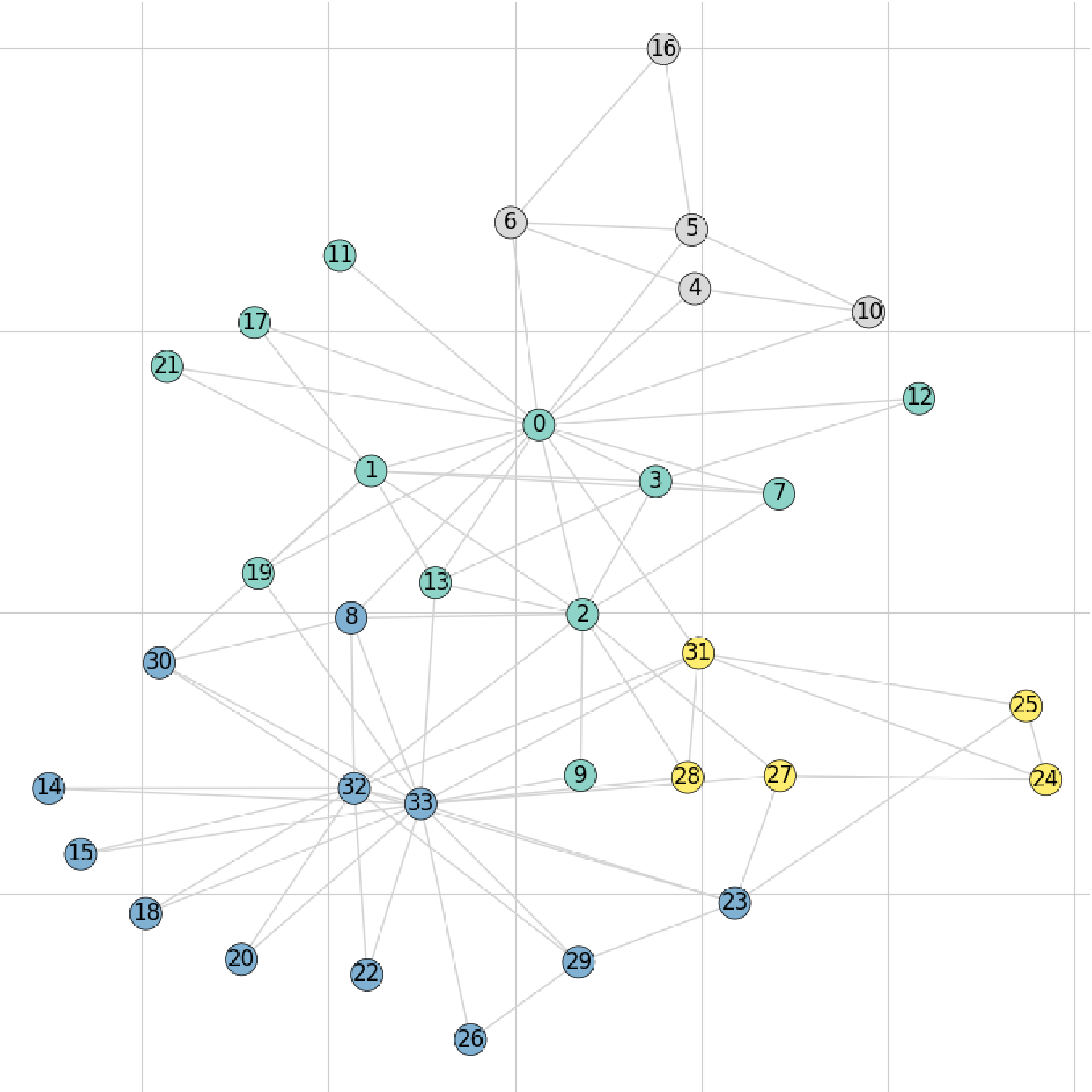}}\\
    \medskip
    \subfloat[]{\includegraphics[width=.39\linewidth]{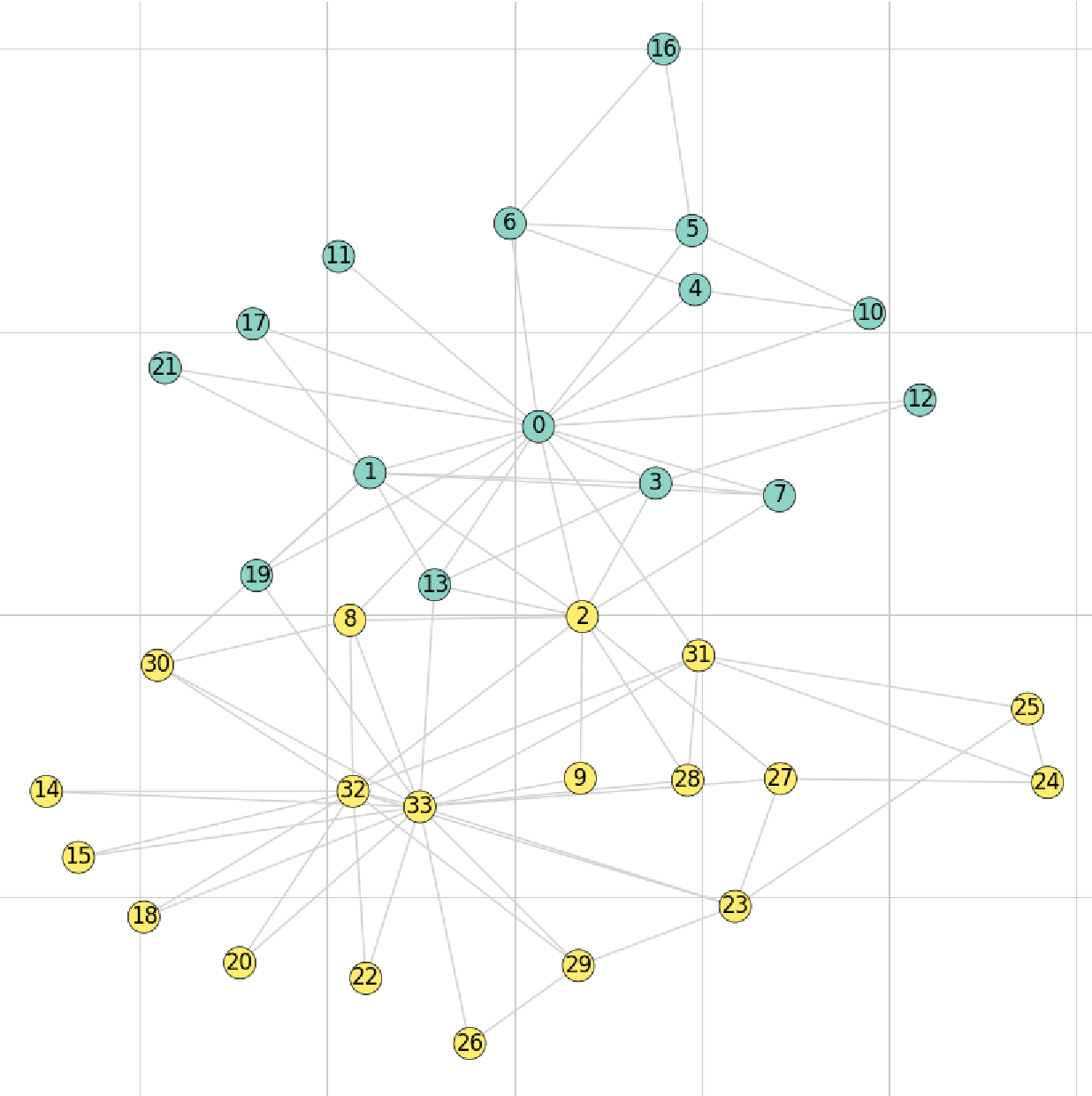}}
    \hfill
    \subfloat[]{\includegraphics[width=.39\linewidth]{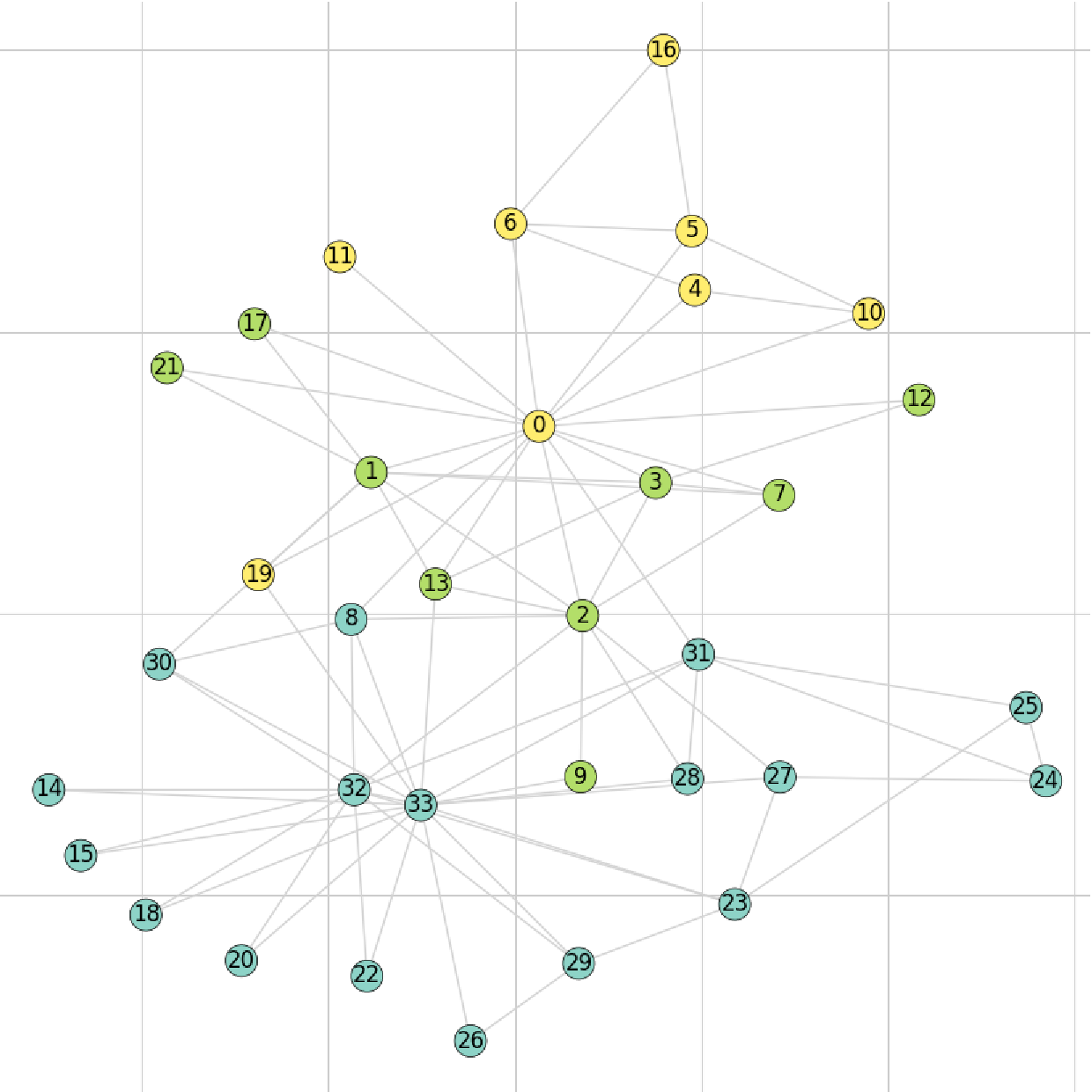}}
    
    \caption{Community detection in the Karate Club network: ground truth communities \textbf{(a)}, and communities detected by sequential edge-deletion using the ORC \textbf{(b)}, $\AF{3}$ \textbf{(c)} and $\AF{4}$ \textbf{(d)}. We also display the communities detected by the Girvan-Newman \textbf{(e)} and the Louvain \textbf{(f)} algorithms.}\label{fig:karate_accuracy}
\end{figure}


\clearpage
\section*{Acknowledgements}

The authors declare that they have no conflict of interest.

\noindent M.W. acknowledges partial support from NSF award 2112085.

\noindent R.L. acknowledges support from the EPSRC Grants EP/V013068/1 and EP/V03474X/1.

\bibliographystyle{unsrt.bst}
\bibliography{references}
\end{document}